\documentclass[letterpaper,10pt,oneside,onecolumn,reqno]{amsart}
\usepackage{amsmath, amssymb}
\usepackage{url}
\usepackage{enumerate}
\usepackage[all]{xy}
\usepackage[left=1in,top=1in,right=1in,bottom=1in]{geometry}
\usepackage{setspace}
\usepackage{hyperref}
\usepackage{wrapfig}
\usepackage{graphicx}
\usepackage{appendix}
\usepackage{color}
\usepackage[usenames,dvipsnames]{xcolor}

\onehalfspace
\newcommand{\A}{\mathcal A}

\newcommand{\B}{\mathcal B}
\newcommand{\C}{\mathcal C}
\newcommand{\DD}{\mathcal D}
\newcommand{\EE}{\mathbb E}

\newcommand{\F}{\mathcal F}
\newcommand{\G}{\mathcal G}

\newcommand{\II}{\mathcal I}
\newcommand{\JJ}{\mathcal J}

\newcommand{\M}{\mathcal M}

\renewcommand{\P}{\mathbb P}
\newcommand{\PP}{\mathbb P}
\renewcommand{\O}{\mathcal O}
\newcommand{\OO}{\mathrm O}
\newcommand{\Q}{\mathbb G}
\newcommand{\R}{\mathbb R}
\newcommand{\T}{\mathcal T}

\newcommand{\V}{\mathcal V}

\newcommand{\Z}{\mathbb Z}
\newcommand{\ZZ}{\mathcal Z}

\renewcommand{\hat}{\widehat}
\renewcommand{\tilde}{\widetilde}
\renewcommand{\bar}{\overline}

\newcommand{\density}{\operatorname{density}}
\newcommand{\jumps}{\operatorname{jumps}}
\newcommand{\ZZcore}{\ZZ_*}
\renewcommand{\div}{\operatorname{div}}
\newcommand{\supp}{\operatorname{supp}}
\newcommand{\SPD}{\operatorname{SPD}}
\newcommand{\Sym}{\operatorname{Sym}}

\newcommand{\deuc}{d_{\operatorname{Euc}}}

\newcommand{\Kmax}{K_{\mathrm{max}}}

\newcommand{\Rshape}{R_{\mathrm{shape}}}
\newcommand{\Rmax}{R_{\mathrm{max}}}
\newcommand{\Tmin}{T_*}

\newcommand{\E}{\mathrm{e}}				% Euler's constant.
\newcommand{\D}{\mathrm{d}}				% Differential operator sans space.
\newcommand{\sD}{\, \mathrm{d}}				% Differential operator with space.
\newcommand{\I}{\mathrm{i}}				% Square root of $-1$.
\newcommand{\oo}{\infty}
\newcommand{\Leb}{\operatorname{Leb}}
\newcommand{\diam}{\operatorname{diam}}

\newcommand{\SO}{\operatorname{SO}}

\theoremstyle{definition}

\numberwithin{equation}{section}
\newtheorem{env_thm}{Theorem}[section]
\newtheorem{env_def}[env_thm]{Definition}
\newtheorem{env_lem}[env_thm]{Lemma}
\newtheorem{env_sublem}[env_thm]{Sublemma}
\newtheorem{env_cor}[env_thm]{Corollary}
\newtheorem{env_pro}[env_thm]{Proposition}

\newtheorem{env_rem}{Remark}

\newtheorem*{env_mainthm}{Main Theorem}

\newtheorem*{env_ack}{Acknowledgements}

\title{Geodesics of Random Riemannian Metrics}

\author[T. LaGatta]{Tom LaGatta}
\address{Courant Institute of Mathematical Sciences \\ New York University \\ 251 Mercer St. \\ New York, New York 10012}
\email{tlagatta@gmail.com}

\author[J. Wehr]{Jan Wehr}
\address{Department of Mathematics\\ The University of Arizona \\ 617 N. Santa Rita Ave. \\ P.O. Box 210089 \\ Tucson, AZ 85721}
\email{wehr@math.arizona.edu}

\date{\today}

\begin{document}

\begin{abstract}
	We analyze the disordered Riemannian geometry resulting from random perturbations of the Euclidean metric. We focus on geodesics, the paths traced out by a particle traveling in this quenched random environment.  By taking the point of the view of the particle, we show that the law of its observed environment is absolutely continuous with respect to the law of the random metric, and we provide an explicit form for its Radon-Nikodym derivative. We use this result to prove a ``local Markov property'' along an unbounded geodesic, demonstrating that it eventually encounters any type of geometric phenomenon. We also develop in this paper some general results on conditional Gaussian measures. Our Main Theorem states that a geodesic chosen with random initial conditions (chosen independently of the metric) is almost surely not minimizing. To demonstrate this, we show that a minimizing geodesic is guaranteed to eventually pass over a certain ``bump surface,'' which locally has constant positive curvature.  By using Jacobi fields, we show that this is sufficient to destabilize the minimizing property.
\end{abstract}

\maketitle
\tableofcontents

%\textbf{Preprint version, please do not disseminate. April 17, 2012}

\part{Random Differential Geometry}

	\section{Introduction} \label{sect_introduction}
	
	In this article, we study the effects of disorder on the standard Euclidean metric in $\R^d$. In particular, we investigate the geodesics of a random Riemannian metric. Our Main Theorem is that, under relatively general conditions, a geodesic with random initial conditions (chosen independently of the metric) is not minimizing. To prove this result rigorously, we draw a number of ideas from different fields, including probability theory, disordered systems, and differential geometry. Instrumental in our argument is our theorem on the point of view of a particle (Theorem \ref{POVthm}), which states that as a particle travels along a geodesic, the law of its observed environment is absolutely continuous with respect to the original law of the random metric. We hope that the tools developed in this article will be useful for future works in random differential geometry.
		 
	The basic assumptions on the random metric are that its law $\PP$ is stationary and ergodic with respect to the action of the Euclidean group,\footnote{i.e., the measure $\PP$ is invariant under translations and rotations.} and also that the local values of the metric are independent when separated by sufficient Euclidean distance. In a previous article \cite{lagatta2009shape}, we proved a Shape Theorem for such models: with probability one, random Riemannian balls grow asymptotically like Euclidean balls. While we do not yet have rigorous control on the fluctuations of balls, we conjecture that our model lies in the Kardar-Parisi-Zhang (KPZ) universality class of growth models. Thus, the fluctuations of the random balls away from their asymptotic limiting shape should be of size $t^\chi$ for some universal exponent $\chi$ depending on the dimension $d$ of the model, and not on the detailed character of the fluctuations. Fluctuations from the limiting shape are related to fluctuations of minimizing geodesics, and in a related model of lattice first-passage percolation, this scaling relation was recently proved by Chatterjee \cite{chatterjee2011universal} and Auffinger-Damron \cite{auffinger2011simplified}. For our model with $d=2$ (random metrics in the plane), we conjecture that $\chi = 1/3$, consistent with other models in the KPZ universality class. Indeed, preliminary numerical results of the first author with physicists Javier Rodr\'{i}guez Laguna, Silvia N. Santalla and Rodolfo Cuerno Rejado suggest that $\chi = 1/3$.%, and that the fluctuations satisfy Tracy-Widom statistics. 

	\subsection{Random Riemannian Metrics} \label{sect_rrm}
	
	Consider the measurable space $\Omega_+ = C^2(\R^d, \SPD)$ of $C^2$-smooth symmetric $2$-tensor fields on $d$-dimensional Euclidean space, equipped with its Borel $\sigma$-algebra $\B(\Omega_+)$.\footnote{$\SPD$ denotes the finite-dimensional space of $d \times d$ symmetric, positive-definite matrices. The space $\SPD$ is the positive cone in the finite-dimensional Banach space $\Sym$ of $d \times d$ symmetric matrices.} The space $\Omega_+$ parametrizes the Riemannian metrics in Euclidean space, and a ``random Riemannian metric'' is described by a probability law $\PP$ on the pair $\big(\Omega_+, \B(\Omega_+)\big)$.\footnote{A note on smoothness assumptions. We use the $C^1$-smoothness of the metric to define geodesics, by way of the Levi-Civita connection (represented by Christoffel symbols). The only place we use $C^2$-smoothness is in Section \ref{sect_bump_short}, where we construct a ``bump metric'' as a function of the curvature tensor of the metric.}
	
	We consider a specific class of random Riemannian metrics generated by Gaussian fields, but most of our arguments hold in wider generality. Let $c : \R \to \R$ be a symmetric, Gaussian covariance function which is non-degenerate ($c(0) > 0$), compactly supported (if $r \ge 1$, then $c(r) = 0$), and $5$-times differentiable. That such covariance functions exist is non-trivial, and we provide an example in Example \ref{exa_gneitingcov}. 
	
	In Appendix \ref{app_gaussian}, we construct a mean-zero, Gaussian random $2$-tensor field with a covariance $4$-tensor $c_{ijkl}$ generated from such a covariance function $c$.\footnote{The Appendices to this article are posted as supplementary material on the web \cite{lagatta2012geodesicsII}.} The assumptions assure that the Gaussian field $\xi$ is almost-surely $C^2$-differentiable, and its law $\Q$ is a Gaussian measure on the Fr\'echet space $\Omega := C^2(\R^2, \Sym)$ of symmetric $2$-tensor fields. We let $\G$ denote the completion of the Borel $\sigma$-algebra $\B(\Omega)$, by including the null sets for the Gaussian measure $\Q$.
	
	Next, let $\varphi : \R \to (0,\oo)$ be a smooth, increasing function satisfying certain growth conditions, which we use to locally transform a symmetric tensor to a symmetric positive-definite tensor (metric).\footnote{The precise growth conditions on $\varphi$ are that there exist constants $C$ and $\eta_1 \le \eta_2$ so that $\tfrac 1 C u^{\eta_1} \le |\varphi(u)|_{C^{2,1}} \le C u^{\eta_2}$ as $u \to \oo$ and $\tfrac{1}{C| u|^{\eta_2}} \le |\varphi(u)|_{C^{1,1}} \le \tfrac{C}{|u|^{\eta_1}}$ as $u \to -\oo$. The notation $|\cdot|_{C^{\alpha,1}}$ denotes the maximum of the function and its first $\alpha$ derivatives at $u$, along with the local Lipschitz constant of the $\alpha$th derivative. Such a function $\varphi$ can be easily defined piecewise.} We define the continuous non-linear operator $\Phi : \Omega \to \Omega_+$ pointwise: $\Phi(\xi)(u) := \varphi( \xi(u) )$ for all $u \in \R^d$. We define the probability measure $\PP := \Q \circ \Phi^{-1}$ as the push-forward of the Gaussian measure onto the space of metrics, and henceforth we let $g$ represent a random Riemannian metric with law $\PP$.\footnote{That is, $g$ is an $\Omega_+$-valued random variable, defined using the measurable function $\Phi$ on the complete probability space $(\Omega, \F, \Q)$.} We let $\F$ denote the completion of the Borel $\sigma$-algebra $\B(\Omega_+)$.
	
	The fundamental property of our model is that the law $\PP$ is invariant under the (orientation-preserving) isometries of Euclidean space: translations and rotations. Owing to the construction by a Gaussian random field with compactly supported covariance, our random metrics also satisfy a strong independence property: if functions of the metric (e.g., the metric tensor, Christoffel symbols or the Riemann curvature tensor) are separated by Euclidean distance at least $1$, then they are independent. While most of the arguments in our article can be extended to more general families of random metrics, we note that Theorem \ref{P_lem} in particular relies on structural properties of Gaussian measures, and Theorem \ref{thm_markov} relies on the independence property.
	
	For every compact $D \subseteq \R^d$, let $\F_D$ denote the $\sigma$-algebra generated by the random metric in an infinitesimal neighborhood of $D$.\footnote{\label{FD} Formally, let $D^\epsilon = \{ x \in \R^d : \deuc(x,D) \le \epsilon \}$ denote the $\epsilon$-neighborhood of $D$, and define the $\sigma$-algebra $\F_D$ to be the completion of the $\sigma$-algebra $\cap_{\epsilon > 0} \sigma\big( g_{ij}(x) : x \in D^\epsilon \big)$, by including the null subsets with respect to the measure $\PP$.} The compact support of the Gaussian covariance function implies that if $\deuc(D, D') \ge 1$, then the $\sigma$-algebras $\F_D$ and $\F_{D'}$ are independent with respect to $\PP$. If $x \in \R^d$ is a single point, we write $\F_x := \F_{ \{x\} }$.  By construction, the random variables $g_{ij}(x)$ are $\F_x$-measurable; since $\F_x$ depends on the metric in an infinitesimal neighborhood of $x$, the pointwise derivatives $g_{ij,k}(x)$ and $g_{ij,kl}(x)$ are also $\F_x$-measurable, as is the curvature tensor.  
	
		Every Riemannian metric $g \in \Omega_+$ generates a norm on the tangent bundle $T \R^d$, given by $\| v \|_g := \sqrt{\langle v, g(x) v \rangle}$ when $v \in T_x \R^d$, as well as a distance function $d_g(x,y) := \inf_\gamma \int \|\dot\gamma\|_g \, \D t$, where the infimum is taken over all smooth curves $\gamma$ connecting $x$ to $y$. The extrema of the length functional are called geodesics, and they solve the geodesic equation $\ddot \gamma^k = -\Gamma_{ij}^k(g, \gamma) \dot \gamma^i \dot \gamma^j$, where the terms $\Gamma_{ij}^k(g,x)$ denote the Christoffel symbols for the metric $g$ at the point $x$.\footnote{We follow the Einstein convention on the right side, summing over repeated indices.}  Without loss of generality, we assume that geodesics are parametrized by Riemannian arc length, i.e., that $\|\dot \gamma(t) \|_g \equiv 1$ for all $t \in \R$. The Shape Theorem (Theorem \ref{shapecor}) implies that with probability one, $g$ is a complete Riemannian metric, so all geodesics are defined for all time. We say that a geodesic $\gamma$ is (forward) minimizing when $d_g( \gamma(t), \gamma(t')) = |t-t'|$ for all times $t, t' \ge 0$. Throughout this article we use concepts and theorems from differential geometry; we detail some notation in Appendix \ref{geombg}, and for a more substantial overview the reader may consult a standard text such as Lee \cite{lee1997rmi}.
	
	We use the notation $\gamma_{x,v} = \gamma_{x,v}(g,\cdot)$ to denote the unit-speed geodesic with initial conditions $\gamma(0) = x$ and $\dot \gamma(0) = v/\|v\|_g$. Since the metric is random, for each initial condition $(x,v) \in T\R^d$, $\gamma_{x,v}$ represents a curve-valued random variable.\footnote{A technical point is that the random variable $\gamma_{x,v} : \Omega_+ \to C^2(\R, \R^d)$ is only a partial function, with domain the set of complete Riemannian metrics. The reason is that if the metric is not complete, the geodesic blows up in finite time. This is not an issue for our model, since our Shape Theorem \cite{lagatta2009shape} implies that a wide class of random metrics is geodesically complete with probability one.}

	We are concerned in this article with almost-sure properties of geodesics with random initial conditions, selected independently from the metric.\footnote{We model the distribution of initial conditions using a probability measure $\beta$ which is absolutely continuous with respect to the Euclidean kinematic measure on the unit tangent bundle $UT\R^d \cong \R^d \times S^{d-1}$. The independence condition means that we consider the product measure $\P' := \beta \times \P$ on the larger probability space $UT\R^d \times \Omega_+$.} By translation- and rotation-invariance of the distribution $\PP$, it suffices to consider the geodesic $\gamma := \gamma_{0,\E_1}$ with deterministic starting conditions $(0,\E_1)$, since almost-sure properties of the random curve $\gamma$ will be the same as for a geodesic with random starting conditions. The random geodesic $\gamma$ represents a typical trajectory for a particle traveling in a random Riemannian environment.

	\begin{env_mainthm}
		Suppose that $d=2$.  Then $\PP(\mbox{$\gamma$ is minimizing}) = 0$.
%		That is, with probability one, $\gamma$ is not minimizing.
	\end{env_mainthm}
	
	The restrictive assumption $d=2$ appears to be a technical artifact of our analysis; we believe that the Main Theorem is true in the general case. The precise location of the obstruction is our Theorem \ref{POVthm}, which describes the law of the geometric environment as seen from the point of view of a particle traveling along a geodesic. Our proof of that theorem extends the Geman-Horowitz argument \cite{geman1975random}, in which we exploit the abelian nature of the two-dimensional rotation group. A version of that result should be true for arbitrary $d$; if such a theorem is proved, then all of the subsequent analysis of this paper will extend easily, \emph{mutatis mutandis}. 

	\subsection{Outline of Proof}

	The proof of the Main Theorem breaks into two cases. On the event $\{\mbox{$\gamma$ is bounded}\}$, the proof follows immediately from the fact that all minimizing geodesics are transient, hence unbounded (Theorem \ref{transientgeodesics}). Owing to the soft nature of the argument, we have no quantitative estimate for when $\gamma$ loses the minimization property in the bounded case.
	
	Curvature plays an essential role in the proof of the unbounded case. If the scalar curvature of the random metric were globally non-positive, then the Cartan-Hadamard theorem \cite{lee1997rmi} would imply that all geodesics are minimizing. Therefore, the presence of positive curvature is a necessary condition for destabilizing the minimization property. We exploit this in our proof of the unbounded case, and construct a model ``bump surface'' which has enough positive curvature to draw geodesics together. In particular, such geodesics must develop conjugate points, which by Jacobi's theorem (Theorem 10.15 of \cite{lee1997rmi}) are an obstruction to minimization.
	
	Conditioned on the event that $\gamma$ is unbounded, our Inevitability Theorem (Theorem \ref{lem_musthappen}) states that, under a certain uniformity condition \eqref{lem_musthappen_estimate}, an unbounded geodesic will eventually encounter any reasonable type of geometric environment. In Section \ref{sect_frontiertimes}, we show that, conditioned on the event $\{\mbox{$\gamma$ is minimizing}\}$, this condition \eqref{lem_musthappen_estimate} is satisfied. Therefore, the geodesic is guaranteed to eventually encounter a bump surface, and thereby develop conjugate points. This contradicts the assumption that $\gamma$ is minimizing. % The proof of the Main Theorem is given in Section \ref{sect_proofofmaintheorem}.
	
	In the unbounded case, we do have an estimate on the time for which $\gamma$ is minimizing. Let $T_* = \sup\{ t > 0 : \mbox{$\gamma$ is minimizing between times $0$ and $t$} \}$ be the maximum such time. Theorem \ref{notminimizing_Wthm} demonstrates that, conditioned on the event that $\gamma$ is unbounded, the random variable $T_*$ has exponential tail decay.
	
	\subsection{Organization}
	
	The article is divided into three main parts. In Part I, we outline some of the general geometric features of random differential geometry. The arguments in this section are robust, and should be easily adaptable to more general settings. In Part II, we present a number of auxiliary theorems necessary to prove the Main Theorem. In Part III, we present the proofs of these auxiliary theorems.
	
	The Appendices to this article are posted as supplementary material on the web \cite{lagatta2012geodesicsII}. In Appendix \ref{app_gaussian}, we present some general results on the construction of Gaussian random fields. In Appendix \ref{app_shape}, we restate our Shape Theorem from \cite{lagatta2009shape}, specialized to the setting of this article. In Appendix \ref{app_geomgeod}, we state some straightforward consequences on the geometry of geodesics for a random metric. In Appendix \ref{geombg}, we provide a rapid introduction to Riemannian geometry for the unfamiliar reader. In Appendix \ref{analytictools}, we present some analytic estimates which we use in the article. In Appendix \ref{proof_mo_lem}, we present the construction of the conditional mean operator for Gaussian measures. In Appendix \ref{fermiproof}, we describe Fermi normal coordinates, which we use in our construction of the bump metric.
	
%	For the remainder of the introduction, we outline the main results of our article and discuss applications.  In Section 2, we define our model formally.  In Section 3, we discuss the method of the point-of-view (POV) of the particle.  In Section 4, we study the exit-time POV process, which describes the environment viewed along a geodesic as it exits large Euclidean balls.  In Section 5, we describe the Local Markov Property, which describes the law of the environment near these exit locations.  In Section 6, we show that the law of our random metric satisfies the continuous disintegration property of \cite{lagatta2010continuous}, and we prove the Inevitability Theorem, which states that along an unbounded geodesic, anything must happen.  In Section 7, we give our proofs of the theorems in this paper.  In Appendix A, we give a brief overview of Riemannian geometry.  In Appendix B, we introduce some analytic tools we use throughout this article.  In Appendix C, we prove Lemma \ref{mo_lem}.  In Part II \cite{lagatta2012geodesicsII}, we use the tools developed in this article to show that a randomly-selected geodesic is almost surely not minimizing.
	
	\begin{env_ack}
		The authors would like to thank Antonio Auffinger, Benjamin Bakker, Robert Bryant, Sourav Chatterjee, Michael Damron, Partha Dey, Krzysztof Gaw\c edzki,  David Glickenstein,  Geoffrey Grimmett, Joey Hirsh, Christopher Hoffman, Alfredo Hubard, Thomas Kennedy, Sun Hyoung Sonya Kim, Michael Marcus, Charles M. Newman, Benjamin Pittman-Polletta, Javier Rodr\'{i}guez Laguna, Jay Rosen, David Sanders, Silvia N. Santalla, Thomas Spencer, John Terilla, Srinivasa Varadhan, Joseph Watkins, Brad Weir, Willie Wong, Maciej Wojtkowski, Deane Yang, and Lai-Sang Young for many helpful discussions on aspects of this project.

%	For generously answering T.L.'s many geometry questions on MathOverflow.net, he thanks Robert Bryant, Willie Wong and Deane Yang.
		
	T.L. was supported by NSF VIGRE Grant No. DMS-06-02173 at the University of Arizona, and by NSF PIRE Grant No. OISE-07-30136 at the Courant Institute.  J. W. was partially supported by the NSF grant DMS 1009508.
	
	\end{env_ack}

	\label{page_introduction_end}

	\section{Minimizing Geodesics} \label{sect_mingeos}
	
	While a geodesic with random initial conditions is a.s. not minimizing (by the Main Theorem), there are many minimizing geodesics starting at any point. Let $\gamma_v := \gamma_{0,v}$ denote the unit-speed geodesic starting at the origin in direction $v \in S^{d-1}$. For any Riemannian metric $g \in \Omega_+$, let $\V_g = \{ v \in S^{d-1} : ~\mbox{$\gamma_v$ is minimizing}~ \}$ denote the set of initial directions which yield (forward) minimizing geodesics. The next result demonstrates that this set is always compact and, for complete metrics, non-empty; we owe the argument to M. Wojtkowski. 
	
	\begin{env_pro} \label{closedV}
		For all complete metrics $g \in \Omega_+$, the set $\V_g$ is compact and non-empty.
	\end{env_pro}
	\begin{proof}
	 	Suppose that $v_n \in \V_g$, and $v_n \to v$ in $S^{d-1}$; we claim that $\gamma_v$ is minimizing. Fix times $t, t' \ge 0$, and let $x = \gamma_v(t)$ and $x' = \gamma_v(t')$ be the corresponding points on the curve. Since the geodesic flow is continuous with respect to the initial velocity, $\gamma_{v_n}(t) \to x$ and $\gamma_{v_n}(t') \to x'$ as $n\to \oo$. The distance function $d_g$ is continuous and the (finite) geodesic segments $\gamma_{v_n}$ are minimizing, so $d_g(x,x') = \lim_{n\to\oo} d_g(\gamma_{n}(t), \gamma_{n}(t')) = |t-t'|,$ which proves that $\gamma_v$ globally minimizes length, hence $v \in \V_g$. Thus $\V_g \subseteq S^{d-1}$ is closed hence compact.
	 	
	 	The argument that $\V_g$ is non-empty is similar. Let $\gamma_n$ denote the minimizing geodesic segment from $0$ to $n\E_1$, which is well-defined since $g$ is complete.  Let $v_n := \dot \gamma_n(0)$ denote the initial direction of $\gamma_n$, so that $\gamma_n = \gamma_{v_n}$.  Since the unit sphere is compact, a subsequence $v_{n_j}$ converges to some direction $v \in S^{d-1}$. For any times $t, t' \ge 0$, consider $x = \gamma_v(t)$ and $x' = \gamma_v(t')$.  As in the previous argument, $d_g(x,x') = \lim_{j\to\oo} d_g(\gamma_{n_j}(t), \gamma_{n_j}(t')) = |t-t'|$, which proves that $\gamma_v$ is minimizing, hence $v \in \V_g$.
	\end{proof}

	When $g$ is the Euclidean metric $\delta \in \Omega_+$, we have $\V_\delta = S^{d-1}$ since all geodesics are minimizing rays. When $g$ denotes a random metric, the structure of $\V_g$ is more interesting; in this case, $g \mapsto \V_g$ denotes a set-valued random variable. A random metric represents a ``perturbation'' of Euclidean geometry, and this perturbation results in a highly non-trivial structure of the (random) set $\V_g$. When $d=2$, the Main Theorem implies that $\PP( v \in \V_g) = 0$ for all starting directions $v \in S^{1}$. The next corollary strengthens this result.

	\begin{env_cor}[Corollary to Main Theorem] \label{Vmeas0}
	 	Suppose that $d=2$.  With probability one, the set $\V_g$ has measure zero on the circle $S^{1}$.  That is, if $\nu$ denotes the uniform measure on $S^{1}$, then
	 		$$\PP\big( \nu(\V_g) = 0 \big) = 1.$$
	\end{env_cor}
	
	\begin{proof}
		For each $v \in S^{1}$, let $M_v = \{g : v \in \V_g\}$ be the event that the geodesic $\gamma_v$ is minimizing.  Since $d=2$, the Main Theorem and rotational invariance imply that $\PP(M_v) = \PP(v \in \V_g) = 0$.  Tonelli's theorem \cite{folland1999real} implies that
			\begin{eqnarray*}
				\EE \nu(\V_g) = \int_{\Omega_+} \nu(\V_g) \sD \PP(g)
				&=& \int_{\Omega_+} \nu(v : M_v ~\mbox{occurs} ) \sD \PP(g) = \int_{\Omega_+} \int_{S^{1}} 1_{M_v}(g) \sD \nu(v) \sD \PP(g) \\
				&=& \int_{S^{1}} \int_{\Omega_+} 1_{M_v}(g) \sD \PP(g) \sD \nu(v) = \int_{S^{1}} \PP(M_v) \sD \nu(v) = \int_{S^{1}} 0 \sD \nu(v) = 0,
			\end{eqnarray*}
		since $\PP(M_v) = 0$.  Since $\nu(\V_g)$ is a real-valued, non-negative random variable with mean zero, it vanishes almost surely.
	\end{proof}
	
	By a similar argument, with probability one, the set of starting conditions $(x,v)$ which yield minimizing geodesics has kinematic measure zero on the unit tangent bundle $UT\R^2 \cong \R^2 \times S^1$. These measure-zero statements are not just technical artifacts of our method:  heuristic arguments suggest that $\V_g$ is a random fractal. We conjecture that, with probability one, $\V_g$ is homeomorphic to the Cantor set.

	\subsection{Transience of Minimizing Geodesics}
	
	It is easy to see that every minimizing geodesic is transient: if $\gamma$ meets a compact set infinitely often, then it must have an accumulation point $x = \lim \gamma(t_k)$.  If $\gamma$ is minimizing and parametrized by Riemannian arc length, this means that the distance from $\gamma(t_k)$ to $x$ is infinite, which is a contradiction. %For the proof of the Main Theorem, we will need uniform estimates on the time it takes for minimizing geodesics to escape compact sets. 
	
	The next theorem is a much stronger version of this result in the context of random Riemannian metrics.  Let $K_g$ be a variable compact set, depending measurably on the metric $g$.\footnote{The space of compact sets of $\R^d$ is a metric space (equipped with the Hausdorff metric), and is thus a measurable space equipped with its Borel $\sigma$-algebra. For more on the Hausdorf metric, see Appendix \ref{analytictools}} The theorem states that with probability one, for all initial directions $v \in \V_g$, the geodesic $\gamma_v$ is guaranteed to exit the set $K_g$ by some time $T_g$, not depending on the initial direction. Our proof makes use of the Shape Theorem in order to get a quantitative estimate on the last exit time of the geodesic, but it is easy to prove such a theorem for general Riemannian metrics (see Remark \ref{rem_trans}).

	\begin{env_thm}[Minimizing Geodesics Are Uniformly Transient] \label{transientgeodesics}
		Let $g \mapsto K_g$ be a compact-set-valued random variable. With probability one, there exists a time $T_g$ so that that for all $v \in \V_g$ and $t > T_g$, $\gamma_v(t) \notin K_g$.  
	\end{env_thm}
	\begin{proof}
		Fix $\epsilon > 0$.  The Shape Theorem (Theorem \ref{shapecor}) implies that with probability one, there exists $\Rshape = \Rshape(g)$ such that if $r \ge \Rshape$, then $B(r) \subseteq B_g\big( (1+\epsilon)\mu r \big)$, where $B$ and $B_g$ denote the Euclidean and random Riemannian balls centered at the origin, respectively, and $\mu$ denotes the shape constant.  %$\tfrac{1}{1+\epsilon} t B(0,r) \subseteq B_g(r/\mu)$.

		Next, let $\hat K_g$ be the smallest Euclidean ball centered at the origin which contains $K_g$, and let $R_{K_g}$ denote its radius (so that $\hat K_g = B(R_{K_g})$). Set $R_g = \max\{ R_{K_g}, \Rshape(g) \}$, and define $T_g = (1+\epsilon) \mu R_g$, so that $K_g \subseteq \hat K_g \subseteq B(R_g) \subseteq B_g(T_g).$ 
		
		Now, suppose that $v \in \V_g$ and that $t > T_g$.  Since $\gamma_v$ is minimizing, $d_g(0, \gamma_v(t)) = t > T_g$.  This means that $\gamma_v(t) \notin B_g(T)$, hence $\gamma_v(t) \notin K_g$.  The time $T_g$ is an upper bound for the last exit time of $\gamma_v$ from the set $K_g$.
	\end{proof}

	Let $\tau_r(g)$ denote the last exit time of the geodesic $\gamma = \gamma_{\E_1}$ from the Euclidean ball of radius $r$. Using the Shape Theorem and a similar argument as in the previous proof, we obtain an upper and lower estimate on the last exit time. For almost every $g$ on the event $\{\mbox{$\gamma$ is minimizing}\}$, there exists $R_g$ so that if $r \ge R_g$, then
		\begin{equation} \label{taur_upperestimate}
			(1-\epsilon) \mu R_g \le \tau_r(g) \le (1+\epsilon) \mu R_g. \end{equation}
	This estimate is a part of our proof of the Main Theorem, and we will revisit it in Section \ref{sect_proofoffrontiertimes}.

	\begin{env_rem} \label{rem_trans}
		Our proof uses the completeness of the metric, by way of the Shape Theorem.  However, a version of Theorem \ref{transientgeodesics} is true for all metrics $g \in \Omega_+$, regardless of completeness.  In that version, we set $T_g = \sqrt{\sup_{\hat K_g} |g(x)| } R_{K_g}$.  Since this involves the maximum value of the metric over the very large random set $\hat K_g$, it is a very poor estimate for the exit time.  Nonetheless, even this weaker estimate implies that $\{ \mbox{$\gamma$ is bounded} \} \subseteq \{ \mbox{$\gamma$ is not minimizing} \}$.
	\end{env_rem}
			
	\subsection{Conjugate Points Along Minimizing Geodesics}

	The following result demonstrates that minimizing geodesics starting from the same point do not meet again. This is a standard theorem in differential geometry. We include its proof for completeness, and also to introduce the concept of a Jacobi field, which will play a major role in the sequel.  
	
	The idea of the proof is that if two minimizing geodesics $\gamma_v$ and $\gamma_w$ do meet at a point $x = \gamma_v(t) = \gamma_w(t)$, then one can take a shorter path to $\gamma_v(t+\epsilon)$ by following a curve near $\gamma_w$, and ``rounding the corner'' at $x$.  This idea is made precise using Jacobi fields; see Chapter 10 of Lee \cite{lee1997rmi} for an overview. We will revisit Jacobi fields when we construct bump surfaces in Section \ref{sect_bump_short}.

	\begin{env_thm} \label{minimizing_geodesics}
		For all metrics $g \in \Omega_+$, and for all $v,w \in \V_g$, the minimizing geodesics $\gamma_v$ and $\gamma_w$ meet only at the origin.
	\end{env_thm}
	\begin{proof}
		Suppose that minimizing geodesics $\gamma_v$ and $\gamma_w$ meet at some point $x \ne 0$.  Since both geodesics are minimizing, they reach $x$ at the same time $t = d_g(0,x)$, and can be extended for a small time $\epsilon$ beyond $t$. Define the variation of geodesics $\Gamma : [0,1] \times [0,t+\epsilon] \to \R^d$ by $\Gamma_\alpha(s) = \exp \big( s( (1-\alpha) v + \alpha w ) \big),$ so that $\Gamma_0$ represents the geodesic $\gamma_v$ and $\Gamma_1$ represents the geodesic $\gamma_w$.  
		
		The vector field $J(s) = \tfrac{\partial}{\partial \alpha} \Gamma_\alpha(s)|_{\alpha = 0}$ is a Jacobi field along $\gamma_v$, and vanishes at $s=0$ and $s=t$. Consequently, the point $x$ is conjugate to the origin along $\gamma_v$.  By Jacobi's theorem (Theorem 10.15 of \cite{lee1997rmi}), the geodesic $\gamma_v$ is not minimizing, a contradiction.
	\end{proof}

	This phenomenon is qualitatively different than what happens in lattice models of first-passage percolation:   minimizing geodesics may meet, and once this occurs, they coalesce. It is likely that the considerable power of Busemann functions \cite{busemann1955geometry} can be deployed in our setting to study minimizing geodesics, as they have been in lattice first-passage percolation \cite{hoffman2005coexistence, hoffman2008geodesics, cator2009busemann, damron2012busemann}.

	\part{Auxiliary Theorems and the Proof of the Main Theorem}

	\section{The Point of View of the Particle} \label{sect_POV}

	Consider a particle traveling in the random Riemannian environment $g$, experiencing no external forces. The geodesic equation $\ddot \gamma^k = - \Gamma_{ij}^k(g, \gamma) \dot \gamma^i \dot \gamma^j$ describes the trajectory of the particle. Instead of propagating the particle forward, an equivalent perspective is to leave its position fixed, and propagate the environment $g$ backwards. The resulting environment $\sigma_t g$ represents the environment as seen from the point of view of the particle.

	We remark that \textbf{no independence or moment assumptions} are needed for any of the results in this section, which rely only on the fact that the measure $\PP$ is invariant under isometries of $\R^d$.  In Section \ref{law_POV}, we will specialize to the case $d=2$; however, some variant of Theorem \ref{POVthm} should be true in arbitrary dimensions $d \ge 3$.

	\subsection{Scenery Along a Geodesic}

	Let $\xi \in \Omega$ be a (possibly random) symmetric tensor field, representing the ``scenery'' of the environment.  If $\xi$ is random, it need not be independent from $g$; it could be, for example, the Ricci curvature tensor $R_{ij}$, or the metric tensor itself.  At time zero, the particle ``sees'' $\xi(0)$, and as it evolves, its perspective remains centered at $\xi\big( \gamma(t) \big)$. As the particle rotates, the scenery tensor rotates accordingly.
	
	% We mention that there is a rich literature on (discrete) random walks in random scenery \cite{den1988mixing, benjamini1996distinguishing, gantert2006deviations}.
	
	The standard basis vectors $\E_1, \cdots\!, \E_d$ form an orthogonal basis of $\R^d$ with respect to the Euclidean metric.  Let $\E_i(g,t)$ denote the parallel translate of the Euclidean basis vector $\E_i$ along $\gamma$, where the parallel translation is with respect to the Levi-Civita connection associated with the random Riemannian metric $g$.  Let $\OO_t \in \SO(d)$ represent the flow on the orthogonal frames sending the standard basis vectors $\E_i$ to the transformed basis vectors $\E_i(t)$ (that is, $\OO_t$ is the random matrix with columns $\E_i(t)$).
	
	We now define the random POV flow $\sigma_t^{(g)} : \Omega \to \Omega$ on a $2$-tensor field; we henceforth suppress the dependence on the metric $g$ and write $\sigma_t = \sigma_t^{(g)}$. The tensor $\sigma_t \xi$ represents the scenery as viewed from the point of the view of the particle at time $t$. The flow $\sigma_t \xi \in \Omega$ is defined using the flow of frames $\OO_t$ by
		\begin{equation} \label{sigmadef}
			(\sigma_{t} \xi)_{ij}(u) = \xi_{ab}\big( \gamma(g,t) + \OO_t u \big) [\OO_t]_i^a [\OO_{t}]_j^b, \end{equation}
	where we use the Einstein convention of summing over the repeated indices $a$ and $b$.  If the tensor field is conformal (i.e., $\xi_{ij}(x) = \lambda(x) \delta_{ij}$), then the above formula simplifies to $(\sigma_t \xi)_{ij}(u) = \lambda\big( \gamma(g,t) + \OO_t u \big) \delta_{ij}$.  We emphasize that the scenery $\sigma_t \xi$ at time $t$ is an $\Omega$-valued random variable, since the random metric enters in the definition of the flow $\sigma_t = \sigma_t^{(g)}$.
			
	\begin{env_lem} \label{sigmatxi}
		The map $(t, g, \xi) \mapsto \sigma_t^{(g)} \xi$ is jointly continuous.
	\end{env_lem}
	\begin{proof}
		By Lemma \ref{geodesics_cts}, the map $(g,t) \mapsto \gamma(g,t)$ is jointly continuous.  Parallel translation is continuous,\footnote{In coordinates, parallel translation is described by the Christoffel symbols.  These are jointly continuous in both the metric and position by Lemma \ref{Lipest_lemma}.} so $(g,t) \mapsto \OO_t(g)$ is jointly continuous.  The proof easily follows from the definition \eqref{sigmadef} of $\sigma_t$.
	\end{proof}

	Since the flow $\sigma_t$ is defined using rigid transformations of the plane, the relative Euclidean distance along the geodesic is preserved.  With probability one,
		\begin{equation} \label{relativedistance}
			\big| \gamma(g, s) - \gamma(g, s') \big| = \big| \gamma(\sigma_t g, s-t) - \gamma(\sigma_t g, s'-t) \big|, \end{equation}
	for all times $s,s', t$.  We will use this fact in Section \ref{sect_exittime}.
	
	The next result implies that if $\gamma$ is a bounded geodesic, then it revisits (essentially) the same scenery infinitely often. The argument is simple and robust.
	
	\begin{env_pro} \label{scenery_limit}
		Let $g \in \Omega_+$ be any Riemannian metric, and let $\xi \in \Omega$ be a scenery $2$-tensor (possibly dependent on $g$).  If $\gamma(g,\cdot)$ is a bounded geodesic, then the family $\{ \sigma_t \xi \}$ is relatively compact.   %For every $g \in B$; family of measures $\{ \PP \circ \sigma_t^{-1} \}$ is tight.
	\end{env_pro}
	\begin{proof}
		For any pair $(x,\OO) \in \R^d \times \SO(d)$, let $\sigma_{x,\OO} \xi$ represent the scenery viewed from the position $x$ with orientation $\OO$.  i.e., $(\sigma_{x,\OO} \xi)_{ij}(u) = \xi_{ab}\big( x + \OO u \big) \OO_i^a \OO_j^b$.  Using this notation, the flow $\sigma_{\gamma(t), \OO_t}$ is $\sigma_t$ as we defined it in \eqref{sigmadef}.  
		
		Let $g$ be a metric for which $\gamma$ is bounded, and let $K_g \subseteq \R^d$ be a compact set which contains the forward trajectory of the geodesic (i.e., $\gamma(t) \in K_g$ for all $t \ge 0$).  Define the family of possible sceneries $\mathcal K_g = \{ \sigma_{x,\OO} \xi : \mbox{$x \in K_g$, $\OO \in \SO(d)$} \} \subseteq \Omega.$ As in Lemma \ref{sigmatxi}, the map $(x,\OO) \mapsto \sigma_{x, \OO} \xi$ is jointly continuous, so the family $\mathcal K_g$ is a compact subset of $\Omega$.  
		
		The geodesic $\gamma$ is trapped in $K$ for all forward time, so $\sigma_t \xi = \sigma_{\gamma(t), \OO_t} \xi \in \mathcal K_g$.  Since $\mathcal K_g$ is compact, the family $\{\sigma_t \xi\}$ is relatively compact.
	\end{proof}
	
%	The limit of $\sigma_t \xi$ should only exist in contrived situations (e.g., if $\xi$ is a constant scenery field).  Nonetheless,   Proposition \ref{scenery_limit} is almost as useful, since a geodesic revisits (essentially) the same scenery infinitely often so $\sigma_t \xi$ has many accumulation points.
		
	\subsection{The POV of the Particle along a Bounded Geodesic}
	
	We now specialize to the case where the scenery tensor $\xi$ is the exactly the random metric $g$ itself, i.e., $\xi = g$.  For each $t$, $g \mapsto \sigma_t^{(g)} g$ is a \emph{random flow} with \emph{random initial conditions}.  By the time-invertibility of the geodesic equation, the transformation $\sigma_t$ is invertible: if $g' := \sigma_t^{(g)} g$, then $g = \sigma_{-t}^{(g')} g'$.
	
	With probability one, the flow $\sigma_t$ preserves the space $\Omega_+$ of Riemannian metrics on $\R^d$.  Consequently, with probability one, $\sigma_t g$ is a random Riemannian metric for all $t$. Clearly, the random flow $\sigma_t$ preserves the geodesic events $\{ \mbox{$\gamma$ is bounded} \}$, $\{ \mbox{$\gamma$ is unbounded} \}$, $\{ \mbox{$\gamma$ is recurrent} \}$, $\{ \mbox{$\gamma$ is transient} \}$ and $\{ \mbox{$\gamma$ is minimizing} \}$.  i.e., if $A$ denotes any of these events, then $g \in A$ if and only if $\sigma_t g \in A$.
	
	Lemma \ref{sigmatxi} implies that, with probability one, the function $t \mapsto \sigma_t g$ is continuous.  This stochastic process $\sigma_t g$, which we call the \emph{(forward) POV process}, has a very complicated correlation structure, and is assuredly not a Markov process.
	
%	Kim and LaGatta \cite{kim2012bounded} show that $\gamma$ is bounded with positive probability.  The environment along a bounded geodesic is quite boring:  the geodesic remains trapped in a compact set, and can only visit a compact set of possible environments (cf. Proposition \ref{scenery_limit}).  It immediately follows that, conditioned on the event $\{\mbox{$\gamma$ is bounded}\}$, the family $\{\sigma_t g\}$ is a.s. relatively compact:	
%		\begin{equation} \label{eqn_rc}
%			\PP\big( \mbox{the family $\{ \sigma_t g \}$ is relatively compact} \big| \mbox{$\gamma$ is bounded} \big) = 1. \end{equation}
%	We remark that \eqref{eqn_rc} holds for any dimension $d \ge 2$.
	
	We remark that \emph{a priori} there could be different behavior in the forward and backward directions.  For example. the forward direction of $\gamma$ could be bounded, while the backward direction remains unbounded.  It is an open question to determine if this phenomenon could occur with a positive probability.

	\subsection{The Law of the POV Process} \label{law_POV}
	
	Henceforth, we restrict ourselves to the case $d = 2$.  Let $v_t = \dot\gamma(t) / |\dot\gamma(t)|$ denote the direction of the tangent vector along the geodesic $\gamma$ (note that $v_0 = \E_1$) and write $v_t^\perp = \left(-v_t^2, v_t^1\right)$ for its perpendicular vector \emph{with respect to the Euclidean metric}.  The vectors $v_t$ and $v_t^\perp$ form a basis of the tangent space $T_{\gamma(t)} \R^2$ with the same orientation as the standard basis $(\E_1,\E_2)$.\footnote{The random vectors $v_t$ and $v_t^\perp$ are orthogonal with respect to the flat Euclidean metric $\delta$, not the random Riemannian metric $g$.}  %Let $\O_{v_t} \in \SO(2)$ be the family of rotation matrices which take $(\E_1, \E_2) \mapsto (v_t, v_t^\perp)$.
	
	Let $\OO_t(g) = \big[ v_t \big| v_t^\perp \big] \in \SO(2)$ represent the orthogonal matrix which changes the orientation of $\R^2$ to point in the direction $v_t$.  We define $\sigma_t g$ as in \eqref{sigmadef}, i.e., $(\sigma_t g)_{ij}(u) = g_{ab}\big( \gamma(t) + \OO_t u \big) [\OO_t]_i^a [\OO_t]_j^b.$  When we apply the Euclidean inner product to two tangent vectors $w, w' \in T_u \R^2$, we have $\langle w, (\sigma_t g)(u) w' \rangle = \langle \OO_{t} w, \, g\big(\gamma(t) + \OO_t u \big) \, \OO_{t} w' \rangle$. i.e., we either use the interpretation that the environment is fixed and the vectors evolve (right side of equation), or that the vectors are fixed and the inner product evolves (left side).  
	
	The variable $u \in \R^2$ represents the displacement from $\gamma(t)$ in the tangent direction $v_t$, so that $(\sigma_t g)(0) = g(\gamma(t))$ always represents the metric at $\gamma(t)$, and the axial directions $\partial/\partial u^1$ and $\partial/\partial u^2$ correspond to the directions $v_t$ and $v_t^\perp$, respectively.  Under the metric $\sigma_t g$, the geodesic of interest is always at the origin pointing in direction $\E_1$.  
	
	Let $\PP_t := \PP \circ \sigma_t^{-1}$ be the push-forward of the law $\PP$ under the transformation $g \mapsto \sigma_t^{(g)} g$.  The measure $\PP_t$ is the law of the random Riemannian metric $\sigma_t g$, and is uniquely defined by the change-of-variables formula:
		\begin{equation}
			\int_\Omega f(g) \sD \PP_t(g) := \int_\Omega f(\sigma_t g) \sD \PP(g) \end{equation}
	for any continuous, bounded function $f : \Omega \to \R$.  By the bounded convergence theorem, it is clear that
		\begin{equation}
			\mbox{the measure-valued flow $t \mapsto \PP_t$ is weakly continuous.} \end{equation}
	That is, if $t_n \to t$, then the measures $\PP_{t_n}$ converge weakly to $\PP_t$. We are now ready to state the main theorem of the section.
	
	\begin{env_thm}[The Law of the POV Process $\sigma_t g$] \label{POVthm}
		Suppose $d=2$.  The law $\PP_t$ of $\sigma_t g$ is absolutely continuous with respect to $\PP$, and its Radon-Nikodym derivative $\tfrac{\D \PP_t}{\D \PP}(g)$ equals
			\begin{equation} \label{rhotg}
				\rho_t(g) := \exp \!\left( \int_{-t}^0 \Big( \big\langle \nabla \log \det g(\gamma(s)), \dot \gamma(s) \big\rangle + 3 \frac{\langle \ddot \gamma(s), \dot \gamma(s) \rangle}{\langle \dot \gamma(s), \dot \gamma(s) \rangle} \Big) \sD s\right) \end{equation}
		almost surely.  This implies that
			\begin{equation} \label{flowformula}
				\EE_t f := \int_\Omega f(\sigma_t g) \sD \PP(g) = \int_\Omega f(g) \rho_t(g) \sD \PP(g) \end{equation}
		for any integrable $f : \Omega \to \R$. The function $(t,g) \mapsto \rho_t(g)$ is jointly continuous in $t$ and $g$.  For each $g \in \Omega_+$, $t \mapsto \rho_t(g)$ is differentiable in $t$.
	\end{env_thm}
	
	We emphasize that \textbf{no independence or moment assumptions} are required for this theorem, only the isometric invariance of the law $\PP$ of the random Riemannian metric $g$. As we saw in the proof of Proposition \ref{scenery_limit}, the infinite-dimensional random flow $\sigma_t$ on $\Omega$ can be reduced to a finite-dimensional random flow $\big(\gamma(t), \OO_t \big)$ on the Lie group of isometries of $d$-dimensional space Euclidean space.  This type of reasoning is the key insight of Geman and Horowitz \cite{geman1975random}. As stated, their argument only applies to the case of an abelian group of transformations. The group of rotations is abelian only in two dimensions, and we are able to extend their argument to the $2$-dimensional Euclidean group $\operatorname{Euc}(2) \cong \R^2 \rtimes O(2)$. The case $d \ge 3$ is still an open question, but we conjecture that a similar theorem should hold in general. We follow Zirbel's outline of the Geman-Horowitz method \cite{zirbel2001lagrangian} to prove Theorem \ref{POVthm}. 
	
	\begin{env_rem}
		A similar theorem holds for other scenery tensors $\xi$ along a geodesic.  For example, the Riemann or Ricci curvature tensors, the scalar curvature, or the difference of two connection forms.
	\end{env_rem}

%	We now sketch the proof; see Section \ref{proof_flowthm} for details.  Let $\nu$ be a probability measure which is absolutely continuous to Haar measure on the isometry group $\SE(2)$, and let $f$ be a real-valued observable of the metric.  We first multiply $\EE_t f = \int f(\sigma_t g) \sD \PP$ by $1 = \int \D \nu(x,\O)$, and inter 
	
%	Liouville's theorem \cite{arnold1978ordinary} gives a formula for the change of measure, involving the divergence of the vector field which generates the flow; this is the basis of the formula \eqref{rhotg} for $\rho_t$.  

	Outline of proof: For each metric $g \in \Omega_+$, we define a vector field $U_g$ on the space of possible locations and orientations of the geodesic; flow lines for this vector field $U_g$ correspond to geodesics of the metric $g$.  We then prove a number of identities about geodesics in these coordinates (cf. Lemma \ref{flowlemma}); it is here that we use the assumption that $d=2$.  The proof then follows from a series of careful calculations.  For full details, see Section \ref{proof_flowthm}.
	
%	\begin{env_rem} \label{rem_SOd} In dimension $d = 2$, it suffices to track only the velocity vector $v_t$; from this, we can uniquely construct an orthonormal basis $(v_t, v_t^\perp)$.  In higher dimensions $d > 2$, however, we have to track the evolution of the entire orthonormal frame $(\E_1, \dots, \E_d)$ of $T_0 \R^d$.  Consequently, the correct phase space to study is $\R^d \times \SO(d)$.  The method of Geman and Horowitz breaks down, since the group $\SO(d)$ is not abelian when $d > 2$.  In particular, equation \eqref{transformXV} does not hold.   \end{env_rem}
	
%	Statement \eqref{eqn_rc} says that, on the event that $\gamma$ is bounded, the POV process $\sigma_t g$ is relatively compact.  Combining this with the exact expression \eqref{rhotg} for the Radon-Nikodym derivative $\rho_t(g)$, one should be able to deduce something about the random geometry along geodesics.
	
	In Appendix \ref{app_geomgeod}, we use Theorem \ref{POVthm} to prove a simple geometric corollary: with probability one, the geodesic $\gamma$ contains neither straight line segments nor circular arcs. This is a formalization of the simple intuition that a random geometry is locally non-Euclidean. The Euclidean character of the random geometry manifests itself only in terms of global symmetries, as expressed by the Shape Theorem.

	\section{The Exit Time Process} \label{sect_exittime}
	
	On the event $\{ \mbox{$\gamma$ is bounded} \}$, any sequence of times $t_n$ has a subsequence $t_{n_k}$ along which the environments $\sigma_{t_{n_k}} g$ converge. On the event $\{ \mbox{$\gamma$ is unbounded} \}$, the behavior of $\sigma_t g$ is quite different.  For any radius $r \ge 0$, let $\tau_r(g)$ denote the exit time of the geodesic $\gamma(g,\cdot)$ from the Euclidean ball $B(0,r)$ (defined formally in Section \ref{subsect_propertiesoftaur}).  The geodesic $\gamma$ is unbounded if and only if $\tau_r < \oo$ for every value of $r$.  We call $r \mapsto \tau_r(g)$ the \emph{exit time process} of $g$.
	
%	\begin{figure}[h!]
%			\includegraphics{figure-exittime.pdf}
%		\caption{$\gamma$ exits the ball $B(0,r)$ at time $\tau_r$.}
%	\end{figure}

%	It is easy to see that for almost every $g$, the process $\tau_r$ is defined for $r$ in a non-empty interval $\big[0,\Rmax(g)\big)$, and that $\tau_r$ is a.s. upper semicontinuous (Lemma \ref{tau_rcontinuous}), hence a stochastic process with jumps.
	
	In Section \ref{subsect_povexittimes}, we will look at the point of view of the environment at exit times, i.e., the metric-valued random variable $\sigma_{\tau_r} g$.  We will examine its law $\PP_{\tau_r}$, and in Theorem \ref{POVthm_taur} we will prove that $\PP_{\tau_r}$ is absolutely continuous with respect to $\PP$.  Unlike in Theorem \ref{POVthm}, its Radon-Nikodym derivative is more complicated than just $\rho_{\tau_r}$:  we have to integrate over all past times $t$ against the ``history measure'' $\ell_t(g,\D t)$, which we define in Section \ref{subsect_povexittimes}.
	
	\subsection{Properties of the Exit Time Process} \label{subsect_propertiesoftaur}
		
	For any $r \ge 0$, let $\F_r := \F_{B(0,r)}$ be the $\sigma$-algebra generated by the metric in an infinitesimal neighborhood of the Euclidean ball $B(0,r)$ (defined formally in footnote \ref{FD}).  It is not hard to see that $\F_r$ defines a right-continuous filtration on the probability space $(\Omega, \F, \PP)$ (cf. \cite[Section 7.2]{durrett1996probability}).  
	
	For every $g \in \Omega_+$, let $\tau_r(g)$ denote the (forward) exit time of the geodesic $\gamma := \gamma_{\E_1}$ from the Euclidean ball $B(0,r)$ of radius $r$ centered at the origin: $\tau_r(g) := \inf \!\left\{ t > 0 : |\gamma(g,t)| > r \right\}$. We follow that the convention that if $\gamma$ is trapped in the ball $B(0,r)$, then we set $\tau_r = \oo$.\footnote{That is, $\tau_r = \oo$ if $|\gamma(t)| \le r$ for all $t \ge 0$.}  
	
	Let $\Rmax(g) = \sup |\gamma(g,t)|$ be the maximum (Euclidean) distance the geodesic reaches from the origin.  Clearly, $\tau_r < \oo$ if and only if $r < \Rmax$.  For every $g \in \Omega_+$, the domain of the process $r \mapsto \tau_r(g)$ is the interval $\big[0,\Rmax(g)\big)$.  The geodesic $\gamma$ is unbounded if and only if $\Rmax = \oo$.  
	
	\begin{env_lem} \label{tau_rcontinuous}
		For every metric $g \in \Omega_+$, we have $\Rmax(g) > 0$, which implies that $\PP( \Rmax > 0 ) = 1$.  Furthermore, the exit time process $r \mapsto \tau_r(g)$ and the exit time POV process $r \mapsto \sigma_{\tau_r} g$ are both upper-semicontinuous on their domain $\big[0,\Rmax(g)\big)$.  Both processes are adapted to the filtration $\F_r$.
	\end{env_lem}
	\begin{proof}
		Since $\gamma$ starts at the origin in the direction $\E_1$, we estimate $\gamma(t) = \E_1 t + O(t^2)$, where the constant on the second-order term depends on the metric $g$.  Consequently, the geodesic $\gamma$ must exit all sufficiently small balls around the origin, so $\tau_r(g) > 0$ for sufficiently small $r$.  This implies that $\Rmax > 0$ for all $g$.
	
		It is clear that the function $r \mapsto \tau_r(g)$ is strictly increasing for every Riemannian metric $g \in \Omega_+$.  Since $r \mapsto \tau_r$ is monotone, limits exist from below and above. Let $\psi_r(t) = |\gamma(t)| - r$, so that $\tau_r = \inf \psi_r^{-1}( (0, \oo) )$.  The upper semi-continuity of $r \mapsto \tau_r$ follows from the following simple real-analytic lemma.
		
		\begin{env_sublem} \label{lem_uppersemicty}
			Let $\psi_r(t)$ be a jointly continuous, real-valued function of $r$ and $t$, and let $B$ be an open set in $\R$.  The function $r \mapsto \inf \psi_r^{-1}(B)$ is upper semi-continuous on its domain.  
		\end{env_sublem}
		\begin{proof}[Proof of Sublemma \ref{lem_uppersemicty}]
			Write $\psi(r,t) = \psi_r(t)$, and consider the open set $\tilde B = \psi^{-1}(B) \subseteq \R^2$.  Suppose that $t(r) := \inf \psi_r^{-1}(B)$ is finite, and let $r_n \to r$.  Let $\hat t > t(r)$.  Since $(r, \hat t) \in \tilde B$, we can find an open neighborhood $U$ of $(r, \hat t)$ contained in $\tilde B$.  Furthermore, $(r_n, \hat t) \in U$ for all but finitely many $n$, so $\limsup_{n \to \oo} t(r_n) \le \hat t$.  Taking the limit $\hat t \downarrow t(r)$ proves that $t(r)$ is upper semi-continuous.
		\end{proof}
		
		For all $g \in \Omega_+$, the function $t \mapsto \sigma_t g$ is continuous in $t$.  By setting $t = \tau_r$, the upper semi-continuity of $r \mapsto \sigma_{\tau_r} g$ follows immediately.  By right-continuity of the filtration $\F_r$, the exit time process $(r,g) \mapsto \tau_r(g)$ is adapted to $\F_r$. This proves Lemma \ref{tau_rcontinuous}.
	\end{proof}
	
	There are plenty of metrics for which $\Rmax = \oo$, including the Euclidean metric $\delta$, though it is an open question to determine whether $\PP( \Rmax = \oo ) = \PP( \mbox{$\gamma$ is unbounded} ) = 0$ or $>0$.
	
%	\begin{env_pro} \label{prop_taurpos}
%		For all $r > 0$, 
%			$$0 < \PP( \tau_r < \oo) = \PP( \Rmax > r) < 1.$$
%	\end{env_pro}
%	\begin{proof}
%		It is clear that $\PP( \tau_r < \oo) = \PP( \Rmax > r)$.  The theorem of Kim and LaGatta \cite{kim2012bounded} says that $\PP(\Rmax \le r) > 0$, so $\PP( \Rmax > r) = 1 - \PP(\Rmax \le r) < 1$.
	
%		For each $r > 0$, the function $g \mapsto \gamma(g, r+1)$ is continuous, so the set $U = \{g : |\gamma(g,r+1)| > r \}$ is an open subset of $\{ \Rmax > r\}$.  The Euclidean metric $\delta$ clearly belongs to $U$, so by strict positivity of the measure $\PP$ (Theorem \ref{PPproperties}.\ref{assumption_totallypositive}), $0 < \PP(U) \le \PP(\Rmax > r)$.  This completes the proof.
%	\end{proof}

	\subsection{The Law of the Exit Time POV Process} \label{subsect_povexittimes}

	For each fixed, deterministic $t$, the environment $\sigma_t g$ is a random Riemannian metric with law $\PP_t$.  Theorem \ref{POVthm} states that this measure $\PP_t$ is absolutely continuous with respect to $\PP$ on $\Omega$, and has a nice Radon-Nikodym derivative $\rho_t(g)$ (defined in \eqref{rhotg}).  
	
	In this section, we investigate the exit time POV process $r \mapsto \sigma^{(g)}_{\tau_r(g)} g$.  Now the radius $r$ is fixed, but the exit time $\tau_r$ is random.  The random Riemannian metric $\sigma_{\tau_r} g$ represents the POV of the particle as it first exits the ball $B(0,r)$. Let $\PP_{\tau_r}(U) := \PP\big( (\sigma_{\tau_r})^{-1} U \big| \tau_r < \oo \big)$ denote the law of $\sigma_{\tau_r} g$ (conditioned on the event that the exit time $\tau_r$ is finite).  Note that by construction, $\PP_{\tau_r}(\tau_r = \oo) = 0$.  We include this conditioning in order to have a well-defined random variable $g \mapsto \sigma_{\tau_r(g)}^{(g)} g$ on the entire space $\Omega_+$.
		
	For each $g \in \Omega_+$ and each $r \ge 0$, the function $t \mapsto \tau_r(\sigma_{-t} g) - t$ is well-defined and upper semi-continuous.  Therefore, for any $r \ge 0$, the ``history measure'' $\ell_r(g,\D t) := \delta\big(\tau_r(\sigma_{-t} g) - t \big) \sD t$ is well-defined.
		
	\begin{env_thm}[The Law of the Exit Time POV Process $\sigma_{\tau_r} g$] \label{POVthm_taur}
		Suppose $d=2$.  The law $\PP_{\tau_r}$ of $\sigma_{\tau_r} g$ is absolutely continuous with respect to $\PP$. For any $r \ge 0$ and integrable $f : \Omega \to \R$,
			\begin{equation} \label{flowformula_taur0}
				\EE_{\tau_r}(f) \cdot \PP(\tau_r < \oo) = \int_{\{\tau_r < \oo\}} f(\sigma_{\tau_r} g) \sD \PP(g) = \int_\Omega f(g) \zeta_r(g) \D \PP(g), \end{equation}
		where $\zeta_r(g) = \int_0^\oo \rho_t(g) \, \ell_r(g,\D t)$, for $\rho_t(g)$ the density from the POV Theorem, and the history measure $\ell_r(g,\D t)$ defined above. For a.e. random Riemannian metric $g$, the function $r \mapsto \zeta_r(g)$ is a stochastic process with jumps.\footnote{That is, $r \mapsto \zeta_r(g)$ is an upper semicontinuous function which is adapted to $\F_r$.}
	\end{env_thm}
	
	More generally, if $f$ and $F$ are any integrable functions, then for any $r \ge 0$,
			\begin{equation} \label{flowformula_taur}
				\int_{\{\tau_r < \oo\}} f(\sigma_{\tau_r} g) F(g) \sD \PP(g) = \int_\Omega f(g) \left( \int_0^\oo F(\sigma_{-t} g) \rho_t(g) \, \ell_r(g, \D t) \right) \D \PP(g), \end{equation}

	We will prove this more general formula below; \eqref{flowformula_taur0} follows from setting $F = 1$. By setting $f = 1$, it follows that $\PP(\tau_r < \oo) = \EE \zeta_r.$  %By Proposition \ref{prop_taurpos}, this number is strictly less than $1$, and is monotone decreasing in $r$.  \newline
%	It is an open question whether this is bounded below by some constant (i.e., $\gamma$ is unbounded with positive probability). \newline
	
\begin{figure}[h!] \label{figure-manyoldorigins}
			\includegraphics{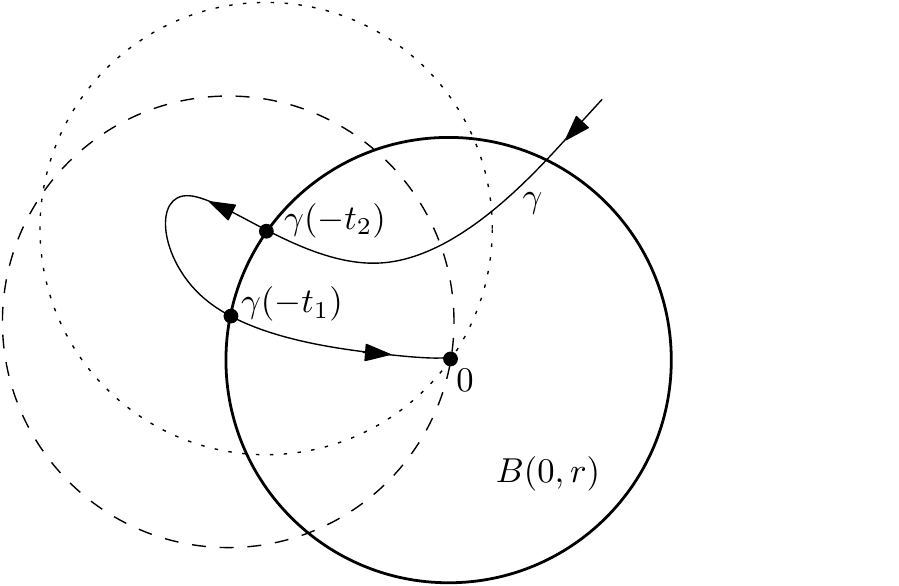}
		\caption{In this typical example, the set $\T_r(g)$ consists of two historical exit times $t_1$ and $t_2$; there are two corresponding historical metrics $h_1 = \sigma_{-t_1} g$ and $h_2 = \sigma_{-t_2} g$.  The curve in this figure is the geodesic $\gamma$.  The dashed and dotted lines indicate the balls $B\big(\gamma(-t_i), r\big)$ centered at the old origins.}
	\end{figure}

	Before we prove this theorem, we introduce the concept of ``historical metrics.''  On the right side of formula \eqref{flowformula_taur}, the metric $g$ should be thought of as the result of a POV transformation.  That is, $g = \sigma_t h = \sigma_t^{(h)} h$, for some Riemannian metric $h$ and exit time $t = \tau_r(h)$.  The metric $h$ is called a historical metric of $g$ at horizon $r$.  This transformation is reversible, in the sense that for this fixed time $t$, $h = \sigma_{-t}^{(g)} g$.  However, neither the history time $t$ nor the historical metric $h$ is necessarily determined by the metric $g$.  See Figure \ref{figure-manyoldorigins} for a typical example.
	
	After the POV transformation, the origin $0$ under the historical metric $h$ is transformed into the ``old origin'' $x = \gamma(g,-t)$.  Similarly, the direction $\E_1$ is rotated to be parallel with $v = -\dot\gamma(g,-t)$.  Under the transformed metric, the geodesic $\gamma_{x,v}(g,\cdot)$ first exits the ball $B(x,r)$ at time $t = \tau_r(h)$.  Since we can write $h = \sigma_{-t} g$, this in fact characterizes the historical metrics.
	
	\begin{env_def}
		Let $g \in \Omega_+$.  Define the closed set of \emph{historical exit times} (at horizon $r$) by $\T_r(g) = \{ t \ge 0 : \tau_r(\sigma_{-t} g) = t \}.$ We say that $h$ is a \emph{historical metric} of $g$ (at horizon $r$) if $h = \sigma_{-t} g$, for some historical exit time $t \in \T_r(g)$. 
	\end{env_def}
	
	The support of the history measure is exactly the set $\T_r(g)$ of historical exit times.  When this set contains isolated points, the history measure has point masses at those times.
	
	The metric $g$ is a result of a POV transformation $\sigma_t h$ if and only if the set $\T_r(g)$ is non-empty.  We remark that it is quite common for the set $\T_r(g)$ to contain multiple values.  For the example in Figure \ref{figure-manyoldorigins}, the set $\T_r(g)$ contains exactly two times $t_1$ and $t_2$, and there is a different historical metric corresponding to each of them.  It is an easy consequence of Proposition \ref{nonconstkappa} that $\T_r(g)$ does not contain any interval, though it may be uncountable.
	
	Having set up the right definitions, the proof of Theorem \ref{POVthm_taur} is quite easy.  On the left side of formula \eqref{flowformula_taur}, for almost every $g$, we introduce an integral over all possible values of $t = \tau_r(g)$.  We then introduce an appropriate approximation to the indicator function $\delta\big(\tau_r(g)-t\big)$, so that we may use Fubini's theorem and interchange the integrals.  In \eqref{flowformula_taur_proof1.5}, for each $t$, we integrate over the metrics for which $t = \tau_r$ is the exit time, where we assign weight $f(\sigma_t g)$ to the metric $g$.
	
	Still keeping $t$ fixed, we make the change of coordinates $g \mapsto \sigma_{-t} g$ in the inner integral.  The function $\rho_{t}(g)$ of Theorem \ref{POVthm} is the Jacobian for this coordinate change.  This is the essential step of the proof, and the remainder is tracing back the original steps we made.

	\begin{proof}[Proof of Theorem \ref{POVthm_taur}]
		We will prove formula \eqref{flowformula_taur} for the case that $f$ and $F$ are non-negative, bounded and continuous functions on $\Omega$.  The general statement follows from standard approximation arguments.  Formula \eqref{flowformula_taur0} follows from setting $F=1$ and writing $\zeta_r(g) = \int_0^\oo \rho_t(g) \ell_r(g,\D t)$.	
		
		We first make an approximation to the Dirac $\delta$-function.  For any $\epsilon > 0$, define $\delta^\epsilon(u) = \tfrac 1 \epsilon 1_{[-\epsilon,0]}(u)$.  Consider the function $\ell_r^\epsilon(g, t) = \delta^\epsilon( \tau_r(\sigma_{-t} g) - t) $, so that
		\begin{equation} \label{ell_approx}
			\mbox{for all $g \in \Omega_+$, the measure $\ell_r^\epsilon(g, t) \sD t$ converges weakly to $\ell_r(g, \D t)$.} \end{equation}
		We do this so that we can work with the density function $\ell_r^\epsilon$, and interchange integrals.\footnote{Technically, the convergence in \eqref{ell_approx} is weak-$*$ convergence. This means that for any continuous, bounded function $f$, $\int f(t) \ell_r^\epsilon(g,t) \sD t$ converges to $\int f(t) \, \ell_r(g,\D t).$} Rewrite the left side of \eqref{flowformula_taur} by integrating over all possible (finite) values of $\tau_r(g)$.  We then introduce the approximation $\delta^\epsilon$, take out the limit, and interchange the integrals:
		 	\begin{eqnarray}
				\int_{\{\tau_r < \oo\}} F(g) f(\sigma_{\tau_r} g) \sD \PP(g) &=& \int_{\Omega} F(g) \int_0^\oo f(\sigma_t g) \, \delta\big(\tau_r(g) - t\big) \sD t \sD \PP(g) \nonumber \\
				&=& \int_\Omega F(g) \, \lim_{\epsilon \downarrow 0} \int_0^\oo f(\sigma_t g) \, \delta^\epsilon\big(\tau_r(g) - t\big) \sD t \sD \PP(g) \label{flowformula_taur_proof0} \\
				&=& \lim_{\epsilon \downarrow 0} \int_\Omega F(g) \int_0^\oo f(\sigma_t g) \, \delta^\epsilon\big(\tau_r(g) - t\big) \sD t \sD \PP(g). \label{flowformula_taur_proof1} \\
				&=& \lim_{\epsilon \downarrow 0} \int_\Omega F(g) \int_0^\oo f(\sigma_t g) \, \delta^\epsilon\big(\tau_r(g) - t\big) \sD \PP(g) \sD t. \label{flowformula_taur_proof1.5}
			\end{eqnarray}
		We must justify each of the steps.  First, Lemma \ref{sigmatxi} implies that $t \mapsto \sigma_{-t} g$ is continuous for each $g$.  By assumption, $f$ is bounded and continuous, so $t \mapsto f(\sigma_t g)$ is bounded and continuous.  Statement \eqref{flowformula_taur_proof0} follows from weak convergence of the functions $\delta^\epsilon$ to the Dirac measure $\delta_r$.  Statement \eqref{flowformula_taur_proof1} follows from the dominated convergence theorem, and \eqref{flowformula_taur_proof1.5} from Fubini's theorem.
		
		Next, we apply the POV Theorem (Theorem \ref{POVthm}).  Making the change of variables $g \mapsto \sigma_{-t} g$ in \eqref{flowformula_taur_proof1.5} yields
			\begin{equation} \label{flowformula_taur_proof2}
				\lim_{\epsilon \downarrow 0} \int_0^\oo \int_\Omega F(\sigma_{-t} g) f(g) \ell^\epsilon(g, t) \, \rho_t(g) \sD \PP(g) \sD t, \end{equation}
		since $\ell^\epsilon(g,t) = \delta^\epsilon(\tau_r(\sigma_{-t} g) - t)$. We now undo the approximations.  We again interchange the integrals in \eqref{flowformula_taur_proof2}, and get
			\begin{equation} \label{flowformula_taur_proof3}
				\int_\Omega f(g) \, \lim_{\epsilon \downarrow 0} \int_0^\oo F(\sigma_{-t} g) \rho_t(g) \, \ell_r^\epsilon(g, t) \sD t \sD \PP(g) = \int_\Omega f(g) \int_0^\oo F(\sigma_{-t} g) \rho_t(g) \, \ell_r(g, \D t) \sD \PP(g). \end{equation}
		This is justified by the dominated convergence theorem, and applying weak convergence of measures with the continuous, bounded integrand $t \mapsto F(\sigma_{-t} g)$.  To complete the proof of \eqref{flowformula_taur}, we approximate general integrable functions $f$ and $F$ by continuous, bounded functions.\footnote{This relies on the fact that $\Omega_+$ is a Polish space, and the measure $\P$ is Radon (hence inner regular).}
	\end{proof}

	\section{The Local Markov Property} \label{sect_markov}

	We want to understand the law $\PP_{\tau_r}$ of the metric $\sigma_{\tau_r} g$ centered at the exit point $\gamma(\tau_r)$.  Unfortunately, the exit time POV process $\sigma_{\tau_r} g$ is a complicated infinite-dimensional object, and little is known about it.  Instead, we focus our attention on local observables of this process, that is, functions which are $\F_D$-measurable for some compact $D \subseteq \R^d$. The goal is to understand conditional expectations of the form $\EE(f \circ \sigma_{\tau_r}^{-1} | \F_r)$, representing the best guess of $f$ after a POV transformation, given the metric information inside the (Euclidean) ball of radius $r$. 
	
	The Local Markov Property (LMP) states that the conditional expectation depends only on the metric in the intersection of two sets: the (random) ball of radius $r$ centered at the old origin from the point of view of the particle, and $D^1$, the Euclidean $1$-ball around the set $D$. The Strong Local Markov Property (SLMP) is a similar theorem, but for random ``stopping radii'' $R$. The Inevitability Theorem formalizes the intuition ``what can happen will happen'' along an unbounded geodesic, provided certain uniformity assumptions are satisfied.
	
	We state and prove our results only in the case $d=2$, since our argument relies on Theorem \ref{POVthm_taur}, which itself relies on the POV Theorem. The key assumption is that the measure $\PP_{\tau_r}$ is absolutely continuous with respect to $\PP$; the precise form of the Radon-Nikodym derivative is irrelevant. Since similar theorems should hold in the case $d \ge 3$, the arguments of this section would generalize \emph{mutatis mutandis}.

	\subsection{Preliminaries} \label{subsect_lagcoordsatexit}

%	Before stating the Local Markov Property, we first set up some preliminary notation.  For every $g \in \Omega_+$, let $\Rmax(g) = \sup\{ |\gamma(t)| : t \ge 0 \}$ denote the maximal distance from the origin which the geodesic reaches (if $\gamma$ is unbounded, then $\Rmax = \oo$).  The initial velocity $\dot\gamma(0)$ is positive, so $\Rmax(g) > 0$ for all $g \in \Omega_+$.  

%	If $r < \Rmax$, then $\gamma$ has not yet reached its maximal distance, so it must exit the ball $B(0,r)$, hence $\tau_r < \oo$.  If $r \ge \Rmax$, then $\gamma$ will never exit the ball $B(0,r)$, so $\tau_r = \oo$.  Lemma \ref{tau_rcontinuous} states that for all $g \in \Omega_+$, the exit time process $r \mapsto \tau_r(g)$ and the exit time POV process $r \mapsto \sigma_{\tau_r} g$ are both upper semi-continuous on the domain $[0, \Rmax)$. 
	
	The geodesic initially starts at the origin in $\R^2$, then reaches the exit location $\gamma(\tau_r)$.  After the POV transformation, the exit location has become the new origin, and the old origin is shifted to the point $o_r(g) := \gamma( \sigma_{\tau_r} g, -\tau_r)$, which we call the \emph{old origin} after the exit time POV transformation.  Note that $o_0(g) = 0$ for all $g$ and $r$.  Since the transformation $\sigma_{\tau_r}$ is defined using isometries of the plane, the original ball $B(0,r)$ is transformed to the ball $B_r(g) := B\big(o_r(g), r \big)$ centered at the old origin.
	
	Let $D \subseteq \R^2$ be some compact set in the plane.  Let $D^1$ denote the $1$-neighborhood of $D$, and define the random \emph{lens-shaped set}
		\begin{equation} \label{Drdef}
			D_r(g) = D^1 \cap B( o_r(g), r) \end{equation}
	on the event $\{ r < \Rmax \}$. The random set $D_r(g)$ is a $\C$-valued random variable, where $\C$ denotes the metric space of compact subsets of $\R^2$ (equipped with the Hausdorff metric).
	
	\begin{env_lem} \label{sigmataur_lemma}
		For all $g \in \Omega_+$, the old origin process $r \mapsto o_r(g)$ and lens set process $r \mapsto D_r(g)$ are upper semicontinuous on their domain $[0, \Rmax(g))$.  Both these processes are adapted to the filtration $\F_r$.
	\end{env_lem}
	\begin{proof}
		The proof of this lemma follows immediately from the upper semicontinuity of $\tau_r$.
	\end{proof}

	\begin{figure}[h!] \label{figure-DDo}
			\includegraphics{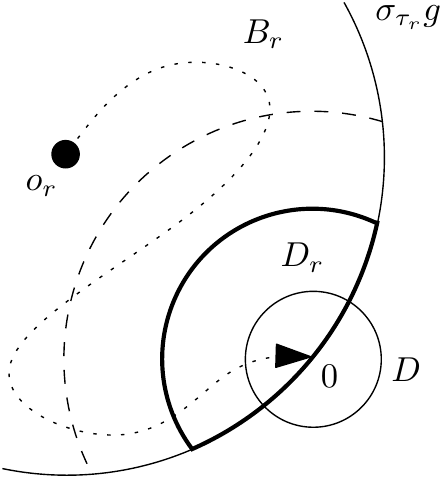}
		\caption{The sets $D$, $D_r(g)$ and $B_r(g) = B( o_r(g), r )$ from the point of view of a particle at the exit location $\gamma(\tau_r)$.}
	\end{figure}

	Let $f : \Omega \to \R$ be some local observable, i.e., an $\F_D$-measurable random variable.  Evaluating this local observable at the exit location is a tricky endeavor, since $\sigma_{\tau_r}g$ may depend on the metric in the entire ball $B(0,r)$.  %When we evaluate $f(\sigma_{\tau_r} g)$, this may also depend on the metric outside that ball.  
	
	If we switch to the point of view of the particle, then this becomes much easier to analyze.  The value of the observable $f$ at the exit location should only depend on the metric information in the neighborhood $D^1$, since the law $\PP$ is $1$-dependent.  Intuitively, the $\sigma$-algebra $\F_r$ now represents the information in the \emph{random} ball $B_r(g) = B(o_r,r)$.  The random lens-shaped set $D_r$ is exactly the restriction of $D^1$ to $B_r$, so the conditional expectation $\EE( f \circ \sigma_{\tau_r} | \F_r )$ only depends on the metric in the random set $D_r$, defined above.
	
	It is not trivial to formalize this intuition.  An important concept will be that of \emph{continuous disintegration} (introduced by LaGatta in \cite{lagatta2010continuous}), which we will carefully investigate in Section \ref{sect_RCP}.  In our context, this is a function $P : \C \times \Omega_+ \to \M$ (written $(C,g) \mapsto P_C(g,\cdot)$)\footnote{Recall that $\C$ is the space of compact sets in $\R^2$ with the Hausdorff topology, and $\M$ is the space of Radon measures on $\Omega$, equipped with the topology of weak convergence of measures.} which satisfies some useful properties (detailed in Theorem \ref{P_lem}). Namely, this function is jointly continuous, and is a regular conditional probability: for each compact set $C \subseteq \R^2$ and almost every $g$, 
		\begin{equation} \label{PCversion}
			\mbox{$P_C(g,\cdot)$ is a version of the conditional probability $\PP(\cdot|\F_C)$.} \end{equation}
	The measure $P_C(g,\cdot)$ is supported on the set of Riemannian metrics which are equal to $g$ on the set $C$:  all the other information $g$ carries is lost.

	\subsection{The Local Markov Property} \label{subsect_LMP}

	For any measurable $f$, we write $P_C(g,f) := \int_\Omega f(g') \, P_C(g,\D g')$ for brevity; \eqref{PCversion} implies that $P_C(g,\D g') = \EE(f|\F_C)$ almost surely.  If $f$ is $\F_C$-measurable, then for almost every $g$, the conditional expectation $P_C(g,f)$ equals $f(g)$.  The (weak) continuity property means that if $f$ is a bounded, continuous function, and $C_n \to C$ and $g_n \to g_n$, then $P_{C_n}(g_n,f) \to P_C(g,f)$.
	
	\begin{env_thm}[Local Markov Property] \label{thm_markov}
		Suppose that $d=2$.  Fix some compact set $D \subseteq \R^2$, and let $f : \Omega \to \R$ be an $\F_{D}$-measurable random variable.  Let $r \ge 0$, and define the random lens-shaped set $D_r(g) = D^1 \cap B_r(g)$ as in \eqref{Drdef}.  Then
			\begin{equation} \label{markovstatement}
				\mbox{$P_{D_r}( \sigma_{\tau_r} g, f)$ is a version of the conditional expectation $\EE(f \circ \sigma_{\tau_r} | \F_r)$.} \end{equation}		
		for $\PP$-almost every $g$ on the event $\{ r < \Rmax \}$.
	\end{env_thm}

	If the POV Theorem is proved for the case $d \ge 3$, then the Local Markov Property extends without difficulty. The first step in the proof is a simple lemma about conditional probabilities.  Lemma \ref{condexp_Bolem} states that if $f$ is $\F_D$-measurable, then for any compact $B$, $\EE(f|\F_B) = \EE(f|\F_{B\cap D^1})$.  This lemma is easy to prove, and relies on the fact that conditional expectations are $L^2$-projections.  The vector $f$ lies in the $D$-subspace, and the left side is the projection of this vector onto the $B$-subspace.  The $\sigma$-algebras $\F_D$ and $\F_{B - D^1}$ are independent, hence their corresponding subspaces are orthogonal.  The rest of the proof is non-trivial, and can be found in Section \ref{sect_proofmarkov}.

	\subsection{Strong Local Markov Property} \label{subsect_SLMP}

%	In this section, we strengthen the Local Markov Property to a Strong Local Markov Property.  This means that we replace the deterministic radius $r$ in Theorem \ref{thm_markov} with a random ``stopping radius'' $R$.  
	
	Let $R = R(g)$ be a non-negative random variable which satisfies the property that for all $r \ge 0$, the event $\{R \le r\} \in \F_r$.  Such an $R$ is a ``stopping radius'' with respect to the filtration $\F_r$.\footnote{We prefer to call such an $R$ a ``stopping radius'' rather than a ``stopping time,'' since we already use the word ``time'' to refer to the parametrization of geodesics.} Define the $\sigma$-algebra $\F_R := \{ A \in \F : \mbox{for all $r \ge 0$, $A \cap \{R \le r\} \in \F_r$} \}$, which represents the information in the ball $B(0,R)$ of (random) radius $R$.  Consider the exit time $\tau_R := \tau_{R(g)}(g)$ from this ball.  On the event $\{R = \oo\}$, we set $\tau_R = \oo$.

	The maximal radius $\Rmax = \sup\{ |\gamma(t)| \}$ is a simple example of a stopping radius, since the event $\{\Rmax \le r\}$ depends only on an infinitesimal neighborhood of the ball $B(0,r)$, hence is $\F_r$-measurable.  If $\Rmax = \oo$, then $\F_{\Rmax} = \F$; if $\Rmax < \oo$, then $\F_{\Rmax}$ is non-trivial, and represents the information given by the metric in the random ball $B(0,\Rmax)$.  The event $\{\mbox{$\gamma$ is bounded}\}$ is $\F_{\Rmax}$-measurable.
	
	\begin{env_lem} \label{RRmax}
		Let $R = R(g)$ be a stopping radius, and consider the exit time $\tau_R$ of the geodesic $\gamma$ from the random ball $B(0,R)$.  The $\Omega$-valued random variable $g \mapsto \sigma_{\tau_{R(g)}(g)}^{(g)}(g) =: \sigma_{\tau_R} g$ and the event $\{R < \Rmax\}$ are both $\F_R$-measurable.
	\end{env_lem}
	
	The proof of this lemma is straightforward.  

	\begin{env_thm}[Strong Local Markov Property] \label{thm_strongmarkov}
		Suppose $d=2$.  Fix some compact set $D \subseteq \R^2$, and let $f : \Omega \to \R$ be an $\F_{D}$-measurable random variable.  Let $R$ be a stopping radius, and define $D_R(g) = D^1 \cap B_R(g)$ as in \eqref{Drdef}.  Then
			\begin{equation} \label{strongmarkovstatement}
				\mbox{$P_{D_R}( \sigma_{\tau_R} g, f)$ is a version of the conditional expectation $\EE(f \circ \sigma_{\tau_R} | \F_R)$.} \end{equation}		
		for $\PP$-almost every $g$ on the event $\{ R < \Rmax \}$.
	\end{env_thm}

	The proof follows easily from the Local Markov Property, and uses a standard approximation argument.  It can be found in Section \ref{sect_proofmarkov}.

\subsection{The Inevitability Theorem} \label{subsect_inevitability}

	We continue with the notation of the previous sections.  Let $D \subseteq \R^2$ be compact, and let $U \in \F_D$ be some event depending only on the metric in set $D$.  Let $U_r = \big\{ g : \sigma_{\tau_r} g \in U \big\}$ be the event that $U$ occurs near the exit location $\gamma(\tau_r)$.\footnote{For a simple example, consider some fixed metric $g_* \in \Omega_+$ and $\epsilon > 0$, and let $U = \{ g : \|g - g_*\|_D < \epsilon \}$ be the event that the random metric $g$ is close to the fixed metric $g_*$ over the domain $D$.  The event $U_r$ is then that the exit time POV metric $\sigma_{\tau_r} g$ is close to $g_*$ near the exit location $\gamma(\tau_r)$.}  Our goal is to find a sufficient condition to guarantee that at least one of the events $U_r$ occurs. 
	
	To do this, we will need some uniform control on the conditional probabilities $p_r(g) := \PP(U_r | \F_r )$.  The Local Markov Property states that $p_r(g) = P_{D_r}(\sigma_{\tau_r} g, U)$ a.s., where $D_r(g)$ is the random lens-shaped set defined in \eqref{Drdef}.  For a typical realization $g$, the function $r \mapsto p_r(g)$ might be very poorly behaved.

	Let $p > 0$.  Suppose that $R_k$ is a sequence of stopping radii satisfying $R_k \ge R_{k-1} + 1$, and suppose that
		\begin{equation} \label{lem_musthappen_estimate}
			P_{D_{R_k}}( \sigma_{\tau_{R_k}}g, U ) \ge p > 0 \end{equation}
	for almost every Riemannian metric $g$ on the event $W := \{ \mbox{$R_k$ is well-defined for all $k$} \}$.\footnote{A priori, condition \eqref{lem_musthappen_estimate} need not be satisfied, of course.  In Section \ref{subsect_unifprobest}, we find a sufficient condition which implies the existence of a sequence of stopping radii $R_k$ satisfying \eqref{lem_musthappen_estimate}.}  %We remark that the theorem of Kim and LaGatta \cite{kim2012bounded} implies that $\PP(W) < 1$.

%	\begin{figure}[h!] \label{fig_manychances}
%			\includegraphics{figure-manychances.pdf}
%		\caption{There is a $p$ chance of the event $U$ occurring at any of the exit locations $\gamma(\tau_{R_k})$.}
%	\end{figure}

	The Inevitability Theorem (Theorem \ref{lem_musthappen}) states that if \eqref{lem_musthappen_estimate} is satisfied, then conditioned on the event $W$, with probability one, one of the events $U_{R_k}$ must occur.  Let $K_0  = 0$, and let 
		\begin{equation} \label{Kj_def}
			K_j = \inf\{ k > K_{j-1} : \mbox{$U_{R_k}$ occurs} \} \end{equation}
	be the $j$th occurrence index.  The strong estimate \eqref{lem_musthappen_estimate} implies that the events $U_{R_k}$ occur infinitely often, and that there is an exponential decay of the distribution of waiting times between occurrences.  The proof of this theorem relies only on the Strong Local Markov Property and some elementary probability.
	
		\begin{env_thm}[The Inevitability Theorem] \label{lem_musthappen}
			Suppose $d=2$.  Let $D \subseteq \R^2$ be compact, and let $U \in \F_D$.  Let $R_k$ be a sequence of stopping radii satisfying $R_k \ge R_{k-1} + 1$ and the uniform probability estimate \eqref{lem_musthappen_estimate}.  Let $W := \{ \mbox{$R_k < \oo$ for all $k$}\}$ be the event that the sequence $R_k$ is well-defined.
			
			For almost every $g$ on the event $W$, the events $U_{R_k}$ occur infinitely often.  Let $p$ denote the same value as in \eqref{lem_musthappen_estimate} and let $K_j$ be the sequence defined in \eqref{Kj_def}.  Then $\PP\big( K_{j+1} - K_j > k \, \big| W_k \big) \le (1-p)^k,$ where $W_k = \{ R_k < \oo \}$ is the event that the $k$th stopping radius $R_k$ is well-defined.
		\end{env_thm}
		\begin{proof}
			We focus on formula \eqref{lem_musthappen_estimate} in the case when $j=0$, and prove that $\PP(K_1 > k | W) \le (1-p)^k$; the general case is similar.  For simplicity of notation, we write $U_k := U_{R_k}$.
			
			Introduce the complementary indicator functions $f_i := 1_{U_i^c}$, and define the product $X_k := \prod_{i=1}^k f_i$.  We must prove that $\EE\big( X_k \big| \F_{R_k} \big) \le (1-p)^k$ a.s..  Applying the Strong Local Markov Property to the function $f(g) = 1_U(g)$, we have
				\begin{equation} \label{lem_musthappen_proof1}
					\EE\big( f_i  \big| \F_{R_i} \big) = 1 - \PP\big( \sigma_{\tau_{R_i}}^{-1} U \big| \F_{R_i} \big) = 1 - \EE\big( 1_{U} \circ \sigma_{\tau_{R_i}} \big| \F_{R_i} \big) = 1 - P_{D_{R_i}}\big( \sigma_{\tau_{R_i}}g, U\big) \le 1 - p \end{equation}
			by the uniform probability estimate \eqref{lem_musthappen_estimate}.
			
			The random variable $X_k$ depends on the metric in the $1$-neighborhood of the ball $B(0,R_k)$.  By construction, $R_{k+1} \ge R_k + 1$, so $X_k$ is $\F_{R_{k+1}}$-measurable.  This and estimate \eqref{lem_musthappen_proof1} imply that
				\begin{equation} \label{lem_musthappen_proof2}
					\EE\big( X_{k+1} \big| \F_{R_{k+1}} \big) = X_k \cdot \EE\big( f_{k+1} \big| \F_{R_{k+1}} \big) \le X_k \cdot (1-p). \end{equation}
			Take the conditional expectation of both sides with respect to the $\sigma$-algebra $\F_{R_k}$.  By the tower property of conditional expectations and formula \eqref{lem_musthappen_proof2},
				$$\EE\big( X_{k+1} \big| \F_{R_{k+1}} \big) = \EE\big[ \EE\big( X_{k+1} \big| \F_{R_{k+1}} \big) \big| \F_{R_k} \big] \le \EE\big( X_{k} \big| \F_{R_{k}} \big) \cdot (1-p).$$
			It follows by induction that $\EE\big( X_{k} \big| \F_{R_{k}} \big) \le (1-p)^k$, which proves the estimate
				\begin{equation} \label{iterated}
					\PP(U^c_1 \cap \cdots \cap U^c_k | \F_{R_k} ) \le (1-p)^k. \end{equation}
					
			Let $W_k = \{R_k < \oo\}$ be the event that the $k$th term of the sequence is well-defined.  Clearly, $W_k$ is $\F_{R_k}$-measurable, and $W = \bigcap W_k$.  We now estimate
				\begin{eqnarray*}
					\PP(K_1 > k | W_k) &=& \tfrac{1}{\PP(W_k)} \PP( \mbox{$K_1 > k$ and $W_k$}) = \tfrac{1}{\PP(W_k)} \EE\big[ \PP( \mbox{$K_1 > k$ and $W_k$} | \F_{R_k} ) \big] \\
					&=& \tfrac{1}{\PP(W_k)} \EE\big[ 1_{W_k} \PP( \mbox{$K_1 > k$} | \F_{R_k} ) \big] \le  \tfrac{1}{\PP(W_k)} \EE\big[ 1_{W_k} (1-p)^k \big] = (1-p)^k 
				\end{eqnarray*}
%			since $\PP(W_k) \to \PP(W)$ by the monotone convergence theorem, and $p > 0$.  This proves that on the event $W$, at least one of the events $U_k$ occurs.  
			
			It is easy to adapt this argument to show that $U_k$ occurs infinitely often.  It is also easy to show that the waiting time between events $U_k$ decays exponentially at the rate at least $-\log(1-p)$. 
	\end{proof}
	
	We will use this theorem to prove the Main Theorem.  In Theorem \ref{prop_unifprobest}, we provide a sufficient condition to guarantee the uniform probability estimate \eqref{lem_musthappen_estimate}.

	\section{Conditioning the Metric and Uniform Probability Estimate} \label{sect_RCP}

	To prove the main theorem, we will need some technical estimates on conditioning Gaussian measures, as well as on fluctuations of random metrics. In Section \ref{subsect_contdis}, we show that the measure $\PP$ exhibits the continuous disintegration property, as introduced by LaGatta in \cite{lagatta2010continuous}. In particular, we show the existence of a measure-valued function $(g,D) \mapsto P_D(g,\cdot)$ which represents the conditional probability $\PP(\cdot|\F_D)$, and we show that this function varies (weakly) continuously with respect to both the metric $g$ and the set $D$. 
	
	In Section \ref{subsect_fluctuations}, we state and prove a simple estimate on the fluctuations of the metric. In Section \ref{subsect_unifprobest}, we use the fluctuation estimate to prove a uniform probability estimate for the law $\PP$.
	
	\subsection{Continuous Disintegrations of $\PP$} \label{subsect_contdis}
	
	By our construction, $\PP$ is the push-forward of a stationary Gaussian measure $\Q$ on the Fr\'echet space $\Omega = C(\R^d, \Sym)$.  The central result in \cite{lagatta2010continuous} is that stationary Gaussian measures on Banach spaces always admit continuous disintegrations, meaning that the Gaussian measure depends (weakly) continuously on the conditioning parameters. In this article, we extend those arguments to the setting of Gaussian tensor fields on $\R^d$, which easily implies that $\PP$ satisfies the continuous disintegration property. We further generalize \cite{lagatta2010continuous} by allowing a continuous dependence on the set $D$.\footnote{Our arguments easily generalize to the setting of Gaussian measures on Fr\'echet spaces.} 
	
	An important property of $\PP$ is strict positivity:  $\PP(U) > 0$ for all open events $U$. The conditional probability $P_D(g,\cdot)$ inherits the strict positivity condition from $\PP$, but with one important caveat: an open set $U$ has positive $P_D(g,\cdot)$-probability only when $U \cap [g]_D \neq \varnothing$, where $[g]_D := \big\{ g' \in \Omega_+ : g'(x) = g(x) \mathrm{~for~all~} x \in D \big\}$ is the equivalence class of metrics which agree with $g$ on the domain $D$.

	\begin{env_thm} \label{P_lem}
		There exists a measure-valued function $P : \C \times \Omega_+ \to \M$, which we write as $(D,g) \mapsto P_D(g,\cdot)$, and which satisfies the following properties:
		 	\begin{enumerate}[a)]
		 		\item \label{P_lem_condprob} (Conditional probability) For any compact $D \subseteq \R^d$ and for $\PP$-almost every $g$, the measure $P_D(g,\cdot)$ is a version of the conditional probability measure $\PP(\cdot | \F_{D} )$.
		 		\item \label{P_lem_fiber} (Support)  For any compact $D \subseteq \R^d$ and every $g \in \Omega_+$, the measure $P_D(g, \cdot)$ is supported on the equivalence class $[g]_D$.\footnote{This implies that $P_D(g, [g]_D) = 1$, so for $P_D(g,\cdot)$-almost every $g'$, $g'(x) = g(x)$ for all $x \in D$.  Furthermore, for any $g' \in [g]_D$, the measures $P_D(g,\cdot)$ and $P_D(g',\cdot)$ are equal.}
		 		\item \label{P_lem_totallypositive}  (Conditional strict positivity)  If $U \in \F$ is an open event which meets the equivalence class $[g]_D$, then $P_D(g,U) > 0$.
		 		\item \label{P_lem_convergence}  (Weak relative compactness) If $D_n \to D$ in the Hausdorff topology and if $g_n \to g$ in $\Omega_+$, then the measures $P_{D_n}(g_n, \cdot)$ converge weakly to $P_D(g, \cdot)$.
		 	\end{enumerate}
	\end{env_thm}
			
		The proof of this theorem is given in Section \ref{sect_proofctsdisint}.
		
%	Recall that the mean-zero Gaussian measure $\Q$ is completely characterized by its covariance operator $K : \Omega^* \to \Omega$.  The conditional measure $Q_D(\xi,\cdot)$ is also a Gaussian measure on $\Omega$, so it is determined by its mean tensor field $m_D(\xi) \in \Omega$ and covariance operator $K_D : \Omega^* \to \Omega$.  
	
%	The most difficult part of the construction is defining the family of \emph{conditional mean operators} $m_D : \Omega \to \Omega$.  This is done in Lemma \ref{mo_lem}; the proof of this lemma is Appendix \ref{proof_mo_lem}.  The definition of the \emph{conditional covariance operators} easily follows:  $K_D = K - K m_D^*$.  Lemma \ref{Ko_lem} states the basic properties of the operators $K_D$.

	\subsection{Fluctuations of Random Riemannian Metrics} \label{subsect_fluctuations}

	We next introduce random variables $Z_D(g)$ which quantify the fluctuations of the random metric $g$.  Each Riemannian metric $g \in \Omega_+$ is strictly positive-definite, so the inverse metric $g^{-1} \in \Omega_+$ (defined by the pointwise matrix inverse) is also a well-defined random Riemannian metric. For each compact set $D \subseteq \R^d$, we define the non-linear functional $Z_D : \Omega_+ \to \R$ by
		\begin{equation} \label{Zdef}
			Z_D(g) = \max\big\{ \|g - \delta\|_{C^{2,1}(D)}, \|g^{-1} - \delta\|_{C^{1,1}(D)} \big\}, \end{equation}
	where $\delta \in \Omega_+$ denotes the flat Euclidean metric on $\R^d$.\footnote{The seminorm $\|\xi\|_{C^{\alpha,1}(D)}$ measures the local Lipschitz fluctuations of the $\alpha$th derivatives of a tensor field $\xi_{ij}(x)$ on an infinitesimal neighborhood of $D$.  We define this precisely in Appendix \ref{analytictools}.}  Clearly, $Z_D(g) = 0$ if and only if $g = \delta$.  

	The quantity $Z_D(g)$ measures the deviation of the metric $g$ from the Euclidean metric on the region $D$.  The Christoffel symbols and scalar curvature are defined pointwise in terms of the metric, its inverse, and their derivatives (cf. \eqref{geoquantitiesdef}).  Consequently, a uniform bound on $Z_D$ translates into Lipschitz bounds on those important geometric quantities.  Curvature is defined using the first two derivatives of the metric and only the first derivatives of its inverse; this is why in the definition of $Z_D$ we take only the $C^{1,1}$ norm of the inverse fluctuations. %We will use this extensively in Part II in our construction of a geometric bump surface.
	
	It is easy to see that the real-valued function $(D,g) \mapsto Z_D(g)$ is jointly continuous (cf. Lemma \ref{ZD_cty}).  That is, if $D_n \to D$ in the Hausdorff topology and $g_n \to g$ in $\Omega_+$, then $Z_{D_n}(g_n) \to Z_D(g)$ in $\R$.
	
	\begin{env_pro}[Estimate on Local Fluctuations] \label{ZD_t}
		Set $m = 3^d$ and fix a compact set $D$.  Let $Z_1, \cdots, Z_m$ be $m$ independent copies of the random variable $Z_D$.  Then
			\begin{equation}
				\EE \min\{ Z_1, \cdots, Z_m \}^{2m+1} < \oo. \end{equation}
	\end{env_pro}
	\begin{proof}
		Fix the domain $D \subseteq \R^d$, and define the random variable $X(\xi) := \max\{ \|\xi\|_{C^{2,1}(D)} \}$ on $\Omega$.  Since $\xi$ is a Gaussian random field, the random variable $X$ satisfies the Gaussian large-deviations estimate $\Q( X > u ) \le \E^{-u^2/2c^2}$ for some constant $c > 0$.\footnote{This follows from a general form of the Borell-TIS inequality \cite{adler07} applied to the Banach space $C^2(D,\Sym)$ of smooth functions.}
	
		Set $Y = \min\{ Z_1, \cdots, Z_m \}^{2m+1}$.  By the pointwise definition of the random metric $g = \varphi \circ \xi$, it is clear that $Y$ is large if and only if $X$ is large.  Since $X$ satisfies a large-deviations estimate, and since the growth of the real-valued functions $|\varphi(u)|_{C^{2,1}}$ and $1/|\varphi(-u)|_{C^{1,1}}$ is dominated above and below by polynomials as $u \to \oo$, this implies the moment estimate $\EE Y < \oo$.
	\end{proof}

	\subsection{Uniform Probability Estimate} \label{subsect_unifprobest}

	Let $U \in \F$ be an open event.  Our goal is to construct a relation $\ZZ \subseteq \C \times \Omega_+$ between compact sets and Riemannian metrics on which the lower bound $\inf_\ZZ P_D(g,U)$ is positive. To construct the set $\ZZ$, we use the fluctuation observables $Z_D(g)$, defined in \eqref{Zdef}, which measure how much $g$ deviates from the Euclidean metric on the set $D$.  %Lemma \ref{ZD_cty} states that the function $(D,g) \mapsto Z_D(g)$ is jointly continuous.  Consequently, any set of the form $\{ Z_D(g) \le h \}$ is closed, a necessary property of compactness.
	
	Fix some compact family $\DD \subseteq \C$ (i.e., the family $\DD$ is compact with respect to the Hausdorff metric).  Fix a number $h \ge 0$, and consider the family
		\begin{equation} \label{ZZdef}
			\ZZ := \ZZ(\DD,h) := \{ (D, g) \in \DD \times \Omega_+ : Z_D(g) \le h \} \end{equation}
	of pairs $(D,g)$ which satisfy the estimate $Z_D(g) \le h$.  This estimate implies that $\| g^{-1} \|_{C(D)} \le (1 + h)$, hence the minimum eigenvalue of metrics on $D$ is uniformly bounded below by $\tfrac{1}{1+h} > 0$.  By Lemma \ref{ZD_cty}, the function $(D,g) \mapsto Z_D(g)$ is continuous.  Consequently, the family $\ZZ$ is a closed subset of $\C \times \Omega_+$. Let 
		\begin{equation} \label{barD_def}
			\bar D := \overline{\bigcup_{D \subseteq \DD} D^1} \subseteq \R^d \end{equation}
	be the \emph{compact cover} of the family $\DD$, where $D^1$ denotes the Euclidean $1$-neighborhood of the compact set $D$.  It is easy to see that the closed set $\bar D$ is bounded, hence a compact subset of $\R^d$.
	
%	If $\Omega_+$ were the open cone $C^2(\bar D, \SPD)$ in the Banach space $C^2(\bar D, \Sym)$, then the Arzel\`a-Ascoli theorem would imply that the family $\ZZ$ is compact.  Instead, $\Omega_+ = C^2(\R^d, \SPD)$ is a subset of the Fr\'echet space $C^2(\R^d, \Sym)$, so the Arzel\`a-Ascoli theorem does not immediately apply.  
	
	For each set $D \in \DD$, we let $\psi_D : \Omega_+ \to \Omega_+$ be a suitable Urysohn operator, precisely defined in \eqref{psiD_def}.  The modification $\psi_D(g)$ preserves the metric on the set $D$, but smoothly interpolates so that at Euclidean distance $1$ away from $D$, the metric is the flat Euclidean metric.  We use the operators $\psi_D$ to define the \emph{compact core}
		\begin{equation} \label{ZZcore_def}
			\ZZcore = \{ (D, \psi_D(g) ) \in \DD \times \Omega_+ : Z_D(g) \le h \}. \end{equation}
	of $\ZZ$.  We can think of $\ZZcore$ as a complete family in the metric space $\DD \times C^2_\delta(\bar D, \SPD)$, where the function space satisfies a flat boundary condition.  The Arzel\`a-Ascoli theorem then implies that the family $\ZZcore$ is compact (this is Sublemma \ref{ZZcore_compact}).
	
	The mollified metric $\psi_D(g)$ is a member of the equivalence class $[g]_D$, so the conditional probabilities $P_D(g,U)$ and $P_D(\psi_D(g), U)$ are equal.  This implies equality of the minimum probabilities 
		\begin{equation} \label{minprobsame}
			\inf_{\ZZ} P_D(g,U) = \inf_{\ZZcore} P_D(g,U). \end{equation}
	Of course, there is still one issue to attend to:  a priori, there is no guarantee that the minimum probability $\inf_{\ZZcore} P_D(g,U)$ is positive.  This will follow from the conditional strict positivity condition, part \eqref{P_lem_totallypositive} of Theorem \ref{P_lem}.
	
	To guarantee conditional strict positivity, we must assume the joint condition \eqref{uniformassumption} on the event $U$ and the family $\ZZ$:  for every pair $(D,g) \in \ZZ$, the open set $U$ must meet the equivalence class $[g]_D$.  Conditional strict positivity then implies that $P_D(g,U) > 0$ for all $(D,g) \in \ZZ$.  This implies the same result for $(D,g) \in \ZZcore$, hence the infimum $\inf_{\ZZcore} P_D(g,U)$ is positive.
	
%	Uses continuous disintegration property.  Essentially we are taking a minimum over a compact set.
			
	\begin{env_thm}[Uniform Probability Estimate] \label{prop_unifprobest}
		Fix some compact family $\DD$ of compact sets and some $h \ne 0$, and define the family $\ZZ$ as in \eqref{ZZdef}.  %Define
%			$$\ZZ = \{ (D, g) \in \DD \times \Omega_+ : Z_D(g) \le h \}.$$
		Let $U \in \F$ be an open subset of $\Omega$ with the property that 
			\begin{equation} \label{uniformassumption}
				\mbox{for all pairs $(D, g) \in \ZZ$, the open set $U$ meets the equivalence class $[g]_D$.} \end{equation}
		Then the probabilities $P_D(g,U)$ are uniformly bounded below by $\inf_{\ZZ} P_D(g, U) > 0$.
	\end{env_thm}
		\begin{proof}		
		For each compact set $D \in \DD$, define some Urysohn operator $\psi_D : \Omega_+ \to \Omega_+$ satisfying
			\begin{equation} \label{psiD_def}
				\psi_D(g)(x) = \begin{cases} g(x), & x \in D \\ \delta, & x\notin D^1, \end{cases},  \end{equation}
		where $D^1$ denotes the $1$-neighborhood of the set $D$, and $\delta$ is the Euclidean metric.  Suppose furthermore that the operators $\psi_D$ satisfy the uniform condition
			\begin{equation} \label{psi_smooth}
				\| \psi_D(g) - \delta \|_{C^{2,1}(\R^d)} \le C \| g - \delta \|_{C^{2,1}(D)} \end{equation}
		for some universal constant $C$. The Urysohn operator $\psi_D$ fixes a metric on the set $D$, flattens it to the Euclidean metric off the larger set $D^1$, and smoothly interpolates in between.  Since $\psi_D(g)$ agrees with $g$ on the set $D$, the metric $\psi_D(g)$ belongs to the equivalence class $[g]_D$.  By property \eqref{P_lem_fiber}, the probabilities $P_D(g,U)$ and $P_D(\psi_D(g),U)$ are equal.  
		
		Define the compact core $\ZZcore$ as in definition \eqref{ZZcore_def}.  The previous paragraph implies that the minimum probabilities $\inf_\ZZ P_D(g,U)$ and $\inf_{\ZZcore} P_D(g,U)$ are equal, which proves statement \eqref{minprobsame}. Let $\bar D := \overline{\bigcup_{\DD} D^1}$ be the compact cover of $\DD$, defined in \eqref{barD_def}.  Since $\DD$ is a compact family of compact sets, it is a simple exercise to show that $\bar D$ is a compact subset of $\R^d$.  The compact core $\ZZcore$ is a complete set in the metric space $\DD \times C^2_\delta(\bar D, \SPD)$, where the function space is with a flat boundary condition.
		
		Let $\G = \{ g \in \Omega_+ : \mbox{$Z_{\bar D}(g) \le Ch$ and $g|_{\partial \bar D} = \delta$} \}$ be the set of metrics which satisfy a certain fluctuation estimate on the compact cover (for $C$ the constant in \eqref{psi_smooth}), and a flat boundary condition.
		
		\begin{env_sublem} \label{ZZcore_compact}
			The set $\G$ is compact.  Consequently, the compact core $\ZZcore$ is compact.
		\end{env_sublem}
		\begin{proof}
			Let $\Omega_{\bar D} = \{ g \in C^2(\bar D, \Sym) : g|_{\partial \bar D} = \delta \}$ be the Banach space of $C^2$-smooth, symmetric quadratic forms on $\bar D$ which satisfy a flat boundary condition.  The estimate $Z_{\bar D}(g) \le h$ implies that the second derivatives of $g$ satisfy a uniform Lipschitz condition.  The Arzel\`a-Ascoli theorem implies that the set $\G$ is precompact in the Banach space $\Omega_{\bar D}$.
			
			The estimate $Z_{\bar D}(g) \le h$ also implies that the minimum eigenvalue of the metric $g$ is uniformly bounded below on $\bar D$ by the constant $\tfrac{1}{1+Ch} > 0$, which means that the set $\G$ is a complete set within the open cone $\Omega_+$.  Since the set $\G$ is precompact and closed, it is compact.
			
			The set $\ZZcore$ is a closed subset of the compact product space $\DD \times \G$, hence compact.
		\end{proof}
		
		By assumption \eqref{uniformassumption}, the set $U$ meets the equivalence class $[g]_D = [\psi_D(g)]_D$ for all $(D,g) \in \ZZ$.  Consequently, $P_D(g,U) = P_D(\psi_D(g),U) > 0$ for all $(D,g) \in \ZZ$ by conditional strict positivity.  The minimum probabilities are equal (by \eqref{minprobsame}), so
			$$\inf_\ZZ P_D(g,U) = \inf_\ZZ P_D( \psi_D(g), U ) = \inf_{\ZZcore} P_D(g, U) > 0.$$
		This complete the proof of Theorem \ref{prop_unifprobest}.		
		\end{proof}
		
	\begin{env_rem}
		Theorem \ref{prop_unifprobest} remains true with the slightly weaker definition $Z_D'(g) = \max\{ \|g\|_{C^{2,1}(D)}, \|g^{-1}\|_{C(D)} \}$ for metric fluctuations.
	\end{env_rem}

	In what follows, we use Theorem \ref{prop_unifprobest} to ensure that the hypotheses of the Inevitability Theorem (Theorem \ref{lem_musthappen}) are satisfied, and we use this to show that $\gamma$ is not minimizing almost surely.

	\section{Frontier Radii} \label{sect_frontiertimes}
	
	In this section, we introduce the notion of a ``frontier radius'': a stopping radius which satisfies additional uniformity properties.  Pick a starting direction $v \in S^{d-1}$, and consider $\gamma_v := \gamma_{0,v}$, the unit-speed geodesic starting at the origin in direction $v$.  The geodesic may be either bounded (so that $|\gamma_v| \le \Rmax$ for some $\Rmax(v,g)$) or unbounded.  
	
	If $\gamma_v$ is unbounded, it will exit arbitrarily large balls.  Let $\tau_{v,r}$ be the exit time of $\gamma_v$ from the ball $B(0,r)$, and let $\sigma_{v,r} g$ denote the environment from the point of view of the exit location $\gamma(\tau_{v,r})$; these quantities are defined in Section \ref{sect_exittime}.  The environment $\sigma_{v,r} g$ is a random Riemannian metric with a complicated law.\footnote{In the case of $d=2$ and deterministic starting direction $v$, a modification of Theorem \ref{POVthm_taur} implies that the law of $\sigma_{\tau_{v,r}} g$ is absolutely continuous with respect to $\PP$, and we give an expression for its Radon-Nikodym derivative.}  It could be the case that as $r \to \oo$, the law of $\sigma_{v,r} g$ concentrates on degenerate or singular metrics.
	
%	For arbitrary Riemannian metrics $g \in \Omega_+$, we showed in Theorem \ref{transientgeodesics} that all minimizing geodesics are transient.  That theorem relied only on a simple geometric argument.  By applying arguments from lattice first-passage percolation to our model, we are able to prove more quantitative statements on the behavior of minimizing geodesics.
	
	\subsection{The Frontier Theorem}  \label{subsect_frontierthm}
	
	In Theorem \ref{frontier_thm}, we show that when $\gamma_v$ is a \emph{minimizing} geodesic (i.e., $v \in \V_g$), the environment as seen along the geodesic is especially well behaved.  In particular, we show that (with probability one) for every $v \in \V_g$, we can find a well-defined sequence of increasing frontier radii $R_k$ such that the fluctuations of the metric (as expressed by the quantity $Z_{D_{v,R_k}}(\sigma_{\tau_{v,R_k}} g)$) are uniformly bounded in $k$. Simultaneously, we prove that the geodesic $\gamma_v$ does not exit the balls $B(0,R_k)$ in a degenerate manner: the (Euclidean) exit angles are uniformly bounded.
	
	To state this theorem precisely, we must introduce some notation.  Let $o_{v,r} := \gamma(\sigma_{\tau_{v,r}}g, -\tau_{v,r})$ denote the ``old origin'' from the point of view of the exit location $\gamma_v(\tau_{v,r})$.  The POV transformation is defined by (random) isometries of $\R^d$, and the old origin $o_{v,r}$ is the image of the origin after these transformations.  Consequently, the (random) ball $B(o_{v,r}, r)$ is of principal importance.  
	
	Define the lens-shaped sets $D_{v,r} := B(0,2) \cap B(o_{v,r}, r)$.\footnote{For an illustration of the old origin $o_{v,r}$ and the lens-shaped set $D_{v,r}$ in the case that $v = \E_1$, consult Figure \ref{figure-DDo}.} Recall that $Z_{D_{v,r}}(\sigma_{\tau_{v,r}} g)$ measures the fluctuations of the POV metric $\sigma_{\tau_{v,r}} g$ on the set $D_{v,r}$. Let $\alpha_{v,r} \in [0, \tfrac \pi 2]$ denote the (Euclidean) exit angle of $\gamma$ from $B(0,r)$:
		\begin{equation} \label{alphadef}
			\cos \alpha_{v,r} := \frac{\langle \gamma_v(\tau_{v,r}), \dot\gamma_v(\tau_{v,r})\rangle}{r|\dot\gamma_v(\tau_{v,r})|}. \end{equation}
	i.e., $\alpha_{v,r}$ equals the angle between the vectors $\gamma_v(\tau_{v,r})$ and $\dot\gamma_v(\tau_{v,r})$.  The geodesic exits the ball tangentially when $\alpha_{v,r} = \tfrac \pi 2$, and its exit vector is normal to the ball when $\alpha_{v,r} = 0$.
	
	The heuristic content of Theorem \ref{frontier_thm} is that there exist uniform constants $h > 0$ and $\theta < \tfrac{\pi}{2}$ such that, with probability one, for all $v \in \V_g$, there exists a sequence $R_k \uparrow \oo$ of frontier radii with
		\begin{equation} \label{frontier_radii_statement}
			\mbox{$\alpha_{v,R_k} \le \theta$ and $Z_{D_{v,R_k}}(\sigma_{\tau_{v,R_k}} g) \le h$.} \end{equation}

	There is of course an issue of measurability, as the random variables $R_k(v,g)$ are themselves defined on the random set $\V_g$.  In this section, we circumvent this difficulty by instead focusing on certain random sets $Q_v \subseteq \R$ (defined in \eqref{def_Qv}).  In Theorem \ref{frontier_thm}, we prove that these sets have uniformly positive (lower) Lebesgue density.  In Section \ref{subsect_repeatedevents}, we focus on the case $v = \E_1$, condition on the event $\{\E_1 \in \V_g\}$, and define the sequence of random variables $R_k(g)$ using $Q_{\E_1}$.
	
	For any parameter choices $\theta$ and $h$, and any metric $g \in \Omega_+$, we define the sets of ``good'' frontier radii
		\begin{equation} \label{def_Qv}
			Q_v := Q_v(\theta, h, g) := \big\{ r \ge 0 : \mbox{$\alpha_{v,r} \le \theta$ and $Z_{D_{v,r}}(\sigma_{\tau_{v,r}} g) \le h$.} \big\}. \end{equation}
	A priori, the sets $Q_v$ may be empty or sparse.  The next theorem demonstrates that for suitable parameter choices $\theta$ and $h$, this is not the case.  Instead, the sets $Q_v$ have uniformly positive Lebesgue density in all directions $v$.

	\begin{env_thm}[Frontier Theorem] \label{frontier_thm}
		There exist non-random constants $\theta \in [0, \tfrac \pi 2)$, $h > 0$ and $\delta > 0$ such that, for $\PP$-almost every $g$ and for every minimizing direction $v \in \V_g$, the (random) sets $Q_v = Q_v(\theta, h, g)$ have positive Lebesgue density bounded below by $\delta$.  
		
		More precisely, there exists a value $r_0$ (independent of $v$) such that if $r \ge r_0$, then $\Leb(Q_v \cap [0,r]) \ge \delta r$ for all $v$.
	\end{env_thm}
	
	This theorem is the only place in this paper where we use methods from first-passage percolation.  The proof is non-trivial, and can be found in Section \ref{sect_proofoffrontiertimes}.  We critically use properties of minimizing geodesics in the proof.  It would be very interesting if one could show that there is a similar estimate along arbitrary unbounded geodesics. 
	
	In the proof of the Main Theorem:  $\E_1 \notin \V_g$ with probability one, we assume otherwise, and construct a sequence of frontier radii $R_k \uparrow \oo$ satisfying the estimates \eqref{frontier_radii_statement}.  We will see later that the existence of such a sequence will imply that $\gamma_{\E_1}$ is not minimizing.  
	
	Let $\theta$ and $h$ be as in the Frontier Theorem.  Define $R_0 = 0$, and
		\begin{equation} \label{Rkdef}
			R_k := \inf\!\big( Q_{\E_1} \cap [R_{k-1}+1, \oo) \big), \end{equation}
%			R_k = \inf\{ r \ge R_{k-1} + 1 : r \in Q_{\E_1} \}. \end{equation}
	setting $R_k = \oo$ if the set on the right side is empty.  By this construction, $R_k \ge k$.  Theorem \ref{frontier_thm} implies that on the event $\{\E_1 \in \V_g \}$, the sequence $R_k$ is well-defined.  By construction, it is easy to verify that each $R_k$ is a genuine stopping radius, i.e., the event $\{ R_k \ge r \} \in \F_r$ for each $r \ge 0$.  
	
	\begin{env_cor} \label{Rk_welldefined}
		For $\PP$-almost every $g$ on the event $\{\E_1 \in \V_g\}$, the sequence of frontier radii $R_k = R_k(g)$ is well-defined.  Writing $C = \tfrac 1 \delta + 1$, we have $k \le R_k \le C k$ for all but finitely many $k$.
	\end{env_cor}
	\begin{proof}
		If $R_k > Ck$, then $\Leb( Q_v \cap [0,Ck] ) \le k$ (otherwise, we could define some $R_{k+1}$ before $Ck$).  However, Theorem \ref{frontier_thm} implies that $\Leb( Q_v \cap [0,Ck] ) \ge \delta Ck$ for large $k$.  Consequently, $1 \ge \delta C = 1 + \delta$, a contradiction. 
	\end{proof}

	While the Corollary will be instrumental in our proof of the Main Theorem, \emph{ex post} it involves conditioning on the measure-zero event $\{\E_1 \in \V_g\}$, hence is logically vacuous.

	\subsection{Repeated Events along a Minimizing Geodesic} \label{subsect_repeatedevents}

	Henceforth, we suppress the subscript $\E_1$ from our notation.  Let $U \in \F_{B(0,1)}$ be an open event depending only on the metric locally near the origin (an example might be the event that the scalar curvature of the metric in the ball $B(0,1)$ is strictly positive).  Let $R_k$ be the sequence of random variables given by Corollary \ref{Rk_welldefined}, and let $U_k$ be the event that the local event $U$ occurs near the point $\gamma(\tau_{R_k})$.  Precisely, $U_k := \{ g : \sigma_{\tau_{R_k}} g \in U \} = (\sigma_{\tau_{R_k}})^{-1} U$. 
	
	Since the events $U_k$ are local, when we condition on the $\sigma$-algebra $\F_{R_k}$, the event $U_k$ should only depend on the part of the random ball $B(o_{R_k}, R_k)$ near the origin of the POV coordinate chart.  That is, the event $U_k$ only depends on the metric on the set $D_{R_k}$, which by definition satisfies the uniform bound \eqref{frontier_radii_statement}.  We then apply Theorem \ref{prop_unifprobest}, which guarantees that the events $U_k$ have a uniform probability $p$ of occurring. 
	
	We next apply the Inevitability Theorem (Theorem \ref{lem_musthappen}), which states that if this uniform probability estimate is satisfied, then the sequence $U_k$ must occur infinitely often.  This theorem also implies that the first occurrence time $K$ is a random variable with exponential tail decay.
	
	\begin{env_pro} \label{pro_infoften}
		Suppose that $d=2$.  Let $W$ be the event that the sequence $R_k$ is well-defined and satisfies the estimate \eqref{frontier_radii_statement} for $v = \E_1$.  Let $U \in \F_{B(0,1)}$ be an open event, and define the events $U_k$ as above.  The events $U_k$ occur infinitely often on the event $W$. Next, let $K = \inf\{ k \ge 0: \mbox{$U_k$ occurs} \}$ be the first occurrence time.  The random variable $K$ has exponential tail decay on the event $W$: $\PP(K > k |W) \le (1-p)^k$.
	\end{env_pro}

	\section{Bump Surface} \label{sect_bump_short}

	In this section, we construct a particular local event $U$ so that if any of the events $U_k$ occur, then the geodesic $\gamma_{\E_1}$ is not minimizing.  Our method involves the construction of a ``bump metric''.  Throughout this section, we assume that a metric $g$ satisfies the estimate $Z_0(g) \le 2h$ at the origin.\footnote{The condition $Z_0(\sigma_{\tau_{R_k}} g) \le h \le 2h$ is guaranteed by the Frontier Theorem.}  Since this is an estimate on the second derivatives (and inverse) of the metric, it implies that a uniform estimate on the scalar curvature at the origin: $|K_0(g)| \le \Kmax$ for some $\Kmax > 0$.  The estimate also gives us a certain length scale $\tau$ for the bump metric.
	
	For every $g \in \Omega_+$ satisfying the estimate $Z_0(g) \le 2h$, we will construct a bump metric $b(g) \in \Omega_+$.  The geodesic starts tracing out the bump surface at the origin, where the curvature equals $K_0(g)$.  As it follows along the bump metric, the curvature continuously transitions to some value $K_+ := \tfrac{4\pi^2}{\tau^2}$ at time $\tfrac \tau 4$.  At this point, the bump metric has constant curvature $K_+$, hence is locally isometric to the sphere with radius $\tfrac{1}{\sqrt{K_+}}$.  At time $\tfrac{\tau}{2}$, the geodesic reaches the antipodal point on the bump, developing a conjugate point to the origin. By Jacobi's Theorem \cite{lee1997rmi}, a geodesic with conjugate points is not minimizing.\footnote{It is easy to see that minimizing geodesics cannot self-intersect (this follows from the argument of Theorem \ref{minimizing_geodesics}).  Consequently,  an alternative proof of the Main Theorem could rely on an event $U'$, manipulating the geodesic $\gamma_{\E_1}$ to self-intersect near the origin.  The event $U_k'$ would then imply that $\gamma_{\E_1}$ self-intersects shortly after time $\tau_{R_k}$.  This is an interesting strategy. We instead opted for the bump metric construction in order to highlight the geometric role of curvature and its fluctuations.} We construct the bump metric in the case $d=2$ (hence the bump metric represents a bump surface), but the construction in higher dimensions is similar.
	
	\begin{figure}[h!] \label{figurebumpsurface}
			\includegraphics{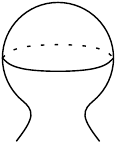}
			\caption{A sketch of a bump surface where $K_0(g)$ is negative. The curvature at the top of the bump is constant and equal to $K_+$, and smoothly transitions to match $K_0(g)$ at the bottom.}
		\end{figure}
	
%	The famous Cartan-Hadamard theorem \cite{ballmann1995lectures} states that for a simply-connected manifold with non-positive (Alexandrov) curvature, there is exactly one geodesic connecting any two points, and all these geodesics are minimizing.  Consequently, for smooth metrics, the presence of positive curvature is a necessary condition for geodesics to lose the minimization property.
	
	To realize the construction of the bump metric, we use Fermi Normal Coordinates, which are a coordinate system adapted along a geodesic.  These coordinates have a canonical form \eqref{gexpansion} which depends only on the curvature of the metric.  Consequently, it is easy for us to define a bump metric with a particular curvature profile. We then convert the Fermi coordinate system back to our original coordinate system.  We show that if we take a sufficiently small perturbation of such a bump metric, the corresponding geodesic is still not minimizing.  Each $g$ gives rise to a bump metric $b(g)$, so we define the open event $U = \{ g : \|g - b(g) \|_{B(0,1)} < \epsilon \}$ for a suitable $\epsilon$.

	\subsection{The Hinterland and Frontier Cones} \label{subsect_hfcones}
	
	We will be describing the construction of the bump surface in a coordinate system centered at the origin.  The reader should think of this as a POV coordinate system, as eventually we plan to show that there is a positive probability of a bump surface near each frontier exit point $\gamma(\tau_{R_k})$.  
	
	As described in Section \ref{subsect_frontierthm}, there are certain uniformity properties which the frontier radii $R_k$ satisfy.  One is a uniformity condition on the metric, which we will return to in Section \ref{subsect_bump}.  The other property is that the geodesic $\gamma$ exits the ball $B(0,R_k)$ at an angle no greater than a fixed constant $\theta < \tfrac{\pi}{2}$.\footnote{The precise statement is that $\alpha_{R_k} \le \theta$, where $\alpha_r := \alpha_{\E_1, r}$ is defined by \eqref{alphadef}.}  
	
	The POV transformation is defined by (random) rigid translations and rotations of the plane.  When we take the POV transformation, the geodesic is sitting at the origin pointing in the horizontal direction.  Consequently, the uniform exit angle translates into a uniform condition on the old origin $o_{R_k}$.  Precisely, (for a.e. $g$ on $\{ \E_1 \in \V_g\}$) the old origin $o_{R_k}$ lies in the \emph{hinterland cone}
		\begin{equation} \label{HCdef}
			HC = \left\{ (y^1, y^2) \in \R^2 : \mbox{$y^1 \le 0$ and $|y^2| \le -\tan \theta \cdot y^1$} \right\} \subseteq \R^2. \end{equation}
	
	The condition $o_{R_k} \in HC$ restricts the form of the lens-shaped sets $D_{R_k} = B(0,2) \cap B(o_{R_k}, R_k)$.  For any point $y \in HC$, we write $D^y = B(0,2) \cap B(y,|y|)$ for the lens-shaped set oriented with old origin $y$, so that $D^{o_{R_k}} = D_{R_k}$.  We then define the compact family of compact sets $\DD := \overline{\{D_y\}}_{y \in HC}.$ The family $\DD$ is compact with respect to the Hausdorff metric on compact subsets of $\R^2$.  As $|y| \to \oo$ along a ray, the sets $D^y$ converge to a half-disk, which is included in the family $\DD$.
	
	Let $\ell_y$ be the tangent line to the ball $B(y,|y|)$ at the origin; equivalently, $\ell_y$ is the tangent line to $D^y$.  The set $D^y$ lies to the left of the line $\ell_y$.  By definition of the hinterland cone $HC$, the line $\ell_y$ meets the vertical-axis at angle less than $\theta$.  	By simple plane geometry, it is easy to see that
		\begin{equation} \label{Dest}
			\mbox{if $D \in \DD$ and $x \in D$, then $x^1 \le \tan \theta \cdot |x^2|$.} \end{equation}

		\begin{figure}[h!] \label{figurecones}
			\includegraphics{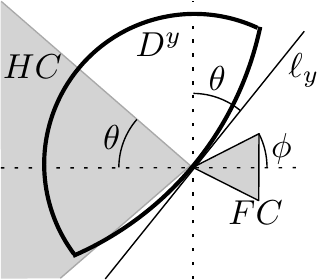}
			\caption{The relationship between the hinterland cone $HC$, the frontier cone $FC$, and a lens-shaped set $D^y$ when $y \in HC$.} %Note that $\theta + \phi < \tfrac \pi 2$, so that the sets $\DD_o$ and $\C$ overlap only at the origin.  Note also that $\C$ is a subset of the unit ball. }
		\end{figure}	
	
	Now, define the angle $\phi := \tfrac 1 2 \big( \tfrac \pi 2 - \theta \big)$.  Since $\theta < \tfrac{\pi}{2}$ by Theorem \ref{frontier_thm}, we have that $\phi > 0$.  We define the \emph{frontier cone}
		\begin{equation} \label{FCdef}
			FC = \left\{ (x^1, x^2) \in \R^2 : \mbox{$0 \le x^1 \le \cos \phi$ and $|x^2| \le \tan \phi \cdot x^1$} \right\} \subseteq \R^2. \end{equation}
	The frontier cone $FC$ is a subset of the ball $B(0,1)$.  
%The next lemma demonstrates that the lens-shaped sets $D^y$ meet the frontier cone $FC$ %only at the origin.  
	
	\begin{env_lem} \label{lem_cones}
		Every set $D \in \DD$ meets the frontier cone $FC$ only at the origin.
	\end{env_lem}
	\begin{proof}
		Let $D \in \DD$, and suppose that $x \in D \cap FC$.  By \eqref{Dest} and the definition \eqref{FCdef} of the set $FC$, $x^1 \le \tan \theta \cdot |x^2|$ and $|x^2| \le \tan \phi \cdot x^1.$ 
		
		If $x^1 = 0$, then $|x^2| \le 0$, so $x=0$.  If $x^1 > 0$, then $x^1 \le \tan \theta \tan \phi \cdot x^1$.  Dividing by $x^1$ and using the sum-of-angles formula for tangent, we have
			$$1 \le 1 -\frac{\tan \theta \tan \phi}{\tan(\theta + \phi)}.$$
		By assumption, $\theta + \phi < \tfrac \pi 2$, so the right side is less than $1$, a contradiction.  Thus $D \cap FC = \{0\}$.
	\end{proof}
	
	This lemma is important in our definition of the bump metric.  For each metric $g$ satisfying the uniformity condition \eqref{Adef}, we will define a bump metric $b(g) \in \Omega_+$ defined on all of $\R^2$.  This bump metric $b(g)$ agrees with $g$ at the origin, and has certain special properties in the frontier cone $FC$. Furthermore, we will show that the function $g \mapsto b(g)$ is continuous on its domain.
	
	\subsection{The Bump Metric} \label{subsect_bump}
	
	We again return to the case that $d=2$, and we are now ready to construct a bump metric $b(g) \in \Omega_+$ for every metric $g$ satisfying the condition $Z_0(g) \le 2h$.  Fix parameters $h > 0$, and define the set of metrics whose fluctuations at the origin are uniformly bounded (as quantified by $Z_0(g) \le 2h$):
		\begin{equation} \label{Adef}
			A_0 = \{ g \in \Omega_+ : Z_0(g) \le 2h \}. \end{equation}
	This is the only place in the paper where we use the assumption that our metrics are $C^2$-smooth.  
	
%	Good explanatory paragraph.  Want to keep it somewhere:  We will construct a surface which has curvature profile $K(t)$ along a geodesic parametrized by time $t$.  The geodesic starts at a point of curvature $K_0$, so that we will be able to match up the metric $g$ with the bump surface at the origin.  Between times $\tau/4$ and $\tau$, the geodesic remains on a spherical region of positive curvature $K_+$.  We will show using Jacobi fields that this implies that the geodesic is not minimizing on the bump surface.
		
	The ``bump map'' is a continuous function $b : A_0 \to \Omega_+$ satisfying a number of nice properties, which are stated precisely in Theorem \ref{bump_thm}.  The bump metric $b = b(g)$ is designed to coincide with $g$ at the origin (up to second derivatives).  It is also designed so that the geodesic $\gamma_b := \gamma(b,\cdot)$ is not minimizing in the frontier cone $FC$.  Furthermore, if $g$ is very close to $b(g)$, then the geodesic $\gamma_g := \gamma(g,\cdot)$ is also not minimizing.
	
%	The bump metric $b(g)$ is an $\Omega_+$-valued random variable, and is measurable with respect to the $\sigma$-algebra $\F_0$, consisting of all the metric information at the origin.

	\begin{env_thm}[Existence of Bump Metrics] \label{bump_thm}
		Suppose $d=2$, fix parameters $h \ge 0$ and $\theta \in [0,\tfrac \pi 2)$, and let $A_0$ be as in \eqref{Adef}.  There exists a continuous function $b : A_0 \to \Omega_+$ such that
			\begin{itemize}
				\item The bump metric $b = b(g)$ agrees with $g$ up to second derivatives at the origin: $\|g - b\|_{C^{2,1}(0)} = 0.$ This implies that their scalar curvatures at the origin agree: $K_0(g) = K_0(b)$.
							
				\item There exists a constant $\tau \in (0,1]$ (independent of $g$) such that for all bump metrics $b \in b(A)$, the geodesic $\gamma_b := \gamma(b,\cdot)$ is not minimizing between times $0$ and $\tau$.  
				
				\item There exists a constant $\epsilon > 0$ (independent of $g$) such that if $\|g - b(g)\|_{C^{2,1}(FC)} < \epsilon$, then $\gamma_g := \gamma(g,\cdot)$ is not minimizing between times $0$ and $\tau$.
			\end{itemize}
		The construction $b(g)$ is $\F_0$-measurable, that is, the bump metric $b(g)$ only depends on the metric $g$ and its derivatives at the origin.
	\end{env_thm}

	We will prove this theorem in Section \ref{sect_bump}.  The condition $g \in A$ implies that the scalar curvature at the origin, $K_0(g)$, satisfies a strong boundedness condition:  $|K_0(g)| \le \Kmax$ for some value $\Kmax$ depending only on the parameter $h$.  We will define a particular curvature profile $K(t)$ which begins at the value $K_0(g)$, then transitions to some value $K_+$.  To realize such a construction, we use Fermi Normal Coordinates adapted to the geodesic starting at the origin in the horizontal direction $\E_1$.  
	
	We first define the curve $\gamma_b$ as a vector-valued polynomial function of $t$, then we construct the bump metric using this curve.  More careful analysis ensures that the bump geodesic $\gamma_b$ lies in the interior of the frontier cone $FC$ for time $(0, \tau]$.  By construction, the geodesic $\gamma_b$ spends time $\tfrac{\tau}{2}$ on a region of constant curvature $K_+ := \tfrac{4\pi^2}{\tau^2}$.  We exactly solve the Jacobi equation \eqref{jacobieqn_b}, and show that it vanishes at times $\tfrac{\tau}{4}$ and $\tfrac{3\tau}{4}$.  Therefore, the points $\gamma(\tfrac{\tau}{4})$ and $\gamma(\tfrac{3\tau}{4})$ are conjugate, hence the geodesic is not minimizing past them.  This argument is essentially a weak form of the Bonnet-Myers theorem \cite{lee1997rmi}.
	
	It is a little trickier to show that this property is preserved under a uniform perturbation of the bump metric.  The key is that the solutions to the Jacobi equation \eqref{jacobieqn_general} vary continuously in the metric parameter $g$.  Thus the solution must change sign, hence vanish.  Again, the geodesic $\gamma_g$ will not be minimizing past critical points.
	
	The value $\tau$ is the natural length scale for the bump metric.  This value is carefully chosen in \eqref{tau} to satisfy multiple technical conditions. We emphasize that the constant $\epsilon$ is non-random and independent of the metric $g$.  This construction uses the fact that the space of bump metrics $b(A_0)$ is compact.
		
	\begin{env_rem}
		There is no mathematical obstruction to extending Theorem \ref{bump_thm} to higher dimensions $d > 2$.  In the general case, the Fermi normal coordinates take the canonical expression \eqref{gexpansion} involving the Riemann curvature tensor $R_{ijkl}$ instead of the scalar curvature $K$.  Under these coordinates, the curvature along the geodesic $\gamma_b$ will start at $R_{ijkl}(g,0)$ at time $t=0$, then transition to constant sectional curvature $K_+$.  The argument involving the Jacobi equation extends without difficulty.
	\end{env_rem}

	Define the open set
		\begin{equation} \label{U_def}
			U = \{ g \in \Omega_+ : \mbox{$Z_0(g) < 2h$ and $\|g - b(g)\|_{C^{2,1}(FC)} < \epsilon$} \} \end{equation}
	of metrics which satisfy the strong regularity estimate at the origin, and which are also close to their associated bump metrics.  Theorem \ref{bump_thm} implies that if $g \in U$, then $\gamma_g$ is not minimizing between times $0$ and $\tau$.  Since $Z_0$ is $\F_0$-measurable, and the frontier cone $FC$ is a subset of the unit ball $B(0,1)$, the event $U$ is $\F_1$-measurable.  It is easy to see that the set $U$ is non-empty (this follows from Lemma \ref{U_lem}).  The set $U$ is non-empty and open, so $\PP(U) > 0$ by strict positivity of the measure $\PP$.  
	
	Consider the family $\DD$ of lens-shaped sets generated by the hinterland cone $HC$ (defined in Section \ref{subsect_hfcones}).  Let $P_D(g,\cdot) = \PP(\cdot|\F_D)$ be the conditional probability defined by Theorem \ref{P_lem} of Part II, and let $[g]_D$ be the equivalence class of metrics which agree with $g$ on the set $D$.\footnote{That is, $g' \in [g]_D$ if and only if $\|g' - g\|_{C^{2,1}(D)} = 0$.}  Part \eqref{P_lem_totallypositive} of Theorem \ref{P_lem} states that if the open set $U$ meets $[g]_D$, then $P_D(g,U) > 0$.
	
	This condition is certainly not satisfied for arbitrary old origins $y$ and metrics $g$.  For example, if $y$ is a point on the positive horizontal axis with $y^1 \ge 1$, then the frontier cone $FC$ is a subset of $D^y$.  Choose any metric $g_0 \in U$, and pick a non-zero point $x \in FC \subseteq D^y$.  Now let $g$ be any metric which equals $g_0$ at the origin (so that $b(g) = b(g_0)$), but for which $|g_{11}(x) - b(g)_{11}(x)| \ge \epsilon$.  Any metric $\tilde g \in [g]_{D^y}$ consequently satisfies $\|\tilde g - b(\tilde g)\|_{C^{2,1}(FC)} \ge \epsilon$, so $U \cap [g]_{D^y}$ is empty.
	
	Again, the crucial condition here is the construction of the hinterland and frontier cones.			
		
	\begin{env_lem} \label{U_lem}
		If $D \in \DD$ and $Z_0(g) < 2h$, then the set $U$ meets the equivalence class $[g]_D$.  %Consequently, $\inf P_D(g,U) > 0$.
	\end{env_lem}
	\begin{proof}
		Since $Z_0(g) < 2h$, Theorem \ref{bump_thm} applies and there exists a well-defined bump metric $b(g)$. By Lemma \ref{lem_cones}, the closed sets $D$ and $FC$ meet only at the origin.  By construction, the metrics $g$ and $b(g)$ agree up to second derivatives at the origin.  Consequently, there exists a Riemannian metric $\tilde g \in \Omega_+$ which is equal to $g$ on the set $D$, equal to $b(g)$ on the set $FC$, and smoothly interpolates between the two.  
		
		By construction, $\tilde g \in [g]_D$.  Since $\tilde g = g$ at the origin, their bump metrics are equal:  $b(\tilde g) = b(g)$.  By construction, $\tilde g = b(g)$ on $FC$, so we have that $\|\tilde g - b(\tilde g)\|_{C^{2,1}(FC)} = 0 < \epsilon$.  Consequently, $\tilde g \in U$.  Since $\tilde g \in [g]_D$, this completes the proof.
	\end{proof}
	
	This lemma allows us to get a uniform lower bound on the conditional probabilities $P_D(g,U)$.  Lemma \ref{U_lem} states that the event $U$ satisfies the hypothesis \eqref{uniformassumption} of the Uniform Probability Estimate (Theorem \ref{prop_unifprobest}).  Consequently, that theorem implies that the lower bound $\inf P_D(g,U)$ is strictly positive.  
	
	\begin{env_pro} \label{U_estimate}
		Let $U$ be the event defined by \eqref{U_def}.  There exists $p > 0$ such that for all $D \in \DD$, if $Z_D(g) \le h$, then $P_D(g,U) \ge p$.
	\end{env_pro}

%	With the construction of the bump metric in hand, and the uniform probability of the event $U$, we are now ready to prove the Main Theorem. 

\section{Proof of Main Theorem} \label{sect_proofofmaintheorem}

	We have set up all the necessary machinery to easily prove the Main Theorem. As throughout, let $\gamma := \gamma_{0,\E_1}(g,\cdot)$ denote the unique unit-speed geodesic starting at the origin in direction $\E_1$.  The Main Theorem states that, with probability one, $\gamma$ is not minimizing.  
	
	\begin{proof}[Proof of the Main Theorem] Let $R_k \uparrow \oo$ be the sequence of frontier radii described in Section \ref{subsect_frontierthm}, and let $W_k = \{ R_k < \oo \}$ be the event that the $k$th frontier radius is well-defined.  Let $W = \bigcap W_k$ be the event that the whole sequence is well-defined.  Corollary \ref{Rk_welldefined} states that for almost every random Riemannian metric $g$ on the event $\{ \mbox{$\gamma$ is minimizing} \}$, the event $W$ is satisfied.  Consequently,
		\begin{equation} \label{notminimizing_Wc}
			\PP( \mbox{$\gamma$ is minimizing} \, | W^c ) = 0. \end{equation}
	
	Define the random variable $\Tmin := \sup \{ t > 0: \mbox{$\gamma$ is minimizing between times $0$ and $t$} \}$, which measures the maximum length of time that the geodesic $\gamma$ is minimizing.  Clearly, $\{ \mbox{$\gamma$ is minimizing} \} = \{ \Tmin = \oo \}$.  On the event $W^c$, $\Tmin < \oo$ almost surely, though we do not have any quantitative estimates on the distribution of $\Tmin$. 
	
	The situation is different on the event $W$.  To prove the Main Theorem, we treat each frontier radius $R_k$ as a new opportunity to see a bump surface.  Let $U$ be the event that a metric is locally like a bump surface, as defined in \eqref{U_def}.  Let $U_k$ be the event that $\sigma_{\tau_{R_k}} g \in U$, defined formally in Section \ref{subsect_repeatedevents}; the event $U_k$ implies that just after the exit time $\tau_{R_k}$, the geodesic $\gamma$ runs over a bump surface and is not length-minimizing.  In particular, the event $U_k$ implies that $\Tmin < \tau_{R_k} + \tau$, where $\tau \le 1$ is the constant described in Theorem \ref{bump_thm}.
	
	By definition, the POV metrics $\sigma_{\tau_{R_k}} g$ each satisfy a strong regularity property and exit angle condition near the origin; this is stated precisely as \eqref{frontier_radii_statement}.\footnote{Equivalently, $g$ satisfies this regularity property near $\gamma(\tau_{R_k})$.  The exit angle condition translates into the condition that the old origin lies in the hinterland cone $HC$.}  Using Proposition \ref{U_estimate}, this gives a uniform probability estimate $P_{D_{R_k}}( \sigma_{\tau_{R_k}}g, U ) \ge p$.  This is the necessary condition \eqref{lem_musthappen_estimate} for the Inevitability Theorem (Theorem \ref{lem_musthappen}) to apply, which then guarantees that the sequence of events $U_k$ occurs infinitely often. This completes the proof of the Main Theorem.
	\end{proof}
	
	Without much difficulty, we can get a quantitative estimate for the time $\Tmin$ conditioned on the event $W$.  Theorem \ref{lem_musthappen} also states that the first occurrence value $K = \inf \{ k : \mbox{$U_k$ occurs} \}$ is a random variable with exponential tail decay on the event $W$.  That is, $\PP( K > k | W) \le (1-p)^k$.  It is not hard to extend this to a similar exponential-decay estimate for the random variable $\Tmin$, which we do in the next and final theorem of the paper.
	
	\begin{env_thm} \label{notminimizing_Wthm}
		There exist positive constants $c$ and $C$ such that 
			\begin{equation} \label{notminimizing_W}
				\PP\big( \mbox{$\gamma$ is minimizing between times $0$ and $t$} \, \big| W \big) \le \PP( \Tmin > t \, | W ) \le C \E^{-ct}. \end{equation}
		Consequently, with probability one, $\gamma$ is not a minimizing geodesic.  
	\end{env_thm}
	\begin{proof}
		Let $T_k = \tau_{R_k}$ be the exit time of the geodesic $\gamma$ from the ball of radius $R_k$, so that $R_k = |\gamma(T_k)|$.  Define the random variable $K := \inf\{ k : \mbox{$U_k$ occurs and $R_k \ge \Rshape$}\},$ where $\Rshape$ is the (random) radius after which the Shape Theorem applies (cf. Theorem \ref{shapecor}).  By definition of the event $U$ (i.e., the construction of the bump metric), $\gamma$ is not minimizing between $0$ and $T_K + \tau \le T_K + 1 \le 2T_K$; the second inequality is a trivial estimate.  By definition of $K$, $R_K \ge \Rshape$, so the Shape Theorem applies and $T_K \le 2\mu R_K$.  By Corollary \ref{Rk_welldefined}, there exists a constant $c_1 \ge 1$ such that $R_k \le c_1 k$.  Thus
		$$\Tmin \le 2 T_K \le 4\mu R_K \le 4\mu c_1 K.$$
	Let $k = \lfloor t / 4\mu c_1 \rfloor$ be the largest integer less than $t/4\mu c_1$, so that trivially, $k \ge t/8\mu c_1$.  By construction, if $\Tmin > t$ then $K > k$, hence
		\begin{equation} \label{blaaah}
			\PP( \Tmin > t \, | W_k ) \le \tfrac{1}{\PP(W_k)} \PP( \mbox{$\Tmin > t$, $K>k$ and $W_k$} ) \le \tfrac{1}{\PP(W_k)} \EE\big[ \PP( U_1^c \cap \dots \cap U_k^c | \F_{R_k}) 1_{W_k} \big] \le (1-p)^k \end{equation}
	by Theorem \ref{lem_musthappen}.
	
	Observe that trivially, $k \ge t/8\mu c_1$.  Combining this with \eqref{blaaah}, we have that $\PP( T_* > t \, | W ) \le \tfrac{1}{\PP(W)} (1-p)^{t/8\mu c_1}.$ Set $C = \tfrac{1}{\PP(W)}$ and $c = -\log(1-p) / 8\mu c_1$.  We have proved statement \eqref{notminimizing_W}, which completes the proof.
	\end{proof}

\part{Proofs of Auxiliary Theorems}

	Having proved the Main Theorem, we now provide proofs of the remaining theorems of the paper.

\section{Proof of the POV Theorem (Theorem \ref{POVthm})} \label{proof_flowthm} %Theorem \ref{POVthm}}

	Suppose that $d \ge 2$; later, we will specialize to the case $d=2$.  The random metric $g$ induces a (random) homogenous Lagrangian $L_g(x,\dot x) = \sqrt{ \langle \dot x, g(x) \dot x \rangle }$.  Consider the tangent bundle $T \R^d$, equipped with \emph{generalized velocity coordinates} $(x, \dot x)$.  Since the manifold $\R^d$ is topologically trivial, we can write $T \R^d \cong \R^d \times \R^d$.  In order to prove the POV Theorem (Theorem \ref{POVthm}), we will need to exploit the rotational invariance of the probability distribution $\PP$.  
	
	As we have seen, geodesics are parametrized proportionally to Riemannian arc length.  This implies that the geodesic flow preserves the Riemannian ``energy shells'' $\{ (x, \dot x) : g_{ij}(x) \dot x^i \dot x^j = \mathrm{constant} \}$.  Unfortunately, these energy shells depend on the metric parameter $g$, and are not invariant under isometries of $\R^d$.  Geodesics are parametrized by arc length; this conserved quantity lets us reduce the dimensionality of the system.  To circumvent the problem of energy shells, we introduce \emph{normalized velocity coordinates} $(x,v)$ on $T \R^d \cong \R^d \times S^{d-1}$, where we define $v = \dot x / |\dot x|$.  To recover the generalized velocity $\dot x$ from $v$ using the metric $g$, we set $\dot x = \lambda_g(x,v) v$,  where $\lambda = \lambda_g(x,v) = 1/\sqrt{\langle v, g(x) v \rangle}$, noting that $|\dot x| = \lambda$.  
		
	We now wish to describe the geodesic flow in the flat velocity coordinates $(x,v)$.  Let $a = a_g(x,v) = -\lambda(x,v) \Gamma_{ij}^k(g,x) v^i v^j \, \E_k$ denote the covariant derivative of $-\lambda v$ in the direction $v$; we use the letter $a$ to denote acceleration.  For each $g \in \Omega_+$, define the vector field $U_g : \R^d \times S^{d-1} \to \R^d \times \R^d$ by $U_g(x,v) := \big( \lambda v, a - \langle a, v \rangle v \big)$. 
	
	The first component of $U_g$ is $\dot x$; the second component is the projection of $a$ onto the hyperplane in $\R^d$ normal to $v$.  If $(X_t, V_t)$ denotes a solution to the differential equation
		\begin{equation} \label{Ug_eqn}
			\tfrac{\D}{\D t} (X_t, V_t) = U_g(X_t, V_t), \end{equation}
	then clearly $\tfrac{\D}{\D t} |V_t|^2 = 2 \langle V_t, \dot V_t \rangle = 0$.  Since $V_0$ is a unit vector, this implies that $V_t \in S^{d-1}$ for all time $t$.  The next lemma states that in these coordinates, the geodesic equations are given by $\dot X_t = \lambda(X_t, V_t) \, V_t$ and $\ddot X_t = \lambda(X_t, V_t) \, a(X_t, V_t)$.  Clearly, $|\dot X_t| = \lambda(X_t, V_t)$.

	\begin{env_lem} \label{XtVtgeodesicflow}
		For any $g \in \Omega_+$ and $(x,v) \in \R^d \times S^{d-1}$, let $(X_t, V_t)(g,x,v)$ denote the solution to the differential equation \eqref{Ug_eqn}, with initial conditions $(X_0, V_0) = (x,v)$.  The curve $X_t$ is the unit-speed geodesic starting at the point $x$ in the direction $v$ for the metric $g$, and $V_t \in S^{d-1}$ is its direction vector at time $t$.		
	\end{env_lem}
	\begin{proof}
		Let $X_t$ be a unit-speed geodesic, and set $V_t = \dot X_t / |\dot X_t|$ and $A_t = a(X_t, V_t)$.  We will show that $(X_t, V_t)$ solves the differential equation \eqref{Ug_eqn}; since solutions to this equation are unique, this will prove the lemma.  Since $X_t$ is parametrized by unit speed, $1 = \langle \dot X_t, g(X_t) \dot X_t \rangle = |\dot X_t|^2 \lambda^{-2}$ hence $|\dot X_t| = \lambda$ and $\dot X_t = \lambda V_t$.  
		
		Applying the quotient rule for differentiation to the formula $V_t = \dot X_t / |\dot X_t|$, and applying the geodesic equation $\ddot X_t =\lambda A_t$, we have
			\begin{eqnarray}
				\dot V_t &=& \tfrac{1}{|\dot X_t|^2} \big( |\dot X_t| \ddot X_t - \dot X_t \tfrac{\D}{\D t} |\dot X_t| \big) = \tfrac{1}{|\dot X_t|^3} \big( |\dot X_t|^2 \ddot X_t - \dot X_t \langle \ddot X_t, \dot X_t \rangle \big) \nonumber \\
				&=& \tfrac{1}{\lambda^3} \big( \lambda^3 A_t - \lambda V_t \langle \lambda A_t, \lambda V_t \rangle \big) = A_t - \langle A_t, V_t \rangle. \label{dotVt}
			\end{eqnarray}
	\end{proof}

	Formula \eqref{dotVt} states that $\dot V_t$ is the projection of $A_t$ onto the hyperplane in $\R^d$ normal to the direction of the geodesic.

	\begin{env_lem}
		The (Euclidean) divergence of the vector field $U_g$ is
			\begin{equation} \label{divUg}
				\div U_g(x,v) = -\langle \nabla \log \det g(x), \dot x \rangle - 3\langle \ddot x, \dot x \rangle / |\dot x|^2, \end{equation}
		where $\dot x = \lambda v$ and $\ddot x = \lambda a$.
	\end{env_lem}
	\begin{proof}
		Using the chain rule, we easily calculate
			$$\tfrac{\partial \lambda}{\partial x^k} = -\tfrac 1 2 \lambda^3 g_{ij,k} v^i v^j \qquad \mathrm{and} \qquad \tfrac{\partial \lambda}{\partial v^k} = - \lambda^3 g_{ik}(x) v^i. $$
		Combining this with the definition of the Christoffel symbols \eqref{geoquantitiesdef}, we have
			\begin{equation} \label{divproof1}
				\tfrac{\partial}{\partial x^k} \big( \lambda v^k \big) = -\tfrac 1 2 \lambda^3 g_{ij,k} v^i v^j v^k = -\lambda^3 g_{lk} \Gamma^l_{ij} v^i v^j v^k = \lambda^2 g_{lk} a^l v^k = \lambda^2 \langle a, g v \rangle, \end{equation}
		since $a = -\lambda \Gamma_{ij}^l(x) v^i v^j$.  We next differentiate the acceleration $a$ with respect to velocity $v$:
			\begin{equation} \label{divproof2}
				\tfrac{\partial a^l}{\partial v^k} = -\lambda^2 g_{ik} v^i a^l - 2 \lambda \Gamma_{ik}^l v^i \qquad \mathrm{and} \qquad \tfrac{\partial a^k}{\partial v^k} = -\lambda^2 \langle a, g v \rangle - \langle \nabla \log \det g(x), \lambda v \rangle, \end{equation}
				%\lambda \tfrac{\partial}{\partial x^i} \big( \log \det g(x) \big) v^i
		where we have used the well-known contraction relation $\Gamma_{ik}^k = \tfrac 1 2 \tfrac{\partial}{\partial x^i} \log \det g(x)$ for the Christoffel symbols.  Using the above, we calculate the divergence of $U$:
			\begin{eqnarray} 
				\div U_g &=& \tfrac{\partial}{\partial x^k} \big( \lambda v^k \big) + \tfrac{\partial}{\partial v^k} \big( a^k - \langle a, v \rangle v^k \big) \nonumber\\
				&=& \lambda^2 \langle a, g v \rangle -\lambda^2 \langle a, g v \rangle - \langle \nabla \log \det g(x), \lambda v \rangle - \big( \delta_{lj} \tfrac{\partial a^l}{\partial v^k} v^j v^k + \delta_{lk} a^l v^k + \delta_{lj} a^l v^j \big) \nonumber \\
				&=& -\langle \nabla \log \det g(x), \lambda v \rangle + \delta_{lj} \big( \lambda^2 g_{ik} v^i a^l + 2\lambda \Gamma_{ik}^l v^i \big) v^j v^k - 2\langle a, v \rangle. \label{divproof3}
			\end{eqnarray}
		Since $\lambda^2 = 1/\langle v, gv \rangle$, we can simplify $\lambda^2 \delta_{lj} g_{ik} v^i a^l v^j v^k = \lambda^2 \langle a, v \rangle \langle v, gv\rangle= \langle a, v \rangle$.  Using the definition of $a$, we have $\delta_{lj} 2\lambda \Gamma_{ik}^l v^i v^j v^k = -2a^l v^j \delta_{lj} = -2\langle a, v \rangle$.  Using these simplifications, we have proved that
			\begin{equation} \label{divproof4}
				\div U_g(x,v) = -\langle \nabla \log \det g(x), \lambda v \rangle - 3\langle a, v \rangle; \end{equation}
		from this, formula \eqref{divUg} immediately follows.
	\end{proof}
		
	We now restrict our attention to the two-dimensional case $d=2$.  The isometries of the plane induce actions on the space $\Omega_+$ of Riemannian metrics.  First, vector addition in $\R^2$ induces a group action on the space $\Omega_+$ of Riemannian metrics.  For every $x \in \R^2$, define the translation $\tau_x : \Omega_+ \to \Omega_+$ by $(\tau_x g)(u) = g(x+u)$.  Similarly, the group $\SO(2)$ of rotations of $\R^2$ also induces a transformation on $\Omega_+$.  Since $\SO(2)$ is isomorphic to $S^1$, we may parametrize rotations by their action on $\E_1$.  For each unit vector $v = (v^1, v^2)$, let $\OO_v : \R^2 \to \R^2$ be the rotation matrix which sends $\E_1 \mapsto v$ and $\E_2 \mapsto v^\perp := (-v^2, v^1)$.  This induces the map $\O_v : \Omega_+ \to \Omega_+$, implicitly defined as the transformation on the metric $2$-tensor $g_{ij}$ by $\left\langle w_1, (\O_v g)(u) w_2 \right\rangle = \left\langle \OO_v w_1, g(\OO_v u) \, \OO_v w_2 \right\rangle$ for any $u, w_1, w_2 \in \R^2$.  Equivalently, $\O_v g(u) = \OO^T_v g(\OO_v u) \OO_v$.  %In coordinates, $(\O_v g)_{ij}(u) = g_{ab}(\OO_v u) [\OO_v]_i^a [\OO_v]_j^b.$
	
	Let $\sigma_t : \Omega_+ \to \Omega_+$ denote the Lagrangian flow defined by \eqref{sigmadef}.  This flow translates the environment by $\gamma(t)$, and rotates it so that $\dot \gamma(t)$ points in the direction $\E_1$.  Since these are just rigid transformations of the plane, we can describe the flow $\sigma_t$ using the transformations $\tau$ and $\O$. This is like our argument in Proposition \ref{scenery_limit}.
	
	Using the fact that $\gamma(t) = X_t(g,0,\E_1)$ and $v_t = V_t(g,0,\E_1)$, we can express the flow $\sigma_t$ in terms of the transformations $\tau_{X_t}$ and $\O_{V_t}$:
		\begin{equation} \label{sigmaOtau}
			\sigma_t g = \left( \O_{V_t(g,0,\E_1)} \circ \tau_{X_t(g,0,\E_1)} \right) g. \end{equation}
	
	For any $(x,v) \in \R^2 \times S^1$, define $g_{xv} := \O_v \tau_x g$ for the environment centered at $x$ in the direction $v$.  The next lemma lets us relate facts about flows in the environment $g_{xv}$ to flows in the environment $g$. 
			
%	We summarize some important facts in the following lemma, whose proof can be found in Appendix \ref{flowlemmasproof}.  

	\begin{env_lem} \label{flowlemma}
		For all $g \in \Omega_+$ and $(x,v) \in \R^2 \times S^1$, the following statements are true. 
		\begin{enumerate}[a)]
			\item The flow $(X_t, V_t)$ beginning at $(x,v)$ in the environment $g$ has a representation as a flow beginning at $(0,\E_1)$ in the environment $g_{xv}$:
				\begin{equation} \label{XVrep}
					X_t(g, x, v) = x + \OO_v X_t(g_{xv}, 0, \E_1) \qquad \mathrm{and} \qquad V_t(g,x,v) = \OO_v V_t (g_{xv}, 0, \E_1). \end{equation}
			\item The POV evolution of the metric $g_{xv}$ along the geodesic $X_t(g_{xv}, 0, \E_1)$ is equal to the POV evolution of the metric $g$ along the geodesic $X_t(g, x, v)$:
				\begin{equation} \label{transformXV}
					\O_{V_t\left(g_{xv},0,\E_1\right)} \tau_{X_t\left(g_{xv},0,\E_1\right)} \big(g_{xv}\big) = \O_{V_t(g,x,v)} \tau_{X_t(g,x,v)} (g). \end{equation}
				This identity relies on the fact that the group $\SO(2)$ is abelian.
			\item \label{flowlemma_rho} The Jacobian of the coordinate change $(\tilde x, \tilde v) = (X_t, V_t)(g, x, v)$ is
				\begin{equation} \label{rho}
					\rho_t(g, \tilde x, \tilde v) := \exp \!\left( - \int_{-t}^0 (\div U_g)\big(X_s(g, \tilde x, \tilde v), V_s(g, \tilde x, \tilde v) \big) \sD s\right), \end{equation}
				where $\div U_g$ is given by the formula \eqref{divUg}.  This function satisfies $\rho_t(g,\tilde x,\tilde v) = \rho_t(g_{\tilde x \tilde v}, 0, \E_1)$ for all $t$, $\tilde x$ and $\tilde v$.
		\end{enumerate}
	\end{env_lem}
	\begin{proof}[Proof of part (a)]
		Fix $g \in \Omega_+$ and $(x,v) \in \R^2 \times S^1$.  By definition, $X_0(g,x,v) = x$ and $V_0(g,x,v) = v$. For notational simplicity in this proof, we write $X_t$ and $V_t$ without arguments to mean $X_t(g_{xv}, 0, \E_1)$ and $V_t(g_{xv}, 0, \E_1)$, respectively.  Define $\tilde X_t := X_t(g, x, v) = x + \OO_v X_t$ and $\tilde V_t := V_t(g,x,v) = \OO_v V_t$.  Clearly, $\tilde X_0 = x$ and $\tilde V_0 = v$.  We will show that $\tilde X_t$ and $\tilde V_t$ solve the differential equation
			\begin{equation} \label{tildeU}
				\tfrac{\D}{\D t}(\tilde X_t, \tilde V_t) = U_g(\tilde X_t, \tilde V_t). \end{equation}
		Lemma \ref{XtVtgeodesicflow} states that $(X_t, V_t)$ solves this equation.  By uniqueness of solutions, this will imply that $X_t(g,x,v) = \tilde X_t$ and $V_t(g,x,v) = \tilde V_t$ for all $t$. 
		
		Before verifying the equation \eqref{tildeU}, we do some preliminary calculations.  Note that
			\begin{equation} \label{tildecalc1}
				\langle V_t, g_{xv} (X_t) V_t \rangle = \langle \OO_v V_t, g(x + \OO_v X_t) \OO_v V_t \rangle = \langle \tilde V_t, g(\tilde X_t) \tilde V_t \rangle \end{equation}
		by the definition of the transformations $\O_v$ and $\tau_x$.  
		
		Define $\tilde \lambda$, $\tilde \Gamma$ and $\tilde a$ for the appropriate quantities using the metric $g_{xv}$.  Equation \eqref{tildecalc1} immediately implies that $\tilde \lambda(X_t, V_t) = \lambda(\tilde X_t, \tilde V_t)$.  Even though the Christoffel symbols are not tensors, they transform like tensors under the linear coordinate change $\OO_v$, so
			$$\tilde \Gamma_{ij}^k(X_t) = [\OO_v^{-1}]_c^k \Gamma_{ab}^c(x + \OO_v X_t) [\OO_v]_i^a [\OO_v]_j^b.$$
		From this it immediately follows that
			\begin{equation} \label{tildecalc2}
				\OO_v \tilde a(X_t, V_t) = a(\tilde X_t, \tilde V_t), \end{equation}
		hence
			\begin{equation} \label{tildecalc3}
				\langle V_t, \tilde a(X_t, V_t) \rangle = \langle \OO_v V_t, \OO_v \tilde a(X_t, V_t) \rangle = \langle \tilde V_t, a(\tilde X_t, \tilde V_t) \rangle, \end{equation}
		since the Euclidean inner product is invariant under the rotation $\OO_v$.
		
		Write $\OO_v \times \OO_v$ for the transformation on $\R^2 \times \R^2$ defined by $(\OO_v \times \OO_v)(x,w) = (\OO_v x, \OO_v w)$.  We now calculate the left hand side of \eqref{tildeU}, applying equations \eqref{tildecalc1}, \eqref{tildecalc2} and \eqref{tildecalc3}:
			\begin{eqnarray*}
				\tfrac{\D}{\D t}(\tilde X_t, \tilde V_t) &=& (\OO_v \times \OO_v) \tfrac{\D}{\D t} (X_t, V_t) = (\OO_v \times \OO_v) U_{g_{xv}}(X_t, V_t) \\
				&=& (\OO_v \times \OO_v) \big( \tilde\lambda V_t, \ \tilde a - \langle \tilde a, V_t \rangle V_t \big) \\
				&=& = \big( \lambda \tilde V_t, \ a - \langle a, \tilde V_t \rangle \tilde V_t \big) = U_g(\tilde X_t, \tilde V_t)
%				&=& \left( \frac{\tilde V_t}{\sqrt{\langle \tilde V_t, g(\tilde X_t) \tilde V_t \rangle}} , \frac{\langle \tilde V_t, a(\tilde X_t, \tilde V_t) \rangle \tilde V_t - a(\tilde X_t,\tilde V_t)}{\sqrt{\langle \tilde V_t, g(\tilde X_t) \tilde V_t \rangle}}  \right) = U(g, \tilde X_t, \tilde V_t),
			\end{eqnarray*}
		which proves \eqref{tildeU}.
	\end{proof}
	
	\begin{proof}[Proof of part (b)]
		Fix $g \in \Omega_+$ and $(x,v) \in \R^2 \times S^1$.  For clarity, we continue to write $X_t = X_t(g_{xv}, 0, \E_1)$ and $V_t = V_t(g_{xv},0,\E_1)$, 
		
		Unlike the previous proof, we now write $\tilde X_t = X_t(g,x,v)$ and $\tilde V_t = V_t(g,x,v)$, so that equation \eqref{XVrep} implies
			\begin{equation} \label{flowlemmaBproof1}
				\tilde X_t = x + \OO_v X_t \qquad \mathrm{and} \qquad \tilde V_t = \OO_v V_t. \end{equation}
		From this, we easily calculate that
			\begin{equation} \label{flowlemmaBproof2}
				\tau_{X_t} g_{xv}(u) = (\tau_{X_t} \O_v \tau_x g)(u) = g(\OO_v u + \OO_v X_t + x ) = g(\OO_v u + \tilde X_t) = (\O_v \tau_{\tilde X_t} g)(u). \end{equation}				
		The group $\SO(2)$ is abelian, so $\OO_v \OO_{V_t} = \OO_{\OO_v V_t} = \OO_{\tilde V_t}$.  This implies that
			\begin{equation} \label{flowlemmaBproof3}
				(\O_{V_t} \O_v g)(u) = g(\OO_v \OO_{V_t} u) = g(\OO_{\tilde V_t} u) = (\O_{\tilde V_t} g)(u). \end{equation}
		Combining \eqref{flowlemmaBproof2} and \eqref{flowlemmaBproof3} yields $\O_{V_t} ( \tau_{X_t} g_{xv} ) = \O_{V_t} (\O_v \tau_{\tilde X_t} g) = \O_{\tilde V_t} \tau_{\tilde X_t} g,$ which proves \eqref{transformXV}.
	\end{proof}
	
	\begin{proof}[Proof of part (c)]
		We follow the same notation as in part (b).  The formula \eqref{rho} for the Jacobian of the flow $(X_t, V_t)$ follows from Liouville's theorem \cite{arnold1978ordinary}; similar computations are made in Sections 9.2 and 9.3 of Zirbel \cite{zirbel2001lagrangian}.  To check that $\rho_t(g,\tilde x, \tilde v) = \rho_t(g_{xv}, 0, \E_1),$ we need only verify that
			\begin{equation} \label{divUtwosides}
				\div U_g(\tilde X_t, \tilde V_t) = \div U_{g_{xv}}(X_s, V_s). \end{equation}
		By applying the differential equation \eqref{tildeU}, the representation \eqref{flowlemmaBproof1} for $\tilde X_s$ and $\tilde V_s$, and the differential equation $\tfrac{\D}{\D s}(X_s, V_s) = U_{g_{xv}}(X_s, V_s)$, we calculate
			$$U_g(\tilde X_s, \tilde V_s) = \tfrac{\D}{\D s}( \tilde X_s, \tilde V_s ) = \tfrac{\D}{\D s}( x + \OO_v X_s, \OO_v V_s) = (\OO_v \times \OO_v) \tfrac{\D}{\D s}(X_s, V_s) = (\OO_v \times \OO_v) U_{g_{xv}}(X_s, V_s).$$
		Since the divergence is invariant under the rotation $\OO_v \times \OO_v$ on $\R^2 \times S^1$, this implies \eqref{divUtwosides}.
	\end{proof}

	Let $f : \Omega \to \R$ be a positive, integrable function.  We are now ready to prove Theorem \ref{POVthm}, which states that $\int f(\sigma_t g) \sD \PP(g) = \int f(g) \rho_t(g) \sD \PP(g)$, for $\rho_t$ defined by \eqref{rhotg}.  We follow the method of Geman \& Horowitz \cite{geman1975random}, which Zirbel describes quite clearly in Section 4 of \cite{zirbel2001lagrangian}.  
	
	\begin{proof}[Proof of Theorem \ref{POVthm}]
	
	Let $\alpha(x) = \tfrac{1}{2\pi} \E^{-|x|^2/2}$ be the standard Gaussian density function on $\R^2$, and let $\nu$ be the uniform measure on $S^{1}$.  The specific form of $\alpha(x)$ is irrelevant; its only role in the proof is that
		\begin{equation} \label{alpha1}
			\int_{S^1} \int_{\R^2} \alpha(x) \sD x \D\nu(v) = 1. \end{equation}
		
		Let $f : \Omega \to \R$ be an integrable function with respect to $\PP$.  Using the representation \eqref{sigmaOtau} for $\sigma_t g$, the trivial identity \eqref{alpha1}, and the translation-invariance of Lebesgue measure, we compute
		\begin{eqnarray}
			\int_\Omega f(\sigma_t g) \sD \PP(g) &=& \int_\Omega f\big(\O_{V_t(g,0,\E_1)} \tau_{X_t(g,0,\E_1)} g\big) \sD \PP(g) \nonumber \\
			&=& \int_\Omega f\big(\O_{V_t(g,0,\E_1)} \tau_{X_t(g,0,\E_1)} g\big) \left( \int_{S^{1}} \int_{\R^2} \alpha(x + \OO_v X_t(g, 0, \E_1) ) \sD x \D \nu(v) \right) \sD \PP(g). \quad \label{flowthmproof1}
		\end{eqnarray}
		
	We interchange the integrals by Fubini's theorem, and make the change of variables $g \mapsto g_{xv}$ so that \eqref{flowthmproof1} equals
		$$\iint_{\R^2 \times S^{1}} \int_\Omega f\big(\O_{V_t(g_{xv},0,\E_1)} \tau_{X_t(g_{xv},0,\E_1)} g_{xv}\big) \,\alpha(x + \OO_v X_t(g_{xv}, 0, \E_1) )  \sD \PP(g) \sD x \D \nu(v)$$
	since, by assumption, the measure $\PP$ is preserved under rotations and translations of the plane.  Using the representations \eqref{transformXV} and \eqref{XVrep}, this simplifies drastically, and is equal to
		\begin{equation} \label{falpha}
			\iint_{\R^2 \times S^{1}} \int_\Omega f\big(\O_{V_t(g,x,v)} \tau_{ X_t(g,x,v)} g\big) \,\alpha( X_t(g,x,v) )  \sD \PP(g) \sD x \D \nu(v). \end{equation}
	To paraphrase Zirbel \cite[p. 817]{zirbel2001lagrangian}:  the idea of \eqref{falpha} is that the evolution $(X_t, V_t)$ from the initial $(0, \E_1)$ is replaced by an integral over possible starting locations $(x,v)$ conditioned (by $\alpha$) on their image under the geodesic flow.  
	
	We again interchange the integrals, and make the change of variables $\tilde x = X_t(g,x,v) \qquad \mathrm{and} \qquad \tilde v = V_t(g,x,v),$ so that $\O_{V_t(g,x,v)} \tau_{ X_t(g,x,v)} g = \O_{\tilde v} \tau_{\tilde x} g = g_{\tilde x \tilde v}$.  
%	Instead of selecting a random starting point $(x,v)$ and considering the flow forward in time, this corresponds to selecting a random ending point $(\tilde x, \tilde v)$ chosen from the distribution weighted by the Jacobian of the coordinate change.  
	By Lemma \ref{flowlemma}.\ref{flowlemma_rho}, the Jacobian of the coordinate change is $\rho_t(g,\tilde x, \tilde v) = \rho_t(g_{\tilde x \tilde v}, 0, \E_1) =: \rho_t(g_{\tilde x \tilde v})$, defined by the formula \eqref{rho}.  Thus \eqref{falpha} is equal to
		\begin{equation} \label{falphatransformed}
			\int_\Omega \iint_{\R^2 \times S^1} f\big( g_{\tilde x \tilde v} \big) \,\alpha( \tilde x ) \rho_t( g_{\tilde x \tilde v}) \sD \tilde x \D \nu(\tilde v) \sD \PP(g). \end{equation}
	We interchange the integrals a third time, and make the change of variable $g_{\tilde x \tilde v} \mapsto g$, so that \eqref{falphatransformed} equals
		$$ \iint_{\R^2 \times S^{1}} \int_\Omega f(g) \,\alpha( \tilde x ) \rho_t(g) \sD \PP(g) \sD \tilde x \D \nu(\tilde v) = \int_\Omega f(g) \rho_t(g) \sD \PP(g),$$
	where the final step follows from interchanging the integrals one more time, and integrating out $\iint \alpha(\tilde x) \sD \tilde x \D\nu(\tilde v) = 1$.  This completes the proof of \eqref{flowformula}.  
	
	The explicit expression \eqref{rhotg} for $\rho_t(g)$ follows from plugging in the expression \eqref{divUg} for the divergence of $U_g$ to the formula \eqref{rho} for $\rho_t(g)$:
		$$\rho_t(g) = \exp \!\left( - \int_{-t}^0 (\div U_g)\big(X_s(g, 0, \E_1), V_s(g, 0, \E_1) \big) \sD s\right) = \exp \!\left( \int_{-t}^0 \Big( \big\langle \nabla \log \det g(X_s), \dot X_s \big\rangle + 3 \frac{\langle \ddot X_s, \dot X_s \rangle}{\langle \dot X_s, \dot X_s \rangle} \Big) \sD s\right).$$
	\end{proof}

%	We can calculate $\rho_t$ explicitly. calculate
%		\begin{eqnarray*}
%			\rho_t(g) &=& \exp \!\left( \int_0^t (\div U)(X_{-t}(g, \tilde x, \tilde v), V_{-t}(g, \tilde x, \tilde v) \sD s\right) \\
%			&=& \exp \!\left( \int_0^t \frac{\left\langle V_{-s}, 3a(X_{-s},V_{-s}) - \nabla \log \det g(X_{-s}) \right\rangle}{\sqrt{\langle V_{-s}, g(X_{-s}) V_{-s}\rangle}} \sD s \right) \\
%			&=& \dots
%		\end{eqnarray*}
%	\end{proof}

% II
\section{Proof of Local Markov Property and Strong LMP} \label{sect_proofmarkov}

	The Local Markov Property (Theorem \ref{thm_markov}) states that we can estimate local observables at the exit location $\gamma(\tau_r)$ (e.g., $K(\gamma(\tau_r))$).  Before we can prove this theorem, we must further investigate the stochastic process $r \mapsto \sigma_{\tau_r}g$.  Lemma \ref{sigmataur_lemma} states that this is a stochastic process with jumps.

	\subsection{A Simple Conditioning Lemma}
	
	Fix some compact sets $D, B \subseteq \R^2$.  Let $f : \Omega \to \R$ be a $\F_D$-measurable function, that is, a random variable which only depends on the metric in the region $D$.  Consider the conditional expectation $\EE(f | \F_B)$ given the metric in the region $B$.  Since the metric has finite-range dependence of length $1$, the conditional expectation will only depend on the metric which belongs to the $1$-neighborhood of $D$ in $B$.
	
	\begin{figure}[h!]
			\includegraphics{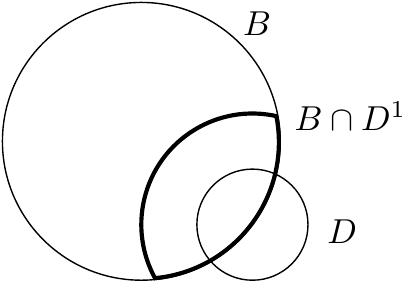}
		\caption{Conditioned on $\F_B$, an $\F_D$-measurable function depends only on the metric in the lens-shaped $B \cap D^1$.}
	\end{figure}
	
	The proof is easy, using standard properties of conditional expectations.
	
	\begin{env_lem} \label{condexp_Bolem}
		Let $f : \Omega \to \R$ be $\F_D$-measurable.  For any compact $B \subseteq \R^d$,
		\begin{equation}
			\EE(f | \F_B) = \EE(f | \F_{B \cap D^1}) \end{equation}
	almost surely.
	\end{env_lem}
	\begin{proof}			
		Fix a compact subset $D \subseteq \R^d$ and an $\F_D$-measurable function.  Since the set $B$ is the union of the sets $A := B \cap D^1$ and $\neg A := \overline{B - D^1}$, the $\sigma$-algebra $\F_B$ has the representation $\F_B = \sigma\big( \F_{A }, \F_{\neg A} \big).$ Define the $\sigma$-algebra $\F_{A }^{\perp} := \sigma \{ f \in \Omega^* : \EE(f|\F_{A }) = \EE f \}$ generated by linear functionals which are uncorrelated with $\F_{A}$.  It is easy to see that
			\begin{equation} \label{condexp_proof1}
				\F_B = \sigma( \F_{A }, \F_{A }^{\perp} \cap \F_{\neg A} ). \end{equation}
		
		The sets $D$ and $\neg A$ are separated by distance at least $1$, so the $\sigma$-algebras $\F_D$ and $\F_{\neg A} $ are independent.  It trivially follows that 
			\begin{equation} \label{condexp_proof2}
				\mbox{the $\sigma$-algebra $\F_{A }^{\perp} \cap \F_{\neg A}$ is uncorrelated with $\sigma(\F_{A }, \F_D)$.} \end{equation}
		
		Now, let $f$ be $\F_D$-measurable.  Taking the conditional expectation with respect to $\F_B$, and using \eqref{condexp_proof1} and \eqref{condexp_proof2}, we have that $\EE(f | \F_B) = \EE\left(f \big| \sigma( \F_{A }, \F_{A }^{\perp} \cap \F_{\neg A} ) \right) = \EE(f | \F_{A } )$ almost surely.\footnote{This follows from Theorem 9.7.k of \cite{williams1991pm}, which is stated in the case of independent $\sigma$-algebras, but Williams's proof only requires the assumption that they are uncorrelated.}
	\end{proof}

	\subsection{Proof of the Local Markov Property}

		Recall that $D_r(g) = D^1 \cap B( o_r(g), r )$.  The Local Markov Property \eqref{markovstatement} is equivalent to
			\begin{equation} \label{markovstatement2}
				\int_{\{\tau_r < \oo\}} F(g) f(\sigma_{\tau_r} g) \sD \PP(g) = \int_{\{\tau_r < \oo\}} F(g) \, P_{D_r}( \sigma_{\tau_r} g, f ) \sD \PP(g), \end{equation}
		for any functions $f$ and $F$ which are $\F_D$-measurable and $\F_r$-measurable, respectively.  Assume furthermore that the functions $f$ and $F$ are non-negative, continuous and bounded functions on $\Omega$.  The general statement will follow by standard approximation arguments.
				
		The Local Markov Property \eqref{markovstatement2} states that if we condition on an arbitrary $\F_r$-measurable random variable $F$ and the event $\{\tau < \oo\}$, then we can express the conditional expectation of $f$ at the exit location $\gamma(\tau_r)$ in terms of the conditional measure $P_{D_r}$.
		
		Let $\ell_r(g,\D t) = \delta\big( \tau_r(\sigma_{-t} g) - t \big)$ be the history measure introduced in Section \ref{sect_exittime}, and let $\ell_r^\epsilon(g,t) := 1_{[-\epsilon,0]}( \tau_r(\sigma_{-t} g) - t)$ be a density function approximation.  
		
	\begin{proof}[Proof of Theorem \ref{thm_markov}: Local Markov Property]
		
		Applying the flow formula \eqref{flowformula_taur} to the left side of \eqref{markovstatement2}, and using the approximation $\ell^\epsilon$ for $\ell$ as in \eqref{ell_approx}, we have				
			\begin{eqnarray}
				\int_{\{\tau_r < \oo\}} F(g) f(\sigma_{\tau_r} g) \sD \PP(g) &=& \int_\Omega f(g) \left( \int_0^\oo F(\sigma_{-t} g) \, \ell_r(g, \D t) \right) \D \PP(g) \nonumber \\
				&=& \lim_{\epsilon \to 0} \int_0^\oo \int_\Omega f(g) \, F(\sigma_{-t} g) \, \ell_r^\epsilon(g, t) \sD \PP(g) \sD t =: \EE(f F_t) \label{markovproof1}
			\end{eqnarray}
		where we use Fubini's theorem, interchange the limit and the integral by the monotone convergence theorem, and write $F_t := F(\sigma_{-t} g) \, \ell_r^\epsilon(g, t)$.
		
		We now consider the expectation $\EE(f F_t)$ with $t$ fixed.  By construction, $\ell_r^\epsilon(g,t)$ is equal to zero unless there is an old origin on the interval $[t, t+\epsilon]$.  When $\ell_r(g,t) \ne 0$, we say that $t$ is an $\epsilon$-approximation to a historical exit time.
		
		In \eqref{Gtdef}, we will introduce a certain $\sigma$-algebra $\A_t$.  By construction, $F_t$ will be $\A_t$-measurable.  Consequently, we can use properties of conditional expectations to decouple $f$ and $F_t$:  $\EE(f F_t) = \EE\big( \EE(f F_t | \A_t) \big) = \EE\big( F_t \EE(f|\A_t) \big)$; from this, the Local Markov Property will easily follow.  
		
		For each $t \ge 0$, define the random variable $T_t(g) = \inf\big( \T_r(g) \cap [t,\oo) \big),$ setting $T_t = \oo$ when $t > \sup \T_r(g)$.  The random variable $T_t$ denotes the first ``old origin time'' at or after $t$.  By upper-semicontinuity of the function $\tau_r(\sigma_{-t} g) - t$, we have that $T_t(g) \in \T_r(g)$.  In particular, if $t \in \T_r(g)$, then $T_t = t$.  Let $o_t(g) = \gamma(g, -T_t(g) )$ be the position of the first old origin at or after time $t$, where we set $o_t = \oo$ if $T_t = \oo$.  %The $\sigma$-algebra $\A_t$ will represent the information contained in the random Euclidean ball $B(o_t, r)$.		
		
		It follows from the relative distance formula \eqref{relativedistance} that for $\PP$-almost every $g$ and every $t \ge 0$, $|o_t(g)| = r$, i.e., the old origin is Euclidean distance $r$ away from the new origin. Define the random set $\tilde D_t(g) := D^1 \cap B\big( o_t(g), r \big).$ For every $g \in \Omega_+$, the functions $t \mapsto o_t(g)$ and $t \mapsto \tilde D_t(g)$ are left-continuous.
		
		We next introduce the $\sigma$-algebra $\A_t$ representing the information contained in the random Euclidean ball $B(o_t(g), r)$.  Of course, since the ball itself is random, our construction has to be a bit delicate.  Heuristically, an event $A$ is $\A_t$-measurable if, conditioned on the event $\{o_t \approx q\}$, it is $\F_{B(q,r)}$-measurable.  To make this construction precise, we introduce a discrete approximation, and let $q$ range over a countable set of values.
		
		Let $Q$ be a countable dense set in the circle of radius $r$ (i.e. $Q$ is dense in $\partial B(0,r)$).  Define
			\begin{equation} \label{Gtdef}
				\A_t = \left\{ A \in F : \forall \, \delta > 0 \mathrm{~and~} q \in Q, \quad A \cap \{ |o_t(g) - q| < \delta \} \in \F_{B(q, r+\delta)} \right\}. \end{equation}
		Let us remark that $\A_t$ is \emph{not} a filtration.  We can determine if $o_t(g)$ is within $\delta$ of an approximating point $q \in Q$ by looking only at the metric in the ball $B(q, r+\delta)$. The next lemma states properties of the $\sigma$-algebra $\A_t$.  
		
		\begin{env_lem} \label{markov_Gtlem}
			~\begin{enumerate}[a)]
				\item The old origin $o_t$ is $\A_t$-measurable, hence so is the random set $\tilde D_t$.
				\item If $F$ is $\F_r$-measurable, then $F_t = F(\sigma_{-t} g) \, \ell_r^\epsilon(g, t)$ is $\A_t$-measurable.  
				\item If $f$ is $\F_1$-measurable, then for $\PP$-almost every $g$ on the event $\{ o_t \ne \oo \}$, $P_{\tilde D_t(g)}(g, f)$ is a version of the conditional expectation $\EE(f | \A_t)$.
			\end{enumerate}
		\end{env_lem}

		The proof of this lemma is relatively straightforward (although technical), and follows from the definitions.  Parts (a) and (b) of this lemma are easy, and just rely on the fact that the geodesic $\gamma$ connecting the origin to $o_t(g)$ lies entirely in the ball $B(o_t, r)$.  Part (c) follows from approximations and Lemma \ref{condexp_Bolem}.  We will prove the lemma after we finish the proof of the Local Markov Property.
		
		Using the $\A_t$-measurability of $F_t$, we use elementary properties of conditional expectations to simplify $\int_0^\oo \EE(f F_t) \sD t$:
			\begin{eqnarray} 
				\lim_{\epsilon \to 0} \int_0^\oo \EE( f F_t ) \sD t &=& \lim_{\epsilon \to 0} \int_0^\oo \EE\big( \EE(f F_t | \A_t) \big) \sD t = \lim_{\epsilon \to 0} \int_0^\oo \EE\big( \EE(f | \A_t) F_t \big) \sD t \nonumber \\
				&=& \lim_{\epsilon \to 0} \int_\Omega \int_0^\oo \EE(f | \A_t) F(\sigma_{-t} g) \, \ell_r^\epsilon(g, t) \sD t \sD \PP(g) \nonumber \\
				&=& \lim_{\epsilon \to 0} \int_\Omega \int_0^\oo F(\sigma_{-t} g) P_{\tilde D_t(g)}(g, f) \, \ell_r^\epsilon(g, t) \sD t \sD \PP(g), \label{markovproof2}
			\end{eqnarray}
		where we plug in the expression for $F_t$, interchange integrals, and replace $\EE(f|\A_t)$ by $P_{\tilde D_t(g)}(g, f)$.
		
		This was the essential step in the proof, and we are now ready to undo all the steps.  We interchange integrals, and make the change of variables $g \mapsto \sigma_t g$, so that \eqref{markovproof2} equals
			\begin{equation} \label{markovproof4}
				\lim_{\epsilon \to 0} \int_\Omega \int_0^\oo F(g) \, P_{\tilde D_t(\sigma_t g)}(\sigma_t g, f) \, \ell_r^\epsilon(\sigma_t g, t) \sD t \sD \PP(g). \end{equation}
		We now wish to undo the approximation $\ell_r^\epsilon$.  We first bring the limit into the integral by the monotone convergence theorem.  	By construction, the function $\ell_r^\epsilon(\sigma_t g, t) = \delta^\epsilon\big( \tau_r(\sigma_{-t} \sigma_t g) - t \big) = \delta^\epsilon \big( \tau_r(g) - t \big)$ is supported on the interval of width $\epsilon$ \emph{behind} the time $t = \tau_r(g)$.  This means that for almost every $t$ (with respect to the measure $\delta^\epsilon \big( \tau_r(g) - t \big) \sD t$), the integrand is constant, and equal to $F(g) P_{\tilde D_{\tau_r}(\sigma_{\tau_r} g)}(\sigma_{\tau_r} g, f)$.  The integrand no longer depends on the variable $t$, so we integrate out $\int_0^\oo \ell_r^\epsilon(\sigma_t g, t) \sD t = 1$ to get that \eqref{markovproof4} equals
			\begin{equation} \label{markovproof5}
				\int_\Omega F(g) P_{\tilde D_{\tau_r}(\sigma_{\tau_r} g)}(\sigma_{\tau_r} g, f) \sD \PP(g). \end{equation}
		
		It is easy to see that $\tau_r \in \T_r(\sigma_{\tau_r} g)$, since $\gamma(\sigma_{\tau_r} g, -\tau_r) = o_{\tau_r}(\tau_r g)$ marks the location of the origin after the point-of-view transformation.  This implies that
			$$\tilde D_{\tau_r(\sigma_{\tau_r} g)}(\sigma_{\tau_r} g, f) = B(0,2) \cap B\big( o_{\tau_r}(\sigma_{\tau_r} g), r \big) = B(0,2) \cap B\big( \gamma(\sigma_{\tau_r} g, -\tau_r), r \big) = D_r,$$
		exactly by the definition \eqref{Drdef} of $D_r$.  Applying this to \eqref{markovproof5} completes the proof of \eqref{markovstatement2} for the case of non-negative, continuous and bounded functions $f$ and $F$.  The general statement follows by standard approximation arguments.  This completes the proof of the Local Markov Property.
	\end{proof}

		\begin{proof}[Proof of Lemma \ref{markov_Gtlem}]
			
			~
			
			\emph{Proof of part (a).}  By construction, $\big| \gamma(g,-s) - o_t(g) \big| \le r$ for all $s \in [0,t]$, so on the event $\{ |o_t - q| < \delta \}$, we have $\big| \gamma(g,-s) - q \big| \le r + \delta$ by the triangle inequality.  The path of the geodesic $\gamma$ is entirely determined by the metric in the ball $B(r,q+r)$, so 
				\begin{equation} \label{gammameas}
					\mbox{the curve $\gamma|_{[0,-t]}$ is $\F_{q,r+\delta}$-measurable on the event $\{ |o_t - q| < \delta \}$.} \end{equation}
			In particular, $o_t(g) = \gamma(g,-t)$ is $\F_{B(q,r+\delta)}$-measurable. % \newline

			\emph{Proof of part (b).}  By \eqref{gammameas}, $\gamma|_{[0,-t]}$ is $\F_{B(q,r+\delta)}$-measurable on the event $\{ |o_t - q| < \delta \}$.  The transformation $g \mapsto \sigma_{-t} g$ only depends on the position and velocity of the geodesic at time $t$, hence is also $\F_{B(q,r+\delta)}$-measurable on that event.  Since $F$ is $\F_r$-measurable, the function $g \mapsto F(\sigma_{-t} g)$ depends only on the metric in the ball $B(o_t, r) \subseteq B(q,r+\delta)$, hence is $\F_{B(q,r+\delta)}$-measurable.  
			
			Let $T_t$ be the old origin time, so that $o_t(g) = \gamma(g, -T_t)$, and $T_t$ is $\F_{B(q,r+\delta)}$-measurable.  By construction, the function $\ell_r^\epsilon(g, t)$ is supported on the set of times $s \in [T_t-\epsilon, T_t]$.  In particular, this is $\F_{B(q,r+\delta)}$-measurable since $T_t$ is. % \newline
			
			\emph{Proof of part (c).}  For each $k \ge 1$, let $Q_k := \{q_{k,1}, \cdots, q_{k,l}\} \subseteq \partial B(0,r)$ be an increasing, non-random, finite set of points on the circle of radius $r$, spaced with distance less than $\delta_k \le 2^{-k}$.  For each $x \in \partial B(0,r)$, assign a unique choice $\psi_k(x) \in Q_k$, satisfying the property that $|\psi_k(x) - x| \le \delta_k \to 0$.
			
			For each $g \in \Omega_+$, define $q_k(g) := \psi_k( o_t(g) )$.  The random variable $q_k$ is a $Q_k$-valued $\delta_k$-approximation of the old origin $o_t$.  This means that for $\PP$-almost every $g$, $q_k(g) \to o_t(g)$ in $\R^2$.
			
			Let $U \in \F$ be any measurable subset of the event $\{o_t \ne \oo\}$, and for each $k \ge 1$ and $q \in Q_k$, define the event $U_{k,q} = U \cap \{ q_k(g) = q \}$.  Write $B_{k,q} = B(q, r+\delta_k)$ for the ball centered at the approximated old origin point $q$.
			
			Fix $k \ge 1$.  On the event $U_{k,q}$, the conditional expectation satisfies $\EE(f | \A_t) = \EE(f | \F_{B_{k,q}})$ almost surely.  Using this property, we compute
				\begin{eqnarray}
					\int_U \EE(f | \A_t) \sD \PP &=& \sum_{q \in Q_k} \int_{U_{k,q}} \EE(f | \F_{B_{k,q}}) \sD \PP = \sum_{q \in Q_k} \int_{U_{k,q}} \EE(f | \F_{B_{k,q} \cap D^1}) \sD \PP  \nonumber \\
					&=& \sum_{q \in Q_k} \int_{U_{k,q}} P_{B_{k,q} \cap D^1}(g,f) \sD \PP(g),  \label{markovlemma_proof1}
				\end{eqnarray}
			where we have applied Lemma \ref{condexp_Bolem}, and the fact that $P_{B_{k,q} \cap D^1}(g,f)$ is a version of the conditional expectation $\EE(f | \F_{B_{k,q} \cap D^1})$ (Theorem \ref{P_lem}.\ref{P_lem_condprob}).  
			
			For simplicity, let us write the (non-random) set $D_{k,q} = B_{k,q} \cap D^1$, and $\tilde D_{k,q}(g) = B(q_k(g), r) \cap D^1$.  Using \eqref{markovlemma_proof1} and the triangle inequality, we compute
				\begin{eqnarray}
					&& \left| \int_U \EE(f | \A_t) - P_{\tilde D_t(g)}(g,f) \sD \PP \right| = \left| \sum_{q \in Q_k} \int_{U_{k,q}} P_{D_{k,q}}(g,f) - P_{\tilde D_t(g)}(g,f) \sD \PP(g) \right| \nonumber \\
					&\le& \sum_{q \in Q_k} \int_{U_{k,q}} \left| P_{D_{k,q}}(g,f) - P_{\tilde D_{k,q}(g)}(g,f) \right| + \left| P_{\tilde D_{k,q}(g)}(g,f) - P_{\tilde D_t(g)}(g,f) \right| \sD \PP(g) \label{markovlemma_proof2}
				\end{eqnarray}
			
			On the event $U_{k,q}$, $q_k(g) = q_k$ almost surely, hence $\tilde D_{k,q}(g) = D_{k,q}$ almost surely, and the first term is trivially equal to zero. The second term is not zero, but vanishes as $k \to \oo$.  The points $q_k(g)$ converge to $o_t(g)$ almost surely on $U$, so the sets $\tilde D_{k,q}(g)$ converge to $\tilde D_t(g)$ almost surely.  Theorem \ref{P_lem}.\ref{P_lem_convergence} implies that the measure-valued function $D \mapsto P_D(g,\cdot)$ is weakly continuous, hence $D \mapsto P_D(g,f)$ is continuous is $D$.  The bounded convergence theorem then implies that the second term vanishes as $k \to \oo$.
			\end{proof}

	\subsection{Proof of the Strong Local Markov Property}

	We now prove the Strong Local Markov Property.  The proof is very similar to part (c) of the above lemma, and features discrete approximations.

	\begin{proof}[Proof of Theorem \ref{thm_strongmarkov}: Strong Local Markov Property]
		By Lemma \ref{RRmax}, $\{ R < \Rmax \}$ is $\F_R$-measurable.  Let $A \in \F_R$ with $A \subseteq \{ R < \Rmax \}$.  We prove formula \eqref{strongmarkovstatement} for $f$ bounded and continuous.  As usual, the standard approximation arguments complete the proof for general $f$.
		
		Consider the dyadic rational numbers $r_{k,l} = \tfrac{l}{2^k}$, and let $Q_k = \{ r_{k,1}, r_{k,2}, \cdots \}$ be the set of $k$-dyadic rationals.  For each $k \ge 0$, let $R_k(g) = \min Q_k \cap [R,\oo)$ be the smallest $k$-dyadic rational number at or above $R$.  The random variables $R_k$ are $\F_R$-measurable, since the event $\{ R_k \ge r_{k,l+1}\}$ is $\F_{B(0, r_{k,l})}$-measurable.
		
		Note that $R_k(g) \downarrow R(g)$ almost surely on $A$.  The map $r \mapsto \sigma_{\tau_r} g$ is upper semi-continuous and $f$ is continuous, so $\lim_{k\to\oo} f(\sigma_{\tau_{R_k}}g) = f(\sigma_{\tau_R}g)$ almost surely.  Thus by the bounded convergence theorem,
			\begin{equation} \label{strongmarkovproperty_proof1}
				\int_A f(\sigma_{\tau_R} g) \sD \PP = \lim_{k \to \oo} \int_A f(\sigma_{\tau_{R_k}} g) \sD \PP = \lim_{k \to \oo} \sum_{l=0}^\oo \int_{A_{k,l}} f(\sigma_{\tau_{r_{k,l}}} g) \sD \PP, \end{equation}
		where we define $A_{k,l} = A \cap \{ R_k(g) = r_{k,l} \}$. By the local Markov property, $f(\sigma_{\tau_{r_{k,l}}} g)$ equals $P_{D_{r_{k,l}}}(g,f)$ on $A_{k,l}$.  As $R_k \downarrow R$,
			$$\lim_{k\to\oo} D_{R_k}(g) = \lim_{k\to\oo} D^1 \cap B\big( \gamma(\sigma_{\tau_{R_k}} g, -\tau_{R_k}), R_k \big) = D^1 \cap B\big( \gamma(\sigma_{\tau_{R}} g, -\tau_{R}), R \big) = D_R(g)$$
		in the Hausdorff topology.  Theorem \ref{P_lem}.\ref{P_lem_convergence} implies that the measure-valued function $D \mapsto P_D(g,\cdot)$ is weakly continuous, hence $D \mapsto P_D(g,f)$ is continuous is $D$.  Applying these facts, \eqref{strongmarkovproperty_proof1} equals
			\begin{equation} \label{strongmarkovproperty_proof2}
				\lim_{k \to \oo} \sum_{l=0}^\oo \int_{A_{k,l}} P_{D_{r_{k,l}}}(g,f) \sD \PP = \lim_{k \to \oo} \int_{A} P_{D_{R_k}}(g,f) \sD \PP = \int_A P_{D_R}(g,f) \sD \PP \end{equation}
		by the bounded convergence theorem.  This completes the proof of the Strong Local Markov Property.		
	\end{proof}

% III
	\section{Continuous Disintegrations of Gaussian Measures} \label{sect_proofctsdisint}

	The goal of this section is to prove Theorem \ref{P_lem}, which quantitatively states that the measure $\PP$ satisfies the continuous disintegration property. The results in this section are stated and proved for general dimension $d \ge 2$. We first develop the theory for the Gaussian measure $\Q$, then everything carries through for the push-forward measure $\PP = \Q \circ \Phi^{-1}$.

	Since $\Omega = C(\R^d, \Sym)$ is a Fr\'echet space of continuous quadratic forms, its dual space $\Omega^*$ has the representation as the space of finite Radon $2$-tensor measures on $\R^d$ with compact support.  More precisely, for every continuous linear functional $f \in \Omega^*$, there exists a Radon $2$-tensor measure $\mu^{kl}_f$ on $\R^d$ with compact support such that for any $\xi \in \Omega$, $f(\xi) = \int_{\R^d} \xi_{kl}(x) \sD \mu^{kl}_f(x)$, where as usual we follow Einstein's convention of summing over repeated indices. 
	
	Let $c : \R \to \R$ be the covariance function described in Section \ref{sect_rrm}, and let $c_{ijkl}(x,y) := c(|x-y|) (\delta_{ik} \delta_{jl}  + \delta_{il} \delta_{jk} )$ be the covariance tensor.  Let $\Q$ be the Gaussian measure on $\Omega$ with mean zero and covariance tensor $c_{ijkl}$.  Define the integral operator $K : \Omega^* \to \Omega$ with kernel $c$:
		\begin{equation} \label{Kdef}
			(Kf)_{ij}(x) := \int_{\R^d} c_{ijkl}(x,y) \sD \mu^{kl}_f(y). \end{equation}
	
	The Structure Theorem (developed by Vakhania et al. \cite{vakhaniya1987probability}, stated as Theorem 1 of \cite{lagatta2010continuous}) completely determines the structure of Gaussian measures on infinite-dimensional topological vector spaces.

	\begin{env_lem} \label{K_lem}  
		The Gaussian measure $\Q$ has mean $0 \in \Omega$ and covariance operator $K : \Omega^* \to \Omega$.  The operator $K$ is continuous, symmetric and injective.  The range of $K$ is dense in $\Omega$, and the support of the measure $\Q$ is the entire space $\Omega$.
	\end{env_lem}
	\begin{proof}
		That $K$ is symmetric follows from the symmetries of the covariance function.  To see that it is continuous, let $f_n \to f$ in $\Omega^*$.  The representation measures $\mu_{f_n} \to \mu_f$ converge weakly.  Since $c$ is continuous and bounded, the quadratic forms $Kf_n$ converge to $Kf$.  
		
		Let $\delta^{ij}_x \in \Omega^*$ be the evaluation functional at $x \in \R^d$, defined by $\delta^{ij}_x(\xi) := \xi_{ij}(x)$; note that $\delta^{ij}_x$ has the representation of a $2$-tensor point-mass measure at $x$.  Clearly, 
			$$\delta^{ij}_x(K \delta^{kl}_y) = c_{ijkl}(x,y) = \EE\big[ \xi_{ij}(x) \xi_{kl}(y) \big] = \int_\Omega \delta^{ij}_x(\xi) \delta^{kl}_y(\xi) \sD \Q(\xi).$$
		Since the evaluation functionals $\delta^{ij}_x$ are linearly dense in $\Omega^*$, this proves that $K$ is the covariance operator of $\Q$. Now, suppose that some $\delta^{ij}_x \in \Omega^*$ is in the kernel of $K$.  Then $0 = (K\delta^{ij}_x)_{ij}(x) = c_{ijij}(x,x),$ a contradiction since this is the variance of $\xi_{ij}(x)$, and is non-zero by assumption.  Thus $K$ is injective.
		
		Recall that the support of the measure $\Q$ is the smallest closed subset of $\Omega$ with $\Q$-probability one.  Vakhania's theorem \cite{vakhania1975topological}) states that $\supp \Q = \overline{K \Omega^*},$ since $\Q$ is a Gaussian measure.  The Gaussian measure is strictly positive:  if $U \in \F$ is any open set, then $\Q(U) > 0$ hence $U \cap \supp \Q$ is non-empty.  The only closed subset of $\Omega$ which meets every open set is the space $\Omega$ itself, so $\overline{K \Omega^*} = \Omega$.
	\end{proof}

%	Let $\C$ denote the space of compact subsets of $\R^d$, equipped with the Hausdorff topology.  
	For any compact set $D \subseteq \R^d$, consider the Banach space $X_D := C(D, \Sym)$ equipped with the supremum norm.  Let $\eta_D : \Omega \to X_D$ denote the restriction map (i.e., $(\eta_D \xi)(x) = \xi(x)$ for all $x \in D$).  The restriction map $\eta_D$ is continuous, linear and surjective onto $X_D$.
	
	Let $\Q_D := \Q \circ \eta_D^{-1}$ denote the push-forward measure of $\Q$ onto $X_D$, so that for any integrable $f : X_D \to \R$, the following change-of-variables formula holds:
		\begin{equation} \label{chvar_gaussian}
			\int_{X_D} f(\hat \xi) \sD \Q_{D}(\hat \xi) = \int_\Omega f(\eta_D \xi) \sD \Q(g). \end{equation}
%	The Gaussian measure is supported on continuous fields which are $C^2$-smooth in the interior of $D$. 
			
	Let $X_D^*$ denote the dual space of $X_D$, and let $\eta_D^* : X_D^* \to \Omega^*$ denote the adjoint map of $\eta_D$, defined by $\eta_D^* f(\xi) = f(\eta_D \xi)$.  For each compact $D \subseteq \R^d$, the Radon measure $\Q_D$ has mean zero in $X_D$ and covariance operator $\eta_D K \eta_D^* : X_D^* \to X_D$. 
	
	Our goal now is to define a continuous, linear operator $m_D : \Omega \to \Omega$, which corresponds to the conditional mean operator for the Gaussian measure. The tensor field $m_D \xi$ represents the best linear predictor of a field $\xi$, given the field information contained over a set $D$. In particular, we define $m_D := \hat m_D \circ \eta_D$, where $\hat m_D : X_D \to \Omega$ extends a tensor field $\hat \xi$ defined only on the set $D$, and extends it to a tensor field $\hat m_D \hat \xi$ defined on all of $\R^d$.  While there are uncountably many ways to do this, it is highly non-trivial that there exists a unique way which is compatible with the probabilistic structure. This is the only place we use the strong assumption of Gaussianity of $\Q$, since the probabilistic structure of a Gaussian measure is characterized by the linear structure of its covariance operator.
	
	%In particular, we will see that for $\Q$-almost every $\xi$, $m_D \xi := \hat m_D( \eta_D \xi) $ is a version of the conditional mean $\EE_{\Q}(\xi | \F_D)$.  The space $\eta_D K \eta_D^* X_D^*$ is dense in $X_D$, and is of particular importance.  
	
	The fundamental fact is that the restriction map $\eta_D$ is invertible on the space $K\eta_D^* X_D^* \subseteq \Omega$.  Furthermore, its inverse $\eta^{-1} : \eta K \eta_D^* X_D^*$ is a continuous, linear operator, with operator norm bounded by $1$ since $\Q$ is stationary. On the space $\eta_D K \eta_D^* X_D^*$, we define $\hat m_D := \eta_D^{-1}$.  Since this operator is continuous and the space $\eta_D K \eta_D^* X_D^*$ is dense in $X_D$, we can extend this to a continuous, linear operator on all of $X_D$.  We then define the map $m_D : \Omega \to \Omega$ by setting $m_D := \hat m_D \circ \eta_D$.\footnote{Since the operator $\hat m_D$ is defined on all of $X_D$, we can actually condition on non-smooth tensor fields, even though they have zero probability of occuring.}

	We now describe properties of the conditional mean operator $m_D$.  For any compact $D \subseteq \R^d$, $\eta_D \circ \hat m_D$ is the identity operator on $X_D$.  The function $\xi \mapsto m_D \xi$ only depends on $\xi$ via the set $D$, hence is $\F_D$-measurable.  The joint continuity property says that if $D_n \to D$ in the Hausdorff metric and $\xi_n \to \xi$ in $\Omega$, then $m_{D_n} \xi_n \to m_D \xi$ in $\Omega$. We will see in Theorem \ref{prop_gaussian} that for $\Q$-almost every $\xi$, $m_D \xi$ is a version of the conditional mean $\EE(\xi|\F_D) \in \Omega$.
	
	\begin{env_lem}[Properties of the conditional mean operator $m_D$] \label{mo_lem}
		Fix a compact set $D \subseteq \R^d$.  There exists a continuous, linear map $\hat m_D : X_D \to \Omega$ (we define $m_D := \hat m_D \circ \eta_D$) such that 
			\begin{enumerate}[a)]
				\item \label{mo_identity} (Identity)  The map $\eta_D \circ \hat m_D$ is the identity operator on $X_D$.  That is, for every $\hat \xi \in X_D$, $\eta_D \hat m_D \hat \xi = \hat \xi$.
				\item \label{mo_finrange}  (Compact support)  For every $\xi \in \Omega$, the function $m_D(\xi) \in \Omega$ is equal to zero on the region $\R^d - D^1$, where $D^1$ denotes the $1$-neighborhood of $D$.  i.e., if $\deuc(x,D) \ge 1$, then $m_D(\xi)(x) = 0$. 
				\item \label{mo_measurability}  ($\F_D$-measurability)  The function $\xi \mapsto m_D \xi$ is $\F_D$-measurable.
				\item \label{mo_continuity}  (Joint continuity)  The function $(D,\xi) \mapsto m_D \xi$ is jointly continuous in $D \in \C$ and $\xi \in \Omega$.
			\end{enumerate}
	\end{env_lem}

	The proof of this lemma can be found in Appendix \ref{proof_mo_lem}.
	
	Let $m_D^* : \Omega^* \to \Omega^*$ denote the formal adjoint of the map $m_D$, defined by the action $m_D^*(f)(\xi) = f \circ m_D(\xi)$ for any $f \in \Omega^*$ and $\xi \in \Omega$.  For each compact $D \subseteq \R^d$, we formally define the \emph{conditional covariance operator} $K_D := K - K m_D^*$. Heuristically, the operator $K_D$ projects onto the randomness generated away from $D$.  For any $x,y \in \R^d$, define the conditional covariance $c_D(x,y) = \delta_x (K_D \delta_y)$, where $\delta_x$ and $\delta_y$ are the evaluation functionals ($c_D$ is a $4$-tensor).  If one of $x$ or $y$ is an element of the set $D \subseteq \R^d$, then $c_D(x,y) = 0$.

	\begin{env_lem}[Properties of the conditional covariance operator $K_D$] \label{Ko_lem}
		For each compact $D \subseteq \R^d$, the linear operator $K_D : \Omega^* \to \Omega$ is a well-defined, symmetric and continuous operator.  As operators, $0 \le K_D \le K$.  The kernel and range of the conditional covariance operator $K_D$ are given by	
		\begin{equation} \label{Ko_kerran}
			\ker K_D = \overline{\eta_D^* X_D^*} \subseteq \Omega^* \qquad \mathrm{and} \qquad \overline{\operatorname{ran}(K_D)} = \eta_D^{-1}(0) \subseteq \Omega. \end{equation}
		The function $(D,f) \mapsto K_D f$ is jointly continuous in $D$ and $f$.
	\end{env_lem}
	\begin{proof}
		The arguments of Lemma 2.3 of \cite{lagatta2010continuous} apply to this situation, and show that $K_D$ is well-defined, continuous and symmetric. We verify that $\ker K_D = \overline{\eta_D^* X_D^*}$.  First, if $e \in X_D^*$, then for any $f \in \Omega^*$, $f( K_D \eta_D^* e) = f\big( K \eta_D^* e - K\eta_D^* m_D^* \eta_D^* e \big) = 0$, so $\eta_D^* e \in \ker K_D$.  Next, suppose that $f \in (K\eta_D^* X_D^*)^\perp$, so $f(K_D f) = f(Kf) - f(K \eta_D^* m_D^* f) = c(0) > 0$, so $f \notin \ker K_D$.  The space $K\Omega^*$ is reflexive \cite{vakhaniya1987probability}, so this proves the claim.  
		
		Next, we verify that the range of $K_D$ is the kernel of the restriction map $\eta_D$ in $\Omega$.  For any $f \in \Omega^*$ and $e \in X_D^*$, $e (\eta_D K_D f) = f(K \eta_D^* e) - f(K \eta_D^* m_D^* \eta_D^* e) = 0$, so $\operatorname{ran} K_D \subseteq \ker \eta_D$.  The other direction follows from similar general arguments as above.  
		
		We next show joint continuity of $(D, f) \mapsto K_D f$.  Suppose that $D_n \to D$ and $f_n \to f$.  For any $f' \in \Omega^*$,
			$$f'(K_{D_n} f_n) = f_n(K_{D_n} f') = f_n(K f') - f'( m_{D_n} \eta_{D_n} K f_n) \to f(Kf') - f'(m_D \eta_D K f) = f'(K_D f),$$
		by the joint continuity of the map $(D,\xi) \mapsto m_D \eta_D \xi$.
	\end{proof}

	Trivially, the empty set $\varnothing$ is a compact subset of $\R^d$.  Clearly, $K_\varnothing = K$, which means that conditioning on nothing yields no information. The next lemma states that the operator-valued function $D \mapsto K_D$ is monotonically decreasing: as we condition on larger sets, the covariance operator gets smaller.

	\begin{env_lem}[Monotonicity of the operators $K_D$] \label{monotoneKD}
		If $D'$ and $D$ are non-empty compact subsets of $\R^d$ with $D' \subseteq D$, then
			\begin{equation}
				0 \le K_{D} \le K_{D'} < K \end{equation}
		as positive operators on $\Omega^*$.  %Equivalently, the operator-valued function $D \mapsto K_D$ is decreasing.
	\end{env_lem}
	\begin{proof}
		The covariance operator $K$ defines a symmetric inner product $\langle f, f' \rangle := f(Kf')$ on $\Omega^*$.  Define an equivalence relation on $\Omega^*$ by $f \sim f'$ if $\| Kf - Kf' \|_{C(\R^d)} = 0$; this is well-defined since the integral kernel $c$ has compact support.  Let $H$ be the Hilbert-space completion of the inner product space $\Omega^* / \sim$, and let $\iota^* : \Omega^* \to H$ be the inclusion map.  Define the unitary map $\iota : H \to \Omega$ first on the dense subspace $\iota^* \Omega^*$ by $\iota(\iota^* f) = Kf$, then extend it continuously to all of $H$.  The operator $K$ factors as $\iota \iota^*$.\footnote{The subspace $\iota H$ of $\Omega$ is called the \emph{Cameron-Martin space} of the Gaussian measure $\Q$.} We summarize this with the following commutative diagram:
		
		\begin{equation} \label{diagram}
			\begin{matrix} \xymatrix{X_D^* \ar@{->}[r]^{\eta_D^*} & \Omega^* \ar@{->}[rr]^{K} \ar@{->}[rd]^{\iota^*} && \Omega\ar@{->}[r]^{\eta_D} & \ar@{->}@/^/[l]^{\hat m_D} X_D \\ && H \ar@{->}[ur]^{\iota} && } \end{matrix} \end{equation}
			
		Let $H_D = \overline{\iota^* \eta_D^* X_D^*}$ be the Hilbert subspace of $H$ generated by the evaluation functionals over the set $D$, and let $\pi_D : H \to H$ be the orthogonal projection onto the subspace $H_D$.  By the same derivation as for equation (2.15) of \cite{lagatta2010continuous}, we have 
			\begin{equation} \label{meta_iota}
				m_D \iota = \hat m_D \eta_D \iota = \iota \pi_D \end{equation}
			on $H$.  Consequently, $K_D = K - K m_D^* = K - K\eta_D^* \hat m_D^* = \iota( \operatorname{Id}_H - \pi_D ) \iota^*$, so $K_D = \iota \pi_D^\perp \iota^*$ where $\pi_D^\perp$ is orthogonal projection onto $H_D^\perp$.  
		
		Using this representation, the lemma follows immediately, since if $D' \subseteq D$, then $K_{D} = \iota \pi_{D}^\perp \iota^* = \iota \pi_{D'}^\perp \pi_D^\perp \iota^* \le \iota \pi_{D'}^\perp  \iota^* = K_{D'}$.
	\end{proof}

	Recall that $\Q$ is the Gaussian measure on $\Omega$ with mean zero and covariance operator $K$.  For any compact $D \subseteq \R^d$, $\Q_D$ is the Gaussian measure on $X_D$ with mean zero and covariance operator $\eta_D K \eta_D^*$, and is the push-forward of $\Q$ onto $X_D$. 
	
	For every compact $D \subseteq \R^d$ and $\xi \in \Omega$, let $Q_D(\xi,\cdot)$ denote the Gaussian measure on $\Omega$ with mean $m_D \xi$ and covariance operator $K_D$. The measure-valued function $Q_D(\xi,\cdot)$ is a continuous disintegration of the Gaussian measure $\Q$.

	\begin{env_thm}[Properties of the conditional Gaussian measures $Q_D$] \label{prop_gaussian}
~%	Fix a compact set $D \subseteq \R^d$, and let $Q_D$ be the Gaussian measure defined in Definition \eqref{QDdef}.  
		 	\begin{enumerate}[a)]
		 		\item \label{prop_gaussian_condprob} (Interpretation as a conditional probability) For $\Q$-almost every $\xi$, the Gaussian measure $Q_D(\xi,\cdot)$ is a version of the conditional probability measure $\Q(\cdot | \F_{D} )$.
		 		\item \label{prop_gaussian_fiber} (Support on the fiber of $\eta_D \xi$)  For every $\xi \in \Omega$, the measure $Q_D(\xi,\cdot)$ is supported on the fiber $\eta_D^{-1}( \eta_D \xi )$ in $\Omega$.\footnote{This means that for $Q_D(\xi,\cdot)$-almost every $\xi'$, $\eta_D \xi' = \eta_D \xi$.  For any $\xi' \in \eta_D^{-1}(\eta_D \xi)$, the measures $Q_D(\xi,\cdot)$ and $Q_D(\xi',\cdot)$ are equal.}
		 		\item \label{prop_gaussian_invariance} (Invariance under $m_D$)  For every $\xi \in \Omega$, the Gaussian measures $Q_D(\xi,\cdot)$ and $Q_D(m_D \xi, \cdot)$ are equal.
		 		\item \label{prop_gaussian_disintegrationeqn} (Disintegration equation) For every integrable $f : \Omega \to \R$,
		 			\begin{equation} \label{disint_gaussian}
		 				\int_{\Omega} f(\xi) \sD \Q(\xi) = \int_{X_D} \int_{\Omega} f(\xi) \, Q_D(\hat m_D \hat \xi, \D \xi) \sD \Q_D(\hat \xi). \end{equation}
%		 		The disintegration equation states that drawing a random $\xi$ with respect to $\Q$ is the same as first drawing a random $\hat \xi$ in $X_D$ from $\Q_D$, then drawing a $\xi$ from $Q_D(m_D \hat \xi, \cdot)$.
		 		\item \label{prop_gaussian_condtotalpositivity} (Conditional strict positivity) Let $U \in \F$ be an open event which meets the fiber $\eta_D^{-1}(\eta_D \xi)$ for some $\xi \in \Omega$.  Then $Q_D(\xi, U) > 0$.		 		
		 		\item \label{prop_gaussian_compactness} (Weak relative compactness) If $D_n \to D$ in the Hausdorff topology and $\xi_n \to \xi$ in $\Omega$, then the Gaussian measures $Q_{D_n}(\xi_n, \cdot)$ converge weakly to $Q_D(\xi, \cdot)$.
		 	\end{enumerate}
	\end{env_thm}
	\begin{proof}
		It follows from Proposition 3.9 of \cite{tarieladze2007disintegration} that since $K_D \le K$, the operator $K_D$ is a Gaussian covariance operator.  Therefore the measures $Q_D(\xi,\cdot)$ are well-defined Gaussian measures.  We will prove part \eqref{prop_gaussian_condprob} as a consequence of the disintegration equation \eqref{disint_gaussian}.

		\emph{Proof of part \eqref{prop_gaussian_fiber}.}  By Vakhania's theorem \cite{vakhania1975topological}, for every $\xi \in \Omega$,
			\begin{equation} \label{suppQ}
				\supp Q_D(\xi, \cdot) = m_D \xi + \overline{K_D \Omega^*},\end{equation}
		Applying the restriction map $\eta_D$ to the right side of this equation, we see that $\eta_D\big( \hat m_D \eta_D \xi + \overline{K_D \Omega^*} \big) = \eta_D \xi,$ since $\eta_D \hat m_D$ is the identity on $X_D$, and $\overline{ \operatorname{ran} K_D } = \ker \eta_D$.  Consequently, the support of $Q_D(\xi,\cdot)$ is $\eta_D^{-1} (\eta_D \xi)$, the fiber of $\xi$ over $D$. If $\eta_D \xi' = \eta_D \xi$, then $Q_D(\xi',\cdot)$ is the Gaussian measure with mean $m_D \xi' = m_D \xi$ and covariance operator $K_D$, hence is equal to $Q_D(\xi,\cdot)$.
	
		\emph{Proof of part \eqref{prop_gaussian_invariance}.}  Part \eqref{mo_identity} of Lemma \ref{mo_lem} states that $\eta_D \hat m_D$ is the identity operator on $X_D$, so $m_D^2 = (\hat m_D \eta_D)^2 = \hat m_D \circ \big( \eta_D \circ \hat m_D\big) \circ \eta_D = \hat m_D \eta_D = m_D.$  The measure $Q_D(m_D  \xi, \cdot)$ is Gaussian with mean $m_D^2 \xi = m_D \xi$ and covariance operator $K_D$, hence is equal to $Q_D(\xi, \cdot)$.
		
		\emph{Proof of part \eqref{prop_gaussian_disintegrationeqn}.}  The disintegration equation for Gaussian measures is the content of Theorem 2 of \cite{lagatta2010continuous}.  We outline the argument here; the reader interested in more details may consult \cite{lagatta2010continuous}.  Let $H$ be the Hilbert space described in the proof of Lemma \ref{monotoneKD}.  Since $\Q$ is a Radon measure, there exists a Gaussian cylindrical measure $\gamma$ on $H$ such the push-forward $\gamma \circ \iota^{-1}$ radonifies $\Q$.  Define the push-forward cylindrical measures $\gamma_D = \gamma \circ \pi_D^{-1}$ and $\gamma_D^\perp = \gamma \circ (\pi_D^\perp)^{-1}$.
				
		The cylindrical measure $\gamma_D$ is supported on the subspace $H_D$ of $H$, and the cylindrical measure $\gamma_D^\perp$ is supported on the orthogonal complement $H_D^\perp$.  A fundamental property of (zero mean) Gaussian measures is that orthogonality implies independence.  In the present context, this means that $\gamma = \gamma_D * \gamma_D^\perp$, so for any measurable function $f : H \to \R$,
			\begin{equation} \label{disint_H}
				\int_H f(h) \sD \gamma(h) = \int_H \int_H f(h + k) \sD \gamma_D^\perp(h) \sD \gamma_D(k). \end{equation}
		It is not hard to see that the measure $\Q_D$ is the radonification of $\gamma_D \circ (\eta_D \iota)^{-1}$ on $X_D$, and $Q_D(0,\cdot)$ is the radonification of $\gamma_D^\perp \circ \iota^{-1}$ on $\Omega$.  By pushing everything forward onto the spaces $\Omega$ and $X_D$ and applying definitions, the disintegration equation \eqref{disint_gaussian} easily follows from \eqref{disint_H}. 		
		
		\emph{Proof of \eqref{prop_gaussian_condprob}.}  Let $f : \Omega \to \R$ be an integrable function and let $U \in \F_D$.  Since $U$ is generated by evaluation functionals on $D$, we have that $\eta_D^{-1}(\eta_D U) = U$.  We calculate
			\begin{eqnarray*}
				\int_U f(\xi) \sD \Q(\xi) &=& \int_{\eta_D (U)} \int_\Omega f(\xi') \, Q_D(\hat m_D \hat \xi, \D  \xi') \sD \Q_D(\hat \xi) \\
				&=& \int_U \int_{\Omega} f( \xi') Q_D(m_D  \xi, \D  \xi') \sD \Q(\xi) = \int_U \int_{\Omega} f(\xi') Q_D(\xi, \D \xi') \sD \Q(\xi),
			\end{eqnarray*}
		where we apply the disintegration equation \eqref{disint_gaussian}, the change of variables formula \eqref{chvar_gaussian}, and the invariance property \eqref{prop_gaussian_invariance}.  This proves that for $\Q$-almost every $\xi$, $\int_{\Omega} f( \xi') Q_D(\xi, \D \xi')$ is a version of the conditional expectation $\EE(f|\F_D)$, as desired. 
		
		\emph{Proof of \eqref{prop_gaussian_condtotalpositivity}.}  Gaussian measures are strictly positive, so if $U$ is an open set which meets $\eta_D^{-1} \eta_D \xi = \supp Q_D(\xi,\cdot)$, then it has positive measure. 
		
		\emph{Proof of \eqref{prop_gaussian_compactness}.}  Fix some compact set $D \subseteq \R^d$ and some $\xi \in \Omega$, and suppose that $(D_n, \xi_n) \to (D, \xi)$.  We first show that the characteristic functionals of $Q_{D_n}(\xi_n, \cdot)$ converge to the characteristic functional of $Q_D(\xi,\cdot)$.  Let $f \in \Omega^*$ be a continuous linear functional.  With respect to the probability measure $Q_D(g,\cdot)$, the random variable $f$ has a real-valued Gaussian distribution with mean $f(m_D \xi)$ and variance $f(K_D f)$.  The characteristic function of $f$ takes the familiar form, and
			\begin{eqnarray*}
				\lim_{n\to\oo} \int \E^{\I f(\xi')} \, Q_{D_n}(\xi_n, \D \xi') &=& \lim_{n\to\oo} \exp\big(\I f( m_{D_n}  \xi_n ) - \tfrac 1 2f(K_{D_n} f) \big) \\
				&=& \exp\big(\I f( m_D \xi ) - \tfrac 1 2 f(K_D f) \big) = \int \E^{\I f(\xi')} \, Q_D(\xi, \D \xi')
			\end{eqnarray*}
		by the continuity of $f$, and the joint continuity of $(D,\xi) \mapsto m_D \eta_D \xi$ and $(D,f) \mapsto K_D f$.
		
		Of course, this is not enough to show weak convergence of probability measures:  we must also show that the sequence $Q_{D_n}(\xi_n, \cdot)$ is tight.  The difficulty is to show that the mean-zero sequence $Q_{D_n}(0,\cdot)$ of measures is tight, as we can recover $Q_{D_n}(\xi_n,\cdot)$ simply by shifting means.  Without loss of generality, we may suppose that the sets $D_n$ are all contained in $D^1$, the $1$-neighborhood of $D$.  
		
		If $x \notin D^2$, then conditioning on $D_n$ or $D^1$ does not affect the distribution of the evaluation functionals $\delta^{ij}_x$:  they are still real-valued Gaussians with mean zero and variance $c(0)$.  Thus we may restrict our attention to the Banach space $\Omega' := C(D^2, \Sym)$ of symmetric quadratic forms on $D^2$, the $2$-neighborhood of $D$.	
		
		For any $h \ge 0$, let $V_h := \{ \xi : \|\xi\|_{C^{0,1}(D^2)} \le h \}$ be the set of functions which satisfy an $h$-Lipschitz estimate on $D^2$.  By the Arzel\`a-Ascoli theorem, the set $V_h$ is compact in the Banach space $\Omega'$.
		
		Let $\epsilon > 0$, and choose $h$ large enough so that $Q_{D^1}(0, V_h) \ge 1 - \epsilon.$  That is, conditioned on the $1$-neighborhood of $D$, the event $V_h$ occurs with probability at least $1-\epsilon$.  By the monotonicity lemma \ref{monotoneKD}, $K_{D_n} \ge K_D$, so the fluctuations for the Gaussian measure $Q_{D_n}$ are greater than those of $Q_D$.  Slepian's inequality (cf. Theorem 2.2.1 of \cite{adler07}) makes this precise, and implies that $Q_{D_n}(0,V_h) \ge Q_{D^1}(0,V_h) \ge 1-\epsilon.$ The set $V_h$ is compact in $\Omega'$, so the sequence of measures $Q_{D_n}(0,\cdot)$ is tight.  
	\end{proof}
	
	We now transfer the disintegration properties of the Gaussian measure $\Q$ to the measure $\PP$.  Let $\varphi : \Sym \to \SPD$ be the function defined in Section \ref{sect_rrm}, and let $\Phi : \Omega \to \Omega$ be the map which acts by $\varphi$ pointwise:  $\Phi \xi(x) := \varphi(\xi(x))$.  We define the measure $\PP$ by pushing forward the Gaussian measure $\Q$ under the map $\Phi$:  $\PP = \Q \circ \Phi^{-1}$.  For any compact $D \subseteq \R^d$ and any Riemannian metric $g \in \Omega_+$, define $P_D(g, U) := Q_D( \Phi^{-1} g, \Phi^{-1} U )$ for any event $U \in \F$. We are now ready to prove that $P$ satisfies all the properties described in Theorem \ref{P_lem}.  Each of these properties follows easily from the corresponding property of the measure $Q_D$.  
	
	\begin{proof}[Proof of Theorem \ref{P_lem}]
		~
		
		\emph{Proof of part \eqref{P_lem_condprob}.}  Let $U \in \F$ be any event, and let $V \in \F_D$.  The operator $\Phi$ acts pointwise, hence preserves the $\sigma$-algebra $\F_D$ so $\Phi^{-1} V \in \F_D$.  Then
			$$\int_V P_D(g,U) \sD \PP(g) = \int_{\Phi^{-1} V} P_D( \Phi \xi, U) \sD \Q(\xi) = \int_{\Phi^{-1} V} Q_D( \xi, \Phi^{-1} U) \sD \Q(\xi) = \Q(\Phi^{-1}V \cap \Phi^{-1}U) = \PP(V \cap U),$$
		by the definition of $P_D$, and property \eqref{prop_gaussian_condprob} of Theorem \ref{prop_gaussian}.
		
		\emph{Proof of part \eqref{P_lem_fiber}.}  By definition, $[g]_D = \Phi \eta_D^{-1} \eta_D \Phi^{-1} g$, so $P_D(g, [g]_D) = Q_D(\Phi^{-1} g, \eta_D^{-1} \eta_D \Phi^{-1} g) = 1$, since $Q_D(\Phi^{-1} g,\cdot)$ is supported on the fiber $\eta_D^{-1} \eta_D \Phi^{-1} g$. If $g' \in [g]_D$, then $\Phi^{-1} g' \in \eta_D^{-1} \eta_D \Phi^{-1} g$, so the measures $Q_D(\Phi^{-1}g', \cdot)$ and $Q_D(\Phi^{-1}g, \cdot)$ are equal, hence so are $P_D(g',\cdot)$ and $P_D(g,\cdot)$. 
		
		\emph{Proof of part \eqref{P_lem_totallypositive}.}  If $U \in \F$ meets $[g]_D$, then $\Phi^{-1} U$ meets $\eta_D^{-1} \eta_D \Phi^{-1} g$, so $Q_D(\Phi^{-1} g, \Phi^{-1} U$ is positive, and $P_D(g, U) = Q_D(\Phi^{-1} g, \Phi^{-1} U) > 0$. 
		
		\emph{Proof of part \eqref{P_lem_convergence}.}  Suppose that $(D_n, g_n) \to (D, g)$.  Then for any open set $U \in \F$, $\liminf P_{D_n}(g_n, U) = \liminf Q_{D_n}(\Phi^{-1} g_n, \Phi^{-1} U) = Q_D(\Phi^{-1} g, \Phi^{-1} U) = P_D(g,U)$, since the function $\Phi^{-1}$ is continuous so $\Phi^{-1} g_n \to \Phi^{-1} g$ and $\Phi^{-1} U$ is open.
	\end{proof}

\section{Proof of Frontier Theorem (Theorem \ref{frontier_thm})} \label{sect_proofoffrontiertimes} %Theorem \ref{frontier}}

	Define $\tau_v(r) := \tau_v(g,r) := \inf \{ t \ge 0 : \gamma_v(t) > r \}$ for the exit time of $\gamma_v$ from the Euclidean ball $B(0,r)$.  It is clear that for all $v \in S^{d-1}$, the random variable $\tau_v(r)$ is $\F_r$-measurable, and the function $r \mapsto \tau_v(r)$ is upper semi-continuous, hence an increasing stochastic process with jumps, adapted to the filtration $\F_r$.  

	\begin{env_lem} \label{lem_taur_estimate}
		Let $\epsilon \in (0,1)$.  With probability one, there exists $r_0$ so that if $r \ge r_0$ and $v \in \V_g$, then $(1-\epsilon) \mu r \le \tau_v(r) \le (1+\epsilon) \mu r.$
	\end{env_lem}
	
	The upper bound is \eqref{taur_upperestimate}; the lower bound is proved similarly following the argument of Theorem \ref{transientgeodesics}.
	
	Define the arccosine of the exit angle $\beta_v(r) = \arccos \alpha_v(r) = \langle \gamma_v, \dot \gamma_v \rangle \big/ r |\dot \gamma_v|,$ where $\gamma_v$ and $\dot \gamma_v$ are evaluated at the exit time $\tau_v(r)$.

	\begin{env_lem}
		The function $r \mapsto \tau_v(r)$ is right-differentiable.  Except at countably many points (corresponding to the jump points of $r \mapsto \tau_v(r)$), we have
		\begin{equation} \label{tauderiv}
			\frac{\D}{\D r} \tau_v(r) = \frac{r}{\langle \gamma_v, \dot \gamma_v \rangle} = \frac{1}{|\dot \gamma_v| \, \beta_v(r)}, \end{equation} %\frac{r}{\big\langle \gamma_v\big(\tau_v(r)\big), \dot \gamma_v\big(\tau_v(r)\big) \big\rangle}. \end{equation}
		where $\gamma_v$ and $\dot \gamma_v$ are evaluated at the exit time $\tau_v(r)$.
	\end{env_lem}
	\begin{proof}
		Let $\rho_v(t) = \sup_{s \le t} |\gamma_v(s)|$ denote the running maximum.  On the set of times where $\rho_v(t)$ is increasing, we have that $\rho_v(t) = |\gamma_v(t)|$.  For such a time $t$, we compute
			\begin{equation} \label{tauderiv_proof1}
				\tfrac{\D}{\D t} \rho_v(t)^2 = 2\rho_v(t) \cdot \tfrac{\D \rho_v}{\D t}(t) = 2 \langle \gamma_v(t), \dot \gamma_v(t) \rangle. \end{equation}
				
		The function $\tau_v(r)$ is the right-continuous inverse of $\rho_v(t)$, in the sense that $(\rho_v \circ \tau_v)(r) = r$ and $(\tau_v \circ \rho_v)(t) \ge t$.  By the chain rule, we have $\tfrac{\D \rho_v}{\D t} \tfrac{\D \tau_v}{\D r} = 1$.  Using the fact that $\rho_v(\tau_v(r)) = r$ and \eqref{tauderiv_proof1}, we have proved \eqref{tauderiv}. Since $\tau_v(r)$ is the exit time from $B(0,r)$, the running maximum increases at $\tau_v(r)$.  Clearly, $(\rho_v \circ \tau_v)(r) = r$.
	\end{proof}
	
	An upper bound on the exit angle $\alpha_v$ corresponds to a lower bound on $\beta_v$, since the arccosine function is decreasing.  Recall that the (lower) density of a set $A \subseteq \R$ is defined by $\density(A) := \liminf_{r \to \oo} \big|A \cap [0,r]\big|$, where the vertical bars denote Lebesgue measure on $\R$.  
	
	Define the random lens-shaped sets $L_v(r) = L_v(g,r) = B\big( \gamma\big(\tau_v(r)\big), 2 \big) \cap B(0,r)$.  We emphasize that these are len-shaped sets \emph{in the initial fixed coordinate chart}; by contrast, the lens-shaped set $D_{v,r} = B(0, 2) \cap B(o_{v,r}, r)$ is the image of $L_v(r)$ after the POV coordinate change.  For all $g \in \Omega_+$ and $v \in \V_g$, the set-valued function $r \mapsto L_v(r)$ is lower-semicontinuous.
	
	Trivially, $\gamma(\tau_v(r)) \in L_r$, so
		\begin{equation} \label{lowercontrol}
			\mbox{if $Z_{L_v(r)}(g) \le h$, then $|\dot \gamma_v\big( \tau_v(r) \big)| \le C$,} \end{equation}
	where $C = 1/\sqrt{1+h}$ is estimated using the minimum eigenvalue of the metric on the set $L_v(r)$.  
	
	Fix some $\epsilon > 0$.  Define the (random) sets of radii
		\begin{equation} \label{Qidef}
			Q_v^1 = Q_v^1(g) = \big\{ r : \beta_v(r) \ge \tfrac{1}{(1 + 2\epsilon) \mu |\dot \gamma_v|} \big\} \qquad \mathrm{and} \qquad Q_v^2 = Q_v^2(g,h) = \big\{ r : Z_{L_r}(g) \le h \big\}. \end{equation}			
	On the set $Q_v^1$, we have a lower bound on $\beta_v$, in terms of the (Euclidean) exit speed.  On $Q_v^2$, the upper bound on $Z_{L_r}$ gives a lower bound on the exit speed.  
	
	Lemma \ref{density_lemma1} states that the density of $Q_v^1$ is bounded below by $\tfrac{\epsilon}{1+2\epsilon}$.  Lemma \ref{density_lemma2} states that for sufficiently large $h$, the density of $Q_v^2$ is bounded below by $1 - \tfrac{\epsilon}{2}$, uniformly in $v \in \V_g$.  By considering the intersection along with the estimate \eqref{lowercontrol}, this gives a uniform lower bound on the density of $Q_v^1 \cap Q_v^2$.
	
	We now prove the Frontier Theorem using these two density estimates.  After the proof, we state and prove Lemmas \ref{density_lemma1} and \ref{density_lemma2}.
	
	\begin{proof}[Proof of Theorem \ref{frontier_thm}]
		Let $\epsilon \in (0,\tfrac 1 2)$.  By Lemma \ref{density_lemma1}, with probability one, $\density(Q_v^1) \ge \tfrac{\epsilon}{1+2\epsilon}$.  By Lemma \ref{density_lemma2}, we may choose $h$ sufficiently large so that, with probability one, $\density(Q_v^2) \ge 1 - \tfrac{\epsilon}{2}$.  By the inclusion-exclusion principle, we have $1 \ge \density(Q_v^1) + \density(Q_v^2) - \density(Q_v^1 \cap Q_v^2)$, hence
			\begin{equation}
				\density(Q_v^1 \cap Q_v^2) \ge \tfrac{\epsilon}{1+2\epsilon} + 1 - \tfrac{\epsilon}2 - 1 > 0 \end{equation}
		since $\epsilon \in (0,\tfrac 1 2)$. Define $Q_v := Q_v^1 \cap Q_v^2$.  Since geodesics are parametrized by constant Riemannian speed, $1 = \langle \dot \gamma, g \dot \gamma \rangle = |\dot \gamma|^2 \langle \tfrac{\dot\gamma}{|\dot\gamma|}, g \tfrac{\dot\gamma}{|\dot\gamma|} \rangle$, hence $|\dot \gamma|^2 \le \|g^{-1}\|$.  For $r \in Q_v$, then, we have that $|\dot \gamma( \tau_v(r) )| \le \|g^{-1}\|^{1/2}_{L_r} \le \sqrt h$, hence $\beta_v(r) \ge 1 / (1 + 2\epsilon) \mu |\dot \gamma_v| \ge 1 / (1+2\epsilon) \mu \sqrt h$.  Let $\theta = \arccos \tfrac{1}{(1+2\epsilon) \mu \sqrt h}$.  This completes the proof of Theorem \ref{frontier_thm}.  
	\end{proof}
	
	We now state and prove the first density lemma.
	
	\begin{env_lem}[First Density Lemma] \label{density_lemma1}
		With probability one, for all $v \in \V_g$, $\density(Q_v^1) \ge \tfrac{\epsilon}{1+2\epsilon}$, uniformly in the direction $v$.
	\end{env_lem}
	\begin{proof}
		Since $\tau_v(r)$ is right-continuous, we can use the fundamental theorem of calculus to write
			$$\tau_v(r) = \int_0^r \tfrac{1}{|\dot \gamma_v| \, \beta_v(r)} \sD r' + \jumps_v([0,r]),$$
		where $\jumps_v([0,r])$ denotes the total height the function $\tau_v(r)$ jumps on the interval $[0,r]$.\footnote{Formally, $\jumps_v([0,r]) = \int_0^r \lim_{h\to 0} \big( \tau_v(r'+h) - \tau_v(r') \big) \sD r'$.}
	
		Write $Q_v^1(r) := Q_v^1 \cap [0,r]$, and $\neg Q_v^1(r) := (Q_v^1)^c \cap [0,r]$.  We will prove the lower bound $|Q_v^1(r)| \ge \tfrac{\epsilon}{1+2\epsilon}$ for large $r$. Choose $r$ large enough so that $\tau_v(r) \le (1+\epsilon) \mu r$ by Lemma \eqref{lem_taur_estimate}.  Using this and the decomposition $[0,r] = Q_v^1(r) \cup \neg Q_v^1(r)$, we have
			\begin{eqnarray*} 
				(1+\epsilon) \mu r \ge \tau_v(r) &=& \int_{Q_v^1(r)} \tfrac{1}{|\dot \gamma_v| \, \beta_v(r')} \sD r' + \int_{\neg Q_v^1(r)} \tfrac{1}{|\dot \gamma_v| \, \beta_v(r')} \sD r' + \jumps_v([0,r]), \\
				&\ge& 0 + (1+2\epsilon) \mu \big| \neg Q_v^1(r) \big| + 0,
			\end{eqnarray*}
		where we trivially estimate the non-negative terms by zero; on the set $\neg Q_v^1(r)$, we use the lower bound $\tfrac{1}{\beta_v(r')} \ge (1+2\epsilon) \mu |\dot\gamma_v|$.  Using the fact that $|Q_v^1(r)| = r - |Q_v^1(r)^c|$ and rearranging the inequality $\tfrac{1+\epsilon}{1+2\epsilon}r \ge r - |Q_v^1(r)|$, we have proved the lemma
	\end{proof}
	
	Before stating the second density lemma, we introduce some discretization methods originally used in \cite{lagatta2009shape}.  These methods are based on first-passage percolation, which is a discrete model of stochastic geometry.  We will tessellate Euclidean space by unit cubes, and consider a dependent first-passage percolation model on the centers of these cubes.
	
	Following \cite{lagatta2009shape}, we define the $*$-lattice to be exactly the graph $\Z^d$, along with all its diagonal edges.  Formally, the vertex set is $\Z^d$, and two points are $*$-adjacent if $|z - z'|_{L^\oo} = 1$.  Note that if $z$ and $z'$ are $*$-adjacent, then the Euclidean distance between $z$ and $z'$ is at most $\sqrt d$.
	
	Let $X : \Z^d \to \R$ be some real-valued random field on the $*$-lattice.  We use the notation $X(\Gamma) := \sum_{z \in \Gamma} X_z$.
	
	\begin{env_thm}[Spatial Law of Large Numbers] \label{SpLLN}
		Let $\{X_z\}$ be a non-negative random field on the $*$-lattice which has a translation-invariant law and satisfies a finite-range dependence estimate.  Write $m = 3^d$, and let $X_1, \cdots, X_m$ be $m$ independent copies of the random variable $X_0$.  Suppose furthermore that $\EE \max\{X_1, \cdots, X_m\}^{2m+1} < \oo$.  Let $\xi = \EE X_z$ denote the mean of $X_z$.  
		
		For all $\epsilon > 0$, with probability one, there exists $N$ such that if $n \ge N$ and $\Gamma$ is a finite $*$-connected set containing the origin with $|\Gamma| \ge N$, then $(1-\epsilon) \xi |\Gamma| \le X(\Gamma) \le (1+\epsilon) \xi |\Gamma|.$
	\end{env_thm}
%	\begin{proof}
		In \cite{lagatta2009shape}, we proved this theorem as Lemmas 2.2 and 2.3 assuming a stronger exponential moment estimate.  By following more closely the argument of Cox and Durrett \cite{cox1981slt}, one can prove the theorem under a finite moment estimate.\footnote{The Cox-Durrett argument for lattice FPP on $\Z^d$ involves constructing $2d$ disjoint paths between the origin and a point $z$ on the $*$-lattice, then estimating the probability that the random passage time along one of these paths is particularly large. If the individual passage times were independent (as they are in usual lattice FPP), then the passage times along paths would be independent, since the paths are disjoint. The argument essentially holds in our model of dependent site FPP, with two modifications. First, since there are $3^d$ neighboring sites in the $*$-lattice (rather than just $2d$), we correspondingly need $3^d$ determined paths. Second, we cannot hope for the passage times of paths to be independent, owing to the finite range of dependence. Even still, for any point $z$ sufficiently far from the origin, it is possible to construct $3^d$ paths between $0$ and $z$ which are separated by more than the range of dependence, except near $0$ and $z$. This way, the passage times of paths are nearly independent, in that correlations decay rapidly as $|z| \to \oo$.}
%	\end{proof}
	
	For all $v \in S^{d-1}$, let $\zeta_v(t) \in \Z^d$ denote the nearest lattice point to the point $\gamma_v(t) \in \R^d$, breaking ties in some uniform way.  For all $v$, the function $t \mapsto \zeta_v(t)$ is a continuous-time process with nearest-neighbor jumps. Let $\tilde \gamma_v(r) = \bigcup_{s \le \tau_v(r)} \zeta_v(s) \subseteq \Z^d$ be the discretization of the curve $\gamma_v$; namely, all the lattice points which it is near.  If we represent $\tilde \gamma_v(r)$ by the union of boxes at the lattice points $z \in \tilde\gamma_v(r)$, then this is a covering of the curve. The next lemma states that the sizes of the sets $\tilde \gamma_v(r)$ are uniformly controlled for directions which yield minimizing geodesics.
	
	\begin{env_lem} \label{Ctildegamma}
		There exists $C \ge 1$ such that with probability one, there exists $r_0$ such that if $r \ge r_0$ and $v \in \V_g$, then $r \le |\tilde \gamma_v(r)| \le C r.$
	\end{env_lem}
	\begin{proof}
		The lower bound $|\tilde \gamma_v(r)| \ge r$ is trivial:  the curve $\gamma_v$ connects the origin to the sphere of radius $r$, so it must meet at least $r$ unit cubes.
	
		The upper bound relies on the Shape Theorem and the Spatial Law of Large Numbers.  Let $B_z = B^\oo(z, 1/2)$ denote the unit cube centered at $z$, and let $\varsigma_{v,z}$ denote the Euclidean length of $\gamma_v$ restricted to the unit cube $B_z$.  If $\varsigma_{v,z} < 1/4$, we say that the curve $\gamma_v$ \emph{barely meets} the cube $B_z$, and if $\varsigma_{v,z} \ge 1/4$, we say that $\gamma_v$ \emph{substantially meets} the cube $B_z$.  Let $\tilde \gamma_v'(r) = \{ z \in \Z^d : \mbox{$\varsigma_{v,z} \ge 1/4$} \}$ represent the unit cubes which $\gamma_v$ substantially meets.  The set $\tilde \gamma_v'(r)$ is $*$-connected; see the discussion following (2.8) of \cite{lagatta2009shape}.  Clearly, $0 \in \tilde \gamma_v'(r)$.
		
		Each time $\gamma_v$ substantially meets some cube $B_z$, it may barely meet up to $3^d-1 \le 3^d$ of its neighbors; this is a worst-case estimate.  This demonstrates that $\tilde \gamma_v'(r)$ is a subset of $\tilde \gamma_v(r)$ with density at least $1/3^d$. That is, $|\tilde \gamma_v'(r)| \ge \tfrac{1}{3^{d}} \, |\tilde \gamma_v(r)|.$
							
		Let $X_z = 1/\|g^{-1}\|_{B_z}$ denote the minimum eigenvalue of the metric $g$ on the unit cube $B_z \subseteq \R^d$.  Write $\gamma_{v}(r) := \gamma_v|_{[0,\tau_v(r)]} \subseteq \R^d$ for the geodesic segment on the time interval $[0, \tau_v(r)]$.  The geodesic segment $\gamma_{v}(r)$ is minimizing, so by the Shape Theorem, with probability one, there exists $r_1$ so that if $r \ge r_1$ and $v \in \V_g$, then $L_g[ \gamma_{v,r} ] \le (1+\epsilon)\mu r$.  Since $\tilde \gamma_v'(r)$ is a subset of $\tilde \gamma_v(r)$, we have
			\begin{equation} \label{Ctildegamma_proof2}
				(1+\epsilon) \mu r \ge L_g[ \gamma_{v,r} ] = \sum_{\tilde\gamma_v(r)} L_g[ \gamma_{v,r} \cap B_z ] \ge \sum_{\tilde\gamma_v'(r)} L_g[ \gamma_{v,r} \cap B_z ] \ge \tfrac 1 4 \sum_{\tilde\gamma_v'(r)} X_z, \end{equation}
		since if $\gamma_v$ substantially meets the cube $B_z$, then the Riemannian length of $\gamma_v$ restricted to that cube must be at least $\tfrac 1 4 X_z$.
		
		We now apply the Spatial Law of Large Numbers to the field $X_z$.  Write $\xi = \EE X_z$ for the mean of $X_z$, and note that by Theorem \ref{ZD_t}, $X_z$ satisfies the moment estimate.  Since the set $\tilde\gamma_v'(r)$ is $*$-connected and contains the origin, the Spatial LLN applies:  with probability one, there exists $r_2$ so that if $r \ge r_2$ and $v \in \V_g$, then $X( \tilde\gamma_v'(r) ) \ge (1-\epsilon) \xi |\tilde \gamma_v'(r)|$.  Combining this with \eqref{Ctildegamma_proof2}, we have
			\begin{equation}
				(1+\epsilon) \mu r \ge \tfrac 1 4 X( \tilde \gamma_v'(r) ) \ge \tfrac 1 4 (1-\epsilon) \xi |\tilde \gamma_v'(r)| \ge \tfrac 1 4  (1-\epsilon) \xi \cdot \tfrac{1}{3^d} |\tilde \gamma_v(r)|. \end{equation}
		Letting $C = 4 \cdot 3^d \tfrac{(1+\epsilon)\mu}{(1-\epsilon)\xi}$ completes the proof that $|\tilde \gamma_v(r)| \le C r$ for large $r$.
	\end{proof}

	We now use Lemma \ref{Ctildegamma} to prove the second density lemma.

	\begin{env_lem}[Second Density Lemma] \label{density_lemma2}
		Let $\epsilon > 0$.  There exists $h \ge 0$ such that, with probability one, for all $v \in \V_g$, $\density(Q_v^2) \ge 1-\tfrac{\epsilon}{2}$, uniformly in the direction $v$. 
	\end{env_lem}
	\begin{proof}
		Define the Euclidean ball $B_z = B(z,2 + \tfrac{1}{2} \sqrt d) \subseteq \R^d$ for each lattice point $z \in \Z^d$.  Let $X_z$ be the indicator function for the event $\{ Z_{B_z} > h \}$, and define $p(h) := \EE X_z = \PP( Z_{B_z} > h )$.  The random variable $Z_{B_z}$ is finite almost surely, so $p(h) \to 0$ as $h \to \oo$.  Let $C$ be as in Lemma \ref{Ctildegamma}, and choose a value of $h$ large enough so that $p(h) \le \epsilon/4C\sqrt d$.  By the Spatial Law of Large Numbers, with probability one, there exists $r_1$ such that if $r \ge r_1$ and $v \in \V_g$, then		
			\begin{equation} \label{density2_proof2}
				X(\tilde \gamma_v(r)) \le 2 p(h) |\tilde \gamma_v(r)| \le \tfrac{\epsilon}{2\sqrt d} r \end{equation}
		using the estimates $|\tilde \gamma_v(r)| \le Cr$ and $p(h) \le \epsilon/4C\sqrt d$.
		
		Let $\hat\zeta_v(r) := \zeta_v\big( \tau_v(r) \big)$ denote the lattice point nearest to the exit location $\gamma\big( \tau_v(r) \big)$.  The process $r \mapsto \hat\zeta_v(r)$ is a continuous ``$r$-time'' jump process on the lattice.  
		
		 Since the lens-shaped set $L_v(r)$ is a subset of the ball $B_{\hat\zeta_v(r)}$, we have that
			\begin{equation} \label{density2_proof3}
				\mbox{if $Z_{L_v(r)}(g) > h$, then $X_{\zeta_v'(r)} = 1$.} \end{equation}	
		
		Let $\ell$ denote Lebesgue measure on $\R$, and let $\mu_v = \ell \circ \hat\zeta_v^{-1}$ denote the push-forward of Lebesgue measure via the map $\hat\zeta_v : \R \to \Z^d$.  By simple plane geometry, the diameter of each set $\hat\zeta_v^{-1}(z) \subseteq \R$ is at most $\sqrt d$.  Consequently, with probability one, $\mu_v(z) \le \sqrt d$ for all $v \in \V_g$.  When $\mu_v(z) \approx \sqrt d$, it means that the geodesic $\gamma_v$ exits many balls near $z$.
		
		Let $\phi_v(r) = 1$ if $Z_{L_r}(g) > h$, and $0$ otherwise.  By \eqref{density2_proof3} and \eqref{density2_proof2}, we have
			\begin{equation} \label{density2_proof4}
				\big| \{ r : Z_{L_r}(g) > h \} \big| = \int_0^r \phi(r') \sD r' \le \int_0^r X_{\zeta_v'(r')} \sD r' = \sum_{\hat\zeta_v([0,r])} X_z \, \mu_v(z) \le \sqrt d \sum_{\tilde \gamma_v(r)} X_z \le \tfrac{\epsilon}{2} r. \end{equation}
		Since $Q_v^2$ is the complement of the set $\{r : Z_{L_r}(g) > h \}$, this completes the proof.
	\end{proof}

% V	
\section{Construction of the Bump Surface} \label{sect_bump}

	To prove Theorem \ref{bump_thm}, we work in a different coordinate system in order to construct bump metrics; finding ideal coordinate systems to work in is a common tactic in both differential geometry and physics.  Normal coordinates are familiar in elementary Riemannian geometry \cite{lee1997rmi}:  at any point $x$ on a Riemannian manifold $(M,g)$ we may change coordinates so that at $x$ the metric is locally flat, i.e., the metric $g_{ij}$ is just the Euclidean metric $\delta_{ij}$ with vanishing Christoffel symbols.  The curvature is an intrinsic geometric invariant, and does not take a canonical form in normal coordinates.
	
	Based on work of Fermi \cite{fermi1922atti}, Manasse and Misner \cite{manasse1963fermi} developed \emph{Fermi normal coordinates}, a coordinate system which is adapted to a particular geodesic.  In this coordinate system $(t,n)$, the geodesic curve traces the $t$-axis, along which the metric $g_{ij}$ takes the form of the Euclidean metric $\delta_{ij}$ and the Christoffel symbols vanish.  Furthermore, the coordinates are normal along the geodesic: to get to the point $(t,n)$ from the origin, we follow the geodesic $\gamma$ for time $t$, then move along a geodesic which is normal to $\gamma$ at time $t$ for time $n$.
	
			\begin{env_thm}[Existence of Fermi Normal Coordinates] \label{fermilemma}
				Let $(M,g)$ be a two-dimensional Riemannian manifold.  Fix a point $x \in M$, as well as a geodesic $\gamma$ starting at $x$.  Let $K(t)$ be the scalar curvature at the point $\gamma(t)$.  There exists an open neighborhood $U$ of the origin in $\R^2$ and a $C^2$-diffeomorphism (coordinate change) $\Phi_g: U  \to M$ such that
				\begin{itemize}
					\item The map $\Phi_g$ sends the $t$-axis in $U$ to the geodesic:  $\Phi_g(t,0) = \gamma(t)$.  It follows that, along the geodesic, the metric is locally flat and the Christoffel symbols vanish:  $g_{ij}(t,0) = \delta_{ij}$ and $\Gamma_{ij}^k(t,0) = 0$.
					\item If we define $\tilde g_{11}(t,n) = 1 - \tfrac 1 2 K(t) n^2$, $\tilde g_{12}(t,n) = 0,$ and $\tilde g_{22}(t,n) = 1$ in a neighborhood of the horizontal axis in $U$, then $(\Phi_* g)_{ab} = \tilde g_{ab} + O(n^3)$. %$\tilde g_{ij} - g_{ij} = O(n^3)$.
				\end{itemize}
			\end{env_thm}

	We outline some of the arguments behind this theorem in Appendix \ref{fermiproof}, following the work of Poisson \cite{poisson2004relativist}.\newline

			Now, we wish to define the bump metric $b = b(g)$ in a manner which depends continuously on the metric $g$ and its first and second derivatives at the origin only.  To formalize this notion, we introduce the equivalence relation $\sim$ on the space $\Omega_+$ of Riemannian metrics, defined by $g \sim g'$ if $\|g - g'\|_{C^2(0)} = 0$.\footnote{Thus $g \sim g'$ if $g_{ij}(0) = g'_{ij}(0)$, $g_{ij,k}(0) = g'_{ij,k}(0)$ and $g_{ij,kl}(0) = g'_{ij,kl}(0)$, for all indices $i,j,k,l$.}  Let $\Gamma_{ij}^k(g,x)$ and $K(g,x)$ denote the Christoffel symbols and scalar curvature of the metric $g$ at the point $x \in \R^2$, as defined by the formulae in equation \eqref{geoquantitiesdef}.  At the origin, these quantities are polynomials in the terms
				\begin{equation} \label{gterms0}
					\mbox{$g_{ij}(0)$, $g_{ij,k}(0)$, $g_{ij,kl}(0)$, $g^{ij}(0)$, and ${g^{ij}}_{,k}(0)$}. \end{equation}
			Thus, if $g \sim g'$ then $\Gamma_{ij}^k(g,0) = \Gamma_{ij}^k(g',0)$ and $K(g,0) = K(g',0)$.			
			
			Let $\Omega_0 = \Omega_+ / \!\sim$ denote the quotient space of $\Omega_+$ by the relation $\sim$, with quotient map $\pi_0 : \Omega_+ \to \Omega_0$.  For each $g \in \Omega_+$, we denote the equivalence class $\pi_0(g)$ by $[g]$.  Let $A = \{ g : Z_0(g) \le 2h \}$ as in \eqref{Adef}, and let $A_0 := \pi_0(A)$ be the image of $A$ under the quotient map $\pi_0$.

			\begin{env_lem}
				$A_0$ is a compact subset of the space $\Omega_0$.
			\end{env_lem}
			\begin{proof}
				Consider the finite-dimensional vector space $\R^{18}$ with the $L^{\oo}$ norm $\|v\| = \max_k \left\{ |v^k| \right\}$, and define a map $\Omega_0 \to \R^{18}$ by sending the equivalence class $[g]$ to the vector $\big(g_{11}(0), g_{12}(0), g_{22}(0), \dots, g_{22,22}(0) \big)$.  This map is an isometry with respect to the $\|\cdot\|_{C^{2}(0)}$ norm on $\Omega_0$, so $\Omega_0$ has the structure of an open cone within a finite-dimensional normed linear space.  To show that $A_0$ is a compact subset of $\Omega_0$, it suffices to show that that the seminorm $\|g\|_{C^2(0)}$ is bounded above and below on $A$:
					$$2h \ge \|g\|_{C^2} \ge \|g\|_{C^1} = \|g\|_{C^1} \, \frac {  \|g^{-1}\|_{C^1} }{ \|g^{-1}\|_{C^1} } \ge \frac { \|gg^{-1}\|_{C^1} }{ \|g^{-1}\|_{C^1} } = \frac{1}{ \|g^{-1} \|_{C^1} }\ge \frac{1}{2h}.$$
			\end{proof}
			
			The compactness of $A_0$ will feature prominently in our analysis.  We will parametrize the bump surface $b(g)$ continuously via the data of $g$ at the origin, i.e., by equivalence classes $[g] \in A_0$.  Since the set $A_0$ is compact, this will mean that quantities of interest are bounded uniformly in the metric $g$.  

			Let $\gamma_g := \gamma_{\E_1}(g,\cdot)$ be the geodesic in the metric $g$ starting at the origin in direction $\E_1$, and let $K(g,x)$ be the scalar curvature of $g$ at the point $x$.  We next introduce Fermi normal coordinates at the origin, adapted along the geodesic $\gamma_g$.  By Theorem \ref{fermilemma}, there exists a neighborhood $U$ of the origin and a map $\Phi_g : U \to \R^2$ (both depending on the metric $g$) such that the pull-back metric $\Phi^{-1}_* g$ takes the form 
				\begin{equation} \label{Phiinverseg}
					(\Phi^{-1}_* g)_{11}(t,n) = 1 - \tfrac 1 2 K(g, \gamma_g(t) ) n^2, \qquad (\Phi^{-1}_* g)_{12}(t,n) = 0, \qquad (\Phi^{-1}_* g)_{22}(t,n) = 1, \end{equation}
			up to $O(n^3)$ on $U$.  The map sends the horizontal axis to the geodesic:  $\Phi_g(t,0) = \gamma_g(t)$.  In particular, $\Phi_g(0) = 0$. Let $\Psi_g : \R^2 \to \R^2$ be the third-order Taylor polynomial of $\Phi_g$ at the origin; note that $\Psi_g(0) = 0$ and that $\Psi_g$ is defined on all of $\R^2$.
			
			%By the uniform bound \eqref{0estimate} on the metric $g \in A$, there is a minimal such neighborhood $U = \bigcap_{[g] \in A_0} U(g)$.			
			
%			The triangular domain $\II$ depends continuously on $\tau$.  We assume henceforth that $\tau$ is sufficiently small so that $\II$ is contained in the neighborhood $U$, hence the coordinate change $\Phi_g$ is defined on $\II$.  

			\begin{env_lem} \label{Psicont}
				The coefficients of the polynomial $\Psi_g$ are algebraic functions in the terms \eqref{gterms0}, hence are continuous functions of the equivalence class $[g]$.  %Consequently, there is a constant $M$ depending on $h$ so that the coefficients of $\Psi_g$ all belong to $[-M,M]$.
			\end{env_lem}
			\begin{proof}
				Write $\tilde g_{ab} = (\Phi^{-1}_* g)_{ab}$ for the pull-back metric defined by \eqref{Phiinverseg}.  In coordinates, the metrics $\tilde g_{ab}$ and $g_{ij}$ are related via the transformation $\Phi_g$ by the change-of-variable equation
					\begin{equation} \label{hatgphig}
						\tilde g_{ab}(\Phi_g(u)) = {\Phi_g^i}_{,a}(u) {\Phi_g^j}_{,b}(u) g_{ij}(u), \end{equation}
			where the subscripts after the commas denote partial derivatives of the components of the function $\Phi_g$.  Plugging in $u = 0$ and using the fact that $\tilde g_{ab}(0) = \delta_{ab}$, we see that the first-order terms ${\Psi_g^i}_{,a}(0) = {\Phi_g^i}_{,a}(0)$ solve a polynomial system of equations with coefficients \eqref{gterms0}, hence are algebraic functions of these terms.  
			
			The analysis of the second- and third-order terms is similar, since formula \eqref{Phiinverseg} implies that $\tilde g_{ab,c}(0) = 0$, $\tilde g_{11,22}(0) = -K_0(g)$, and $\tilde g_{ab,cd}(0) = 0$ for other values of $a$, $b$, $c$ and $d$.  We take the first derivative of \eqref{hatgphig} using the chain rule, plug in $u = 0$, and use the fact that $\tilde g_{ab,c}(0) = 0$ to see that ${\Phi_g^i}_{,ab}(0)$ is a algebraic function of the terms $g_{ij}(0)$ and $g_{ij,k}(0)$.  
			
			We take another derivative of \eqref{hatgphig} to analyze the third-order terms.  The second derivatives of $\tilde g_{ab}$ are not quite canonical, due to the presence of the scalar curvature $K_0(g)$.  Nonetheless, this is no obstruction, since $K_0(g)$ is a polynomial in the terms \eqref{gterms0}, hence $\tilde g_{ab}$ a polynomial in the terms \eqref{gterms0}.
			\end{proof}
			
			Next, we wish to define the number $\tau$, described in Theorem \ref{bump_thm}.  The constant $\tau$ represents a uniform length scale imposed on all the bump surfaces $b(g)$ near the origin.  
			
			As a consequence of Lemma \ref{Psicont}, both the functions $\Psi_g$ and $\Psi_g^2/\Psi_g^1$ are locally Lipschitz maps, with Lipschitz constants varying continuously in $[g] \in A_0$.  Let $L_1(g)$ be the Lipschitz constant for $\Psi_g$ on the Euclidean ball $B(0,\sqrt 2)$, and let $L_2(g)$ be the Lipschitz constant for $\Psi_g^2/\Psi_g^1$ on the Euclidean ball $B(0,\sqrt 2)$.   Let 
				\begin{equation} \label{Ldef}
					L = \sup_{[g] \in A_0} \left\{ 1, L_1(g), L_2(g) \right\} \end{equation}
			be the largest such Lipschitz constant on the set $B(0, \sqrt 2)$.

			Since $\Phi_g : U \to \R^2$ is a local $C^2$-diffeomorphism at the origin, there exists $\delta(g) > 0$ so that the polynomial $\Psi_g$ is a $C^2$-diffeomorphism on the closed Euclidean ball $B(0,\delta(g))$. We may choose this function $g \mapsto \delta(g)$ to vary continuously in $[g] \in A_0$, since the coefficients of $\Psi_g$ are continuous functions of $[g] \in A_0$ by the previous lemma.  Since $A_0$ is compact, there is a minimum such 
				\begin{equation} \label{deltadef}
					\delta := \inf_{[g] \in A_0} \delta(g) > 0. \end{equation}
			
			By assumption, the geodesic $\gamma_g$ satisfies $\gamma_g(0) = 0$ and $\dot \gamma_g(0) = \E_1$.  By the geodesic equation \eqref{geoeqn_geombg}, the second and third derivatives $\ddot \gamma_g(0)$ and $\dddot \gamma_g(0)$ of the geodesic at the origin are polynomial functions in $\Gamma_{ij}^k(0)$ and $\Gamma_{ij,l}^k(0)$, hence vary continuously in $[g] \in A_0$.  Define the constant
				\begin{equation} \label{gammagest}
					M = \sup_{[g] \in A_0} \max_k \left\{ 1, |\ddot \gamma^k_g(0)|, |\dddot \gamma^k_g(0)| \right\} < \oo. \end{equation}
			The constant $M$ lets us uniformly control the fluctuations of the plane curve $\gamma_g$ near the origin.  The assumption that $M \ge 1$ is by no means essential to the analysis, but it does make various calculations simpler. Let $\theta \in [0, \tfrac \pi 2)$ be the parameter assumed in Section \ref{sect_bump_short}, and define $\phi := \tfrac 1 2 \big( \tfrac \pi 2 - \theta \big)$ as in that section. Choose $\tau > 0$ to satisfy
				\begin{equation} \label{tau}
					\tau < \min \left\{\frac{\delta}{\sqrt{2}}, \frac{1}{2M}, \frac{\cos \phi}{L\sqrt 2 + 3M}, \frac{\tan \phi}{L \sqrt{2}  + 10M^2} \right\}. \end{equation}
%					\tau < \min \left\{\frac{\delta}{\sqrt{2}}, \frac{1}{2M}, \frac{1}{L\sqrt 2 + 3M}, \frac{\rho}{L \sqrt{2} + 10M^2} \right\}. \end{equation}
			Since $M \ge 1$ by assumption, it follows that $\tau \le \tfrac 1 2$.\newline
			
			Now that we have a natural length scale $\tau$, we are ready to define the curvature of the bump surface.  Recall that curvature is measured in units of $1/\operatorname{length}^2$.  Define $K_+ := 4\pi^2 / \tau^2$, which represents the positive curvature at the ``top'' of the bump. We are going to construct the bump surface so that geodesic transitions from the origin, where curvature equals $K_0(g)$, to a region of constant positive curvature $K_+$.  Even though the curvature at the origin is a random variable, it is uniformly bounded by 
				\begin{equation} \label{Kmaxdef}
					\Kmax = \max \left\{1, K_+, \sup_{[g] \in A_0} |K_0(g)| \right\}. \end{equation}
					
			In the Fermi coordinate chart, define the compact triangular region
				\begin{equation} \label{IIdef}
					\II = \left\{ (t,n) \in \R^2 ~:~ 0 \le t \le \tau \mathrm{~and~} |n| \le \frac{t}{\sqrt{\Kmax}} \right\} \end{equation}
			along the horizontal axis $(t,0)$. Note that the polynomial $\Psi_g$ is well-defined on $\II$ for all $g \in A$, and is identical for all metrics in the equivalence class $[g]$.  If $u \in \II$, then
				\begin{equation} \label{udelta}
					|u| \le \tau \sqrt{1 + \frac{1}{\Kmax}} \le \tau \sqrt{2} \le \delta, \end{equation}
			since $\Kmax \ge 1$ and $\tau \le \delta/\sqrt{2}$ by assumption.  This implies by the definition of the constant $\delta$ that the polynomial $\Psi_g$ is a $C^2$-diffeomorphism on the region $\II$.  Furthermore, since $\tau \le 1$, the region $\II$ is entirely contained in the Euclidean ball $B(0, \sqrt 2)$, so the polynomial $\Psi_g$ is Lipschitz on $\II$ with constant less than $L$.\newline

			We next define the curvature profile of the geodesic along the bump surface.  For each $[g] \in A_0$, define the piecewise-linear function $K^{(g)} : [0,\tau] \to \R$ by
				\begin{equation} \label{Kdef_curvature}
					K(t) := K^{(g)}(t) = \begin{cases} K_0(g) + (K_+ - K_0(g)) \tfrac{t}{\tau/4}, & 0 \le t \le \tfrac{\tau}{4} \\ K_+, & \tfrac{\tau}{4} \le t \le \tau. \end{cases} \end{equation}
			By the definition of the constant $\Kmax$, it is readily apparent that
				\begin{equation} \label{KtKmax}
					\sup_{0 \le t \le \tau} |K(t)| \le \Kmax. \end{equation}

			We now consider $\II$ as a closed coordinate chart, and define a ``bump surface'' metric $f_{ab}(g)$ in Fermi normal coordinates on $\II$.  Fermi coordinates are canonical up to the choice of curvature profile along the horizontal geodesic, which we take to be the function $K(t)$.  Define the symmetric $2$-tensor $f_{ab}$ by
				\begin{equation} \label{bumpdef}
					f_{11}(t,n) = 1 - \tfrac 1 2 K(t) n^2, \qquad f_{12}(t,n) = 0, \qquad f_{22}(t,n) = 1. \end{equation}
			We easily verify that $f(u)$ is positive-definite, hence a Riemannian metric:
				$$\inf_{u \in \II} f_{11}(u) \ge \inf_{t \le \tau} f_{11} \left(t, \tfrac{t}{\sqrt{\Kmax}} \right) = \inf_{t \le \tau} \left(1 - \frac 1 2 K(t) \frac{t^2}{\Kmax} \right) \ge 1 - \frac 1 2 \Kmax \frac{\tau^2}{\Kmax} \ge \frac 1 2 > 0,$$
			by the estimates $K(t) \le \Kmax$ and $\tau \le 1$.  Thus for every $[g] \in A_0$, $f$ is a Riemannian metric in Fermi normal coordinates on the coordinate chart $\II$, and its curvature profile along the $t$-axis is the function $K(t)$
			
			Define $\JJ_g := \Psi_g(\II) \subseteq \R^2$ to be the image of $\II$ under the diffeomorphism $\Psi_g$.  The dependence on $g$ in this definition arises from the components of the polynomial $\Psi_g$.  Since the coefficients of $\Psi_g$ are continuous in $g$, the function $g \mapsto \JJ_g(g)$ is continuous in the Hausdorff topology on closed sets in $\R^2$.  Clearly, $\JJ_g$ is a simply-connected compact set with piecewise-smooth boundary.
			
			\begin{env_lem} \label{JJCC}
				For all $[g] \in A_0$, the compact set $\JJ_g$ contains the origin, and is a subset of the frontier cone $FC$ defined in \eqref{FCdef}.  The set $\JJ_g$ is in the interior of $B(0,1)$.
			\end{env_lem}
			\begin{proof}
				The origin is contained in the set $\II$, and mapped to itself under $\Psi_g$.  Thus $0 \in \JJ_g$ for all $[g] \in A_0$.
				
%				To show that $\JJ_g(g)$ is contained in the set $\C$ for all $g \in A$, we first use the 
				
				Since the Fermi coordinate change $\Phi_g$ sends the horizontal axis to the geodesic $\gamma_g$, the polynomial $\Psi_g$ sends the horizontal axis to the third-order Taylor approximation to $\gamma_g$, defined by
					\begin{equation} \label{gammabdef}
						\gamma_b(t) := \Psi_g(t,0) = \E_1 t + \tfrac 1 2 \ddot\gamma_g(0) t^2 + \tfrac 1 6 \dddot\gamma_g(0)t^3 \end{equation}
				for $t \in [0,\tau]$.  This is a vector-valued polynomial in $t$, and its coefficients are uniformly bounded by the constant $M$ defined by \eqref{gammagest}.

				The curve $\gamma_b$ remains in the right half-plane:  if $t > 0$, then
					$$\gamma_b^1(t) \ge t - Mt^2 - Mt^3 \ge t - 2Mt^2 \ge t(1 - 2M\tau) > 0,$$
				since $t \le \tau < 1/2M < 1$ by assumption.
				
				To prove that $\JJ_g = \Psi_g(\II)$ is a subset of the frontier cone $FC$, it suffices to show that $\Psi_g^1(u) \le \cos \phi$ and $| \Psi_g^2(u) / \Psi_g^1(u) | \le \tan \phi$ for all $u \in \II$.  By construction, the maps $\Psi_g^1$ and $\Psi_g^2/\Psi_g^1$ are Lipschitz on $\II$ with Lipschitz constant less than $L$.  We use this fact, along with some simple estimates on the curve $\gamma_b$.
		
				For any $u \in \II$ and $t \in [0,\tau]$,
					$$\Psi_g^1(u) \le |\Psi_g^1(u) - \Psi_g^1(t,0)| + \gamma_b^1(t) \le L \diam \II + (\tau + M \tau^2 + M \tau^3) \le L\sqrt 2 \tau + 3M\tau < \cos \phi,$$
				since $\diam \II < \sqrt 2 \tau$ by \eqref{udelta}, $1 \le M$, and $\tau^3 \le \tau^2 \le \tau \le \tfrac{\cos \phi}{L\sqrt 2 + 3M}$.
				
				Similarly, the function $\Psi_g^2/\Psi_g^1$ has Lipschitz constant $L$, so
					\begin{eqnarray*}
						\left| \frac{\Psi_g^2(u)}{\Psi_g^1(u)} \right| &\le& \left| \frac{\Psi_g^2(u)}{\Psi_g^1(u)} - \frac{\Psi_g^2(t,0)}{\Psi_g^1(t,0)} \right| + \left| \frac{\gamma_b^2(t)}{\gamma_b^1(t)} \right| \le L \diam \II + \frac{0 + M\tau^2 + M\tau^3}{1 - M\tau^2 - M\tau^3} \\
						&\le& L \sqrt{2} \tau + \frac{2M\tau}{1 - 2M\tau} \le L \sqrt{2} \tau + 2M\tau (1 + 4M\tau) \le L \sqrt{2} \tau + 10M^2 \tau < \tan \phi,
					\end{eqnarray*}
				since $x/(1-x) \le x(1 + 2x)$, $1 \le M$, and $\tau^2 \le \tau < \tfrac{\tan \phi}{L\sqrt 2 + 10M^2}$.  This completes the proof that $\JJ_g$ is a subset of the frontier cone $FC$.  In fact, we have shown that $\JJ_g - \{0\}$ is in the interior of $FC$, hence $\JJ_g$ is in the interior of $B(0,1)$.
			\end{proof}
			
			Consider $\JJ_g$ as a closed manifold with piecewise-smooth boundary, and let $\Psi_*$ be the map which pushes forward a metric in Fermi coordinates from $\II$ to a metric on $\JJ_g$.  In the next lemma, we define the bump metric $b(g)$ on all of $\R^2$.  On the set $\JJ_g$, the metric $b(g)$ agrees with $\Psi_* f$, the push-forward of metric $f$ defined in the Fermi coordinate system, defined in \eqref{bumpdef}.  Away from the unit ball $B(0,1)$, the bump metric is equal to $\delta$, the Euclidean metric.  The content of the next lemma is that we can $C^2$-smoothly interpolate between the two metrics in a manner which varies continuously in the parameter $g$. 
			
			\begin{env_lem} \label{bumplemma_definition}
				There exists a continuous map $b : A \to \Omega_+$ such that for all $g \in A$, $b(g)$ is a $C^2$-smooth Riemannian metric on $\R^2$ satisfying
					\begin{equation}
						\mbox{$b(g)(x) = (\Psi_* f)(x)$ for $x \in \JJ_g$, \qquad and \qquad $b(g)(x) = \delta$ for $x \notin B(0,1)$.} \end{equation}
				The function $b : A \to \Omega_+$ is $\F_0$-measurable.
			\end{env_lem}
			\begin{proof}
				By construction, the metric $f$ is $C^{2,1}$-smooth, and satisfies the uniform bound $\|f\|_{C^{2,1}(\II)} \le C_1 := \Kmax^{3/2} / \tau$.  Since the map $\Psi_g$ is a polynomial with coefficients varying continuously in $g$, the operator $(\Psi_g)_* : C^{2,1}(\II, \SPD) \to C^{2,1}(\JJ_g, \SPD)$ has operator norm bounded by some constant $C_2$, independently of $g \in A$. Define $f_g := (\Psi_g)_* f$ for the image of the bump metric on $\JJ_g$. Multiplying the operator bounds, we have that the bump metric (in Euclidean coordinates) satisfies $\|f_g\|_{C^{2,1}(\JJ_g)} \le C_3 := C_1 C_2$. The construction is $\F_0$-measurable since the bump metric only depends on $g$ using the information at the origin.
				
				In this proof, we first extend the metric $f_g$ to a metric on the unit disk $B(0,1)$, then use a mollifying function which forces the metric to equal the Euclidean metric along the boundary. Applying Lemma \ref{mo_lem} with a modified covariance structure, we can define an operator $m_{\JJ_g} : C^2(\JJ_g, \Sym) \to C^2(\R^2, \Sym)$ which extends a tensor field on $\JJ_g$ to one defined on all of $\R^2$, and which also satisfies the property that the tensor field $m_g := m_{\JJ_g}(f_g)$ is positive-definite on the ball $B(0,1)$ for every $g \in A$.
				
				Define the set $\overline{\JJ} := \overline{\bigcup_{g \in A} \JJ_g}$, the smallest compact set which contains each set $\JJ_g$. The set $\overline{\JJ}$ is a compact subset of the interior of $B(0,1)$. Let $\rho \in (0,1)$ denote the minimal radius so that $\overline{\JJ} \subseteq B(0,\rho)$. Let $\varphi : [0,\oo) \to [0,1]$ be a monotone smooth function which equals $1$ on the interval $[0,\rho]$, and $0$ on the interval $[1,\oo)$. We define the bump metric $b(g)(x) = \big(1 - \varphi(|x|) \big) m_g(x) + \varphi(|x|) \delta$. 
				
				We now show that the function $g \mapsto b(g)$ is continuous. Suppose that $g^n \to g$ is a convergent sequence in $\Omega_+$, meaning that $\|g^n - g\|_{C^2(K)} \to 0$ for each compact $K \subseteq \R^2$. Let $f^n$ and $f$ be the corresponding bump metrics defined in the set $\II$, so that $f^n \to f$ in the Banach space $C^2(\II)$. The sets $\JJ_{g^n}$ converge to $\JJ_g$ in the Hausdorff topology, and $(\Psi_{g^n})_* f^n \to (\Psi_g)_* f$ on compact sets.\footnote{More precisely, if $K' \subseteq \JJ_{g^n}$ eventually, then $(\Psi_{g^n})_* f^n \to (\Psi_g)_* f$ in $C^2(K')$.} Consequently, $m_{g^n} (\Psi_{g^n})_* f^n$ converges to $m_g (\Psi_{g})_* f$ on compact sets. Since the function $\varphi$ doesn't depend on $g$, we have that $b(g^n) \to b(g)$. Thus $g \mapsto b(g)$ is continuous.
				
				The $\F_0$-measurability of $g \mapsto b(g)$ is preserved since the function $\Psi_g$ depends only on $g$ via the metric information at the origin, hence so does the set $\JJ_g$ and the operator $m_{\JJ_g}$. This completes the proof.

			\end{proof}
			
			Now that we have constructed the bump metric $b(g)$, we are ready to prove that it satisfies the geometric properties stated in Theorem \ref{bump_thm}:  $b(g)$ is equal to $g$ up to second derivatives at the origin; the central geodesic $\gamma_b$ on the bump metric is not minimizing on the time interval $[0,\tau]$; and most crucially, if $g$ is sufficiently close to its bump metric $b(g)$, then also $\gamma_g$ too is not minimizing on the time interval $[0,\tau]$.

			\begin{env_lem} \label{bumplemma_origin}
				For all $g \in A$, the bump metric $b(g)$ agrees with $g$ up to second derivatives at the origin in $\JJ_g$: $\|g - b\|_{C^2(0)} = 0.$ This includes the fact that their respective scalar curvatures $K_0(g)$ and $K_0(b)$ at the origin are equal.
			\end{env_lem}
			\begin{proof}
				Let $\tilde g = \Phi^{-1}_* g$ denote the metric $g$, changed into Fermi normal coordinates.  Since these coordinates take the canonical form \eqref{gexpansion}, they are determined up to the scalar curvature $K_0(g)$ at the origin.  By our construction of the bump metric, $f = \Psi^{-1}_* b$ also has scalar curvature $K_0(g)$ at the origin.  Consequently, the metrics $\tilde g$ and $f$ are equal at the origin.  The map $\Phi_g \Psi_g^{-1} : \JJ_g \to \R^2$ is equal to the identity up to second derivatives at the origin, so the metrics $g$ and $b$ are also equal up to second derivatives at the origin.	
			\end{proof}

			To show that geodesics on the bump surface are not minimizing, we will make use of the method of Jacobi fields; for a good overview, see Chapter 10 of \cite{lee1997rmi}.
			
			For any metric $g$, pick a tangent vector $n = n(g)$ which is orthogonal to $\E_1$ at the origin (i.e. $\langle n, g(0) \E_1 \rangle = 0$), and let $\dot \gamma_g^\perp(t)$ be the parallel translation of $n$ along the geodesic $\gamma_g$.  Note that for all $t$, the vector field $\dot \gamma_g^\perp$ is normal to $\dot \gamma_g$ with respect to $g$.  As before, let $K(g,x)$ denote the scalar curvature of $g$ at a point $x \in \R^2$.  Let $j(g,t)$ be a solution to the \emph{Jacobi equation}
				\begin{equation} \label{jacobieqn_general}
					j''(g,t) + K(g,\gamma_g(t)) \, j(g,t) = 0, \end{equation}
			and define the \emph{Jacobi field}
				$$J(g,t) = j(g,t) \, \dot \gamma_g^\perp(t)$$
			along the geodesic $\gamma_g(t)$.  The Jacobi field $J$ measures the second-order variations of the geodesic $\gamma_g$.  
			
			If $j(g,t_1) = 0$ and $j(g, t_2) = 0$ for two different times $t_1$ and $t_2$, then the points $\gamma_g(t_1)$ and $\gamma_g(t_2)$ are called \emph{conjugate points} along the geodesic $\gamma_g$.  A consequence is that the geodesic $\gamma_g$ is not minimizing beyond the time interval $[t_1, t_2]$; this is Jacobi's Theorem (cf. Theorem 10.15 of \cite{lee1997rmi}).
			
			Let $b \in b(A)$ denote any bump metric, and consider the unit-speed geodesic $\gamma_b$ starting at the origin in direction $\E_1$ (the explicit form of the curve $\gamma_b$ is given by \eqref{gammabdef}).  By our construction of the bump metric, the scalar curvature along the geodesic $\gamma_b$ is constant and equal to $K_+ = \tfrac{4\pi^2}{\tau^2}$ on the time interval $[\tfrac{\tau}{4}, \tau]$.  In this case, we can solve the Jacobi equation \eqref{jacobieqn_general} explicitly.  
			
			Let $j(b,t)$ be the solution to the equation
					$$j''(b,t) + \tfrac{4\pi^2}{\tau^2} j(b,t) = 0$$
			subject to the initial conditions $j(b,\tfrac{\tau} 4) = 0$ and $j'(b,\tfrac{\tau} 4) = \tfrac{2\pi}{\tau}$.  This has the explicit solution
					\begin{equation} \label{jacobieqn_b}
						j(b,t) = \sin\!\big( \tfrac{2\pi}{\tau} (t - \tfrac{\tau}{4})\big) \end{equation}
			on the interval $t \in [\tfrac{\tau} 4, \tau]$, so that $j(b,\tfrac{3\tau} 4) = 0$.  Thus the points $\gamma_b(\tfrac{\tau}{4})$ and $\gamma_b(\tfrac{3\tau}{4})$ are conjugate along $\gamma_b$, so Jacobi's Theorem implies that $\gamma_b$ is not minimizing.  We record this as the following lemma:
				
			\begin{env_lem} \label{gammab_notmin}
				For any bump surface $b \in b(A)$, the geodesic $\gamma_b$ is not minimizing between times $0$ and $\tau$.
			\end{env_lem}

			As a consequence of the explicit solution \eqref{jacobieqn_b} for $j(b,t)$, we have that
				\begin{equation} \label{jbtau}
					j(b,\tau) = -1. \end{equation} 
			Let $j(g,t)$ be the solution to the equation
					\begin{equation} \label{jacobieqn_g}
						j''(g,t) + K(g, \gamma_g(t)) j(g,t) = 0 \end{equation}
			subject to the initial conditions $j(g,\tfrac{\tau} 4) = 0$ and $j'(g,\tfrac{\tau} 4) = \tfrac{2\pi}{\tau}$.  We will show that if $g$ is sufficiently close to its bump metric $b(g)$, then $j(g,\tau)$ will be close to $j(b,\tau) = -1$.  This implies that $j(g,t)$ changes sign on the interval $[0,\tau)$, hence there is some point $\gamma_g(t)$ conjugate to $\gamma_g(\tfrac{\tau}{4})$.  By Jacobi's Theorem, this implies that $\gamma_g$ is not minimizing.
			
			\begin{env_lem} \label{gammag_notmin}
				There exists a constant $\epsilon > 0$ so that if $\|g - b(g)\|_{C^{2,1}(FC)} < \epsilon$, then the geodesic $\gamma_g$ is not minimizing between times $0$ and $\tau$.
			\end{env_lem}
			\begin{proof}
				By the estimates \eqref{Lipest_K}, \eqref{Lipest_gamma} and \eqref{Lipest_Gamma}, there exist constants $\epsilon_1(b)$, $C_1(b)$ and $L(b)$ (varying continuously in the bump metric $b$) such that if $\|g - b\|_{C^{2,1}(FC)} < \epsilon_1$, then 
					\begin{eqnarray}
						\big| K(g, \gamma_g(t) ) - K(b, \gamma_b(t) ) \big| &\le& L \, \|g - b\|_{C^{2,1}(FC)} \cdot |\gamma_g(t) - \gamma_b(t)| \label{gammag_proof0} \\
						&\le& L \, \|g - b\|_{C^{2,1}(FC)} \cdot C_1 \, \big\| \Gamma(g,\cdot) - \Gamma(b,\cdot) \big\|_{C^{0,1}(FC)} \nonumber \\
						&\le& L^2 C_1 \, \|g - b\|_{C^{2,1}(FC)}^2. \label{gammag_proof1}
					\end{eqnarray}
				
				The Jacobi equation \eqref{jacobieqn_g} is a second-order ODE, featuring the coefficient $K(g,\gamma_g(t))$.  The function $(g,t) \mapsto \gamma_g(t)$ is locally Lipschitz; this and \eqref{gammag_proof0} implies that $(g,t) \mapsto K(g, \gamma_g(t))$ is locally Lipschitz.  Consequently, a theorem of smoothness of solutions of ODEs and \eqref{gammag_proof1} imply that
					\begin{eqnarray} 
						\sup_{t \in [0, \tau]} \big| j(g,t) - j(b,t) \big| &\le& C_2 \sup_{t \in [0,\tau]} \big| K(g, \gamma_g(t) ) - K(b, \gamma_b(t) ) \big| \nonumber \\
						 &\le& C_2 L^2 C_1 \, \|g - b\|_{C^{2,1}(FC)}^2 \label{gammag_proof2}
					\end{eqnarray}
				for some constant $C_2(b)$ varying continuously in $b$.
				
				Since the constants $C_1$, $C_2$ and $L$ vary continuously in $b$, we may define $C_3 = \sup \{C_2 L^2 C_1\} < \oo$, where the supremum is taken over the compact set of bump metrics $b(A)$.  Similarly, define $\epsilon = \inf\{ \epsilon_1, \tfrac{1}{\sqrt{2 C_3}} \} > 0$.  
				
				If $\|g - b\| < \epsilon$, then \eqref{gammag_proof2} implies that
					$$j(g,\tau) \le -1 + \big| j(g,\tau) - j(b,\tau) \big| \le -1 + C_3 \epsilon^2 \le -1 + \tfrac{1}{2} < 0.$$
				The Jacobi field changes sign on the interval $(0,\tau)$, hence there are conjugate points, so Jacobi's Theorem implies that $\gamma_g$ is not minimizing.
			\end{proof}

			This completes the proof of Theorem \ref{bump_thm}.

%\bibliographystyle{alpha}
%\bibliography{biblio}

\include{geodesics_supplement_mar13_content}

\setcounter{page}{56}
\end{document}